\input amstex
\let\myfrac=\frac
\input eplain
\let\frac=\myfrac
\input epsf




\loadeufm \loadmsam \loadmsbm
\message{symbol names}\UseAMSsymbols\message{,}

\font\myfontdefault=cmr10

\font\mytdmchapfont=cmb10 at 14pt
\font\mytdmheadfont=cmb10 at 10pt
\font\mytdmsubheadfont=cmr10

\magnification 1200
\newif\ifinappendices
\newif\ifundefinedreferences
\newif\ifchangedreferences
\newif\ifloadreferences
\newif\ifmakebiblio
\newif\ifmaketdm

\undefinedreferencesfalse
\changedreferencesfalse


\loadreferencestrue
\makebibliofalse
\maketdmfalse

\def\headpenalty{-400}     
\def\proclaimpenalty{-200} 

%
%

\def\alphanum#1{\ifcase #1 _\or A\or B\or C\or D\or E\or F\or G\or H\or I\or J\or K\or L\or M\or N\or O\or P\or Q\or R\or S\or T\or U\or V\or W\or X\or Y\or Z\fi}
\def\gobbleeight#1#2#3#4#5#6#7#8{}

\newwrite\references
\newwrite\tdm
\newwrite\biblio

\newcount\chapno
\newcount\headno
\newcount\subheadno
\newcount\procno
\newcount\figno
\newcount\citationno

\def\setcatcodes{%
\catcode`\!=0 \catcode`\\=11}%

\ifloadreferences
    {\catcode`\@=11 \catcode`\_=11%
    \global\def\_@citation@And{1}
\global\def\_@citation@BallGromSch{2}
\global\def\_@citation@Frib{3}
\global\def\_@citation@GallKapMard{4}
\global\def\_@citation@LabA{5}
\global\def\_@citation@LeRoux{6}
\global\def\_@citation@SmiB{7}
\global\def\_@citation@SmiE{8}
\global\def\_@citation@SmiF{9}
\global\def\_@proc@GallKapMarden{1.2}
\global\def\_@proc@PresentationChIVExistence{1.4}
\global\def\_@proc@PresentationChIVConvexRealisation{1.5}
\global\def\_@proc@NastyAlgebraI{2.2}
\global\def\_@proc@NastyAlgebraII{2.3}
\global\def\_@proc@NastyAlgebraIV{2.4}
\global\def\_@proc@NastyAlgebraV{2.5}
\global\def\_@proc@EmptyCylinder{2.6}
\global\def\_@subhead@ConeControl{2.3}
\global\def\_@proc@LemmaConeControl{2.7}
\global\def\_@proc@LemmaWeakConvexity{2.8}
\global\def\_@proc@LemmaSchottkyI{2.9}
\global\def\_@proc@SchottkyGroups{2.10}
\global\def\_@proc@SchottkyTrouserDecomposition{3.1}
\global\def\_@proc@DecompositionStepI{3.2}
\global\def\_@proc@ObtainingANonElemHandle{3.3}
\global\def\_@fig@FigureDehnTwistI{3.1}
\global\def\_@proc@MakingOtherHandlesHyperbolic{3.4}
\global\def\_@proc@CorMakingOtherHandlesHyperbolic{3.5}
\global\def\_@proc@DecompositionStepII{3.6}
\global\def\_@proc@MakingTrousersSchottky{3.7}
\global\def\_@proc@CorMakingTrousersSchottky{3.8}
\global\def\_@proc@DecompositionStepIII{3.9}
\global\def\_@proc@TheFinalCut{3.10}
\global\def\_@proc@PrincipalObstructionToJoin{4.1}
\global\def\_@proc@QuotientOfInvDomainExists{4.2}
\global\def\_@proc@OrderOfFixedPoints{4.3}
\global\def\_@proc@ThmFriberg{4.4}
\global\def\_@proc@FundamentalGroupOfHomeo{4.5}
\global\def\_@subhead@HomologicalClasses{4.4}
\global\def\_@proc@ConditionForEquivalenceOfLifts{4.6}
\global\def\_@proc@LemmaBraguette{4.7}
\global\def\_@fig@FigureDistributionOfFixedPoints{4.2}
\global\def\_@proc@InverseZipper{4.8}
\global\def\_@proc@ExistenceForTrousers{5.1}
\global\def\_@proc@JoiningDistinctTrousers{5.2}
\global\def\_@proc@JoiningSameTrouser{5.3}
\global\def\_@proc@ExistenceOfConvexSets{5.4}
\global\def\_@proc@SolutionToSimpleConvexPlateauProblem{5.5}
\global\def\_@head@MakingHomeomorphism{A}
\global\def\_@proc@HomeoOfFixedPointSet{A.1}
\global\def\_@proc@HomeoEquivalence{A.2}
\global\def\_@proc@LimitsOfEndPtsOfGeos{A.3}
\global\def\_@proc@AdvancedPingPong{A.5}
\global\def\_@proc@FormulaForInteriorOfCircles{A.6}
\global\def\_@proc@CirclesAreNested{A.7}
\global\def\_@proc@LengthOfEvalTendsToInfinity{A.8}
\global\def\_@proc@CtyOfFPSPartI{A.9}
\global\def\_@proc@CtyOfFPSPartII{A.10}
    }%
\else
    \openout\references=references.tex
\fi

\newcount\newchapflag 
\newcount\showpagenumflag 

\global\chapno = -1 
\global\citationno=0
\global\headno = 0
\global\subheadno = 0
\global\procno = 0
\global\figno = 0

\def\resetcounters{%
\global\headno = 0%
\global\subheadno = 0%
\global\procno = 0%
\global\figno = 0%
}

\global\newchapflag=0 
\global\showpagenumflag=0 

\def\chinfo{\ifinappendices\alphanum\chapno\else\the\chapno\fi}%
\def\headinfo{\ifinappendices\alphanum\headno\else\the\headno\fi}%
\def\subheadinfo{\headinfo.\the\subheadno}
\def\procinfo{\headinfo.\the\procno}
\def\figinfo{\headinfo.\the\figno}
\def\citationinfo{\the\citationno}%
\def\nextheadno{\global\advance\headno by 1 \global\subheadno = 0 \global\procno = 0}
\def\nextsubheadno{\global\advance\subheadno by 1}
\def\nextprocno{\global\advance\procno by 1 \procinfo}
\def\nextfigno{\global\advance\figno by 1 \figinfo}

{\global\let\noe=\noexpand%
%
%
\catcode`\@=11%
\catcode`\_=11%
\setcatcodes%
!global!def!_@@internal@@makeref#1{%
!global!expandafter!def!csname #1ref!endcsname##1{%
!csname _@#1@##1!endcsname%
!expandafter!ifx!csname _@#1@##1!endcsname!relax%
    !write16{#1 ##1 not defined - run saving references}%
    !undefinedreferencestrue%
!fi}}%
!global!def!_@@internal@@makelabel#1{%
!global!expandafter!def!csname #1label!endcsname##1{%
!edef!temptoken{!csname #1info!endcsname}%
!ifloadreferences%
    !expandafter!ifx!csname _@#1@##1!endcsname!relax%
        !write16{#1 ##1 not hitherto defined - rerun saving references}%
        !changedreferencestrue%
    !else%
        !expandafter!ifx!csname _@#1@##1!endcsname!temptoken%
        !else
            !write16{#1 ##1 reference has changed - rerun saving references}%
            !changedreferencestrue%
        !fi%
    !fi%
!else%
    !expandafter!edef!csname _@#1@##1!endcsname{!temptoken}%
    !edef!textoutput{!write!references{\global\def\_@#1@##1{!temptoken}}}%
    !textoutput%
!fi}}%
!global!def!makecounter#1{!_@@internal@@makelabel{#1}!_@@internal@@makeref{#1}}%
!unsetcatcodes%
}
\makecounter{ch}%
\makecounter{head}%
\makecounter{subhead}%
\makecounter{proc}%
\makecounter{fig}%
\makecounter{citation}%
\def\newref#1#2{%
\def\temptext{#2}%
\edef\bibliotextoutput{\expandafter\gobbleeight\meaning\temptext}%
\global\advance\citationno by 1\citationlabel{#1}%
\ifmakebiblio%
    \edef\fileoutput{\write\biblio{\noindent\hbox to 0pt{\hss$[\the\citationno]$}\hskip 0.2em\bibliotextoutput\medskip}}%
    \fileoutput%
\fi}%
\def\cite#1{%
$[\citationref{#1}]$%
\ifmakebiblio%
    \edef\fileoutput{\write\biblio{#1}}%
    \fileoutput%
\fi%
}%
%
%
%

\let\mypar=\par


\def\raggedleft{\leftskip=0pt plus 1fil \parfillskip=0pt}


\font\lettrinefont=cmr10 at 28pt
\def\lettrine #1[#2][#3]#4%
{\hangafter -#1 \hangindent #2
\noindent\hskip -#2 \vtop to 0pt{
\kern #3 \hbox to #2 {\lettrinefont #4\hss}\vss}}

\font\mylettrinefont=cmr10 at 28pt
\def\mylettrine #1[#2][#3][#4]#5%
{\hangafter -#1 \hangindent #2
\noindent\hskip -#2 \vtop to 0pt{
\kern #3 \hbox to #2 {\mylettrinefont #5\hss}\vss}}


\edef\Pagetitle={Blank}

\headline={\hfil\Pagetitle\hfil}

\footline={\hfil\myfontdefault\folio\hfil}

\def\nextoddpage
{
\newpage%
\ifodd\pageno%
\else%
    \global\showpagenumflag = 0%
    \null%
    \vfil%
    \eject%
    \global\showpagenumflag = 1%
\fi%
}


\def\newchap#1#2%
{%
%
%
\global\advance\chapno by 1%
\resetcounters%
%
%
\newpage%
\ifodd\pageno%
\else%
    \global\showpagenumflag = 0%
    \null%
    \vfil%
    \eject%
    \global\showpagenumflag = 1%
\fi%
\global\newchapflag = 1%
\global\showpagenumflag = 1%
%
%
{\font\chapfontA=cmsl10 at 30pt%
\font\chapfontB=cmsl10 at 25pt%
\null\vskip 5cm%
{\chapfontA\raggedleft\hfil%
{%
\ifnum\chapno=0
    \phantom{%
    \ifinappendices%
        Annexe \alphanum\chapno%
    \else%
        \the\chapno%
    \fi}%
\else%
    \ifinappendices%
        Annexe \alphanum\chapno%
    \else%
        \the\chapno%
    \fi%
\fi%
}%
\par}%
\vskip 2cm%
{\chapfontB\raggedleft%
\lineskiplimit=0pt%
\lineskip=0.8ex%
\hfil #1\par}%
\vskip 2cm%
}%
\edef\Pagetitle{#2}%
%
%
\ifmaketdm%
    \def\temp{#2}%
    \def\tempbis{\nobreak}%
    \edef\chaptitle{\expandafter\gobbleeight\meaning\temp}%
    \edef\mynobreak{\expandafter\gobbleeight\meaning\tempbis}%
    \edef\textoutput{\write\tdm{\bigskip{\noexpand\mytdmchapfont\noindent\chinfo\ - \chaptitle\hfill\noexpand\folio}\par\mynobreak}}%
\fi%
\textoutput%
}


\def\newhead#1%
{%
\ifhmode%
    \mypar%
\fi%
\ifnum\headno=0%
\ifinappendices
    \nobreak\vskip -\lastskip%
    \nobreak\vskip .5cm%
\fi
\else%
    \nobreak\vskip -\lastskip%
    \nobreak\vskip .5cm%
\fi%
\nextheadno%
\ifmaketdm%
    \def\temp{#1}%
    \edef\sectiontitle{\expandafter\gobbleeight\meaning\temp}%
    \edef\textoutput{\write\tdm{\noindent{\noexpand\mytdmheadfont\quad\headinfo\ - \sectiontitle\hfill\noexpand\folio}\par}}%
    \textoutput%
\fi%
\font\headfontA=cmbx10 at 14pt%
{\headfontA\noindent\headinfo\ -\ #1.\hfil}%
\nobreak\vskip .5cm%
}%


\def\newsubhead#1%
{%
\ifhmode%
    \mypar%
\fi%
\ifnum\subheadno=0%
\else%
    \penalty\headpenalty\vskip .4cm%
\fi%
\nextsubheadno%
\ifmaketdm%
    \def\temp{#1}%
    \edef\subsectiontitle{\expandafter\gobbleeight\meaning\temp}%
    \edef\textoutput{\write\tdm{\noindent{\noexpand\mytdmsubheadfont\quad\quad\subheadinfo\ - \subsectiontitle\hfill\noexpand\folio}\par}}%
    \textoutput%
\fi%
\font\subheadfontA=cmsl10 at 12pt
{\subheadfontA\noindent\subheadinfo\ #1.\hfil}%
\nobreak\vskip .25cm %
}%

%
%


\font\mathromanten=cmr10
\font\mathromanseven=cmr7
\font\mathromanfive=cmr5
\newfam\mathromanfam
\textfont\mathromanfam=\mathromanten
\scriptfont\mathromanfam=\mathromanseven
\scriptscriptfont\mathromanfam=\mathromanfive
\def\mathroman{\fam\mathromanfam}


\font\sansseriften=cmss10
\font\sansserifseven=cmss7
\font\sansseriffive=cmss5
\newfam\sansseriffam
\textfont\sansseriffam=\sansseriften
\scriptfont\sansseriffam=\sansserifseven
\scriptscriptfont\sansseriffam=\sansseriffive
\def\mathsf{\fam\sansseriffam}


\font\boldten=cmb10
\font\boldseven=cmb7
\font\boldfive=cmb5
\newfam\mathboldfam
\textfont\mathboldfam=\boldten
\scriptfont\mathboldfam=\boldseven
\scriptscriptfont\mathboldfam=\boldfive
\def\mathbf{\fam\mathboldfam}


\font\mycmmiten=cmmi10
\font\mycmmiseven=cmmi7
\font\mycmmifive=cmmi5
\newfam\mycmmifam
\textfont\mycmmifam=\mycmmiten
\scriptfont\mycmmifam=\mycmmiseven
\scriptscriptfont\mycmmifam=\mycmmifive

\def\hexa#1{\ifcase #1 0\or 1\or 2\or 3\or 4\or 5\or 6\or 7\or 8\or 9\or A\or B\or C\or D\or E\or F\fi}
\mathchardef\mathi="7\hexa\mycmmifam7B
\mathchardef\mathj="7\hexa\mycmmifam7C


\font\mymsbmten=msbm10 at 8pt
\font\mymsbmseven=msbm7 at 5.6pt
\font\mymsbmfive=msbm5 at 4pt
\newfam\mymsbmfam
\textfont\mymsbmfam=\mymsbmten
\scriptfont\mymsbmfam=\mymsbmseven
\scriptscriptfont\mymsbmfam=\mymsbmfive

\mathchardef\mybeth="7\hexa\mymsbmfam69
\mathchardef\mygimmel="7\hexa\mymsbmfam6A
\mathchardef\mydaleth="7\hexa\mymsbmfam6B


\def\placelabel[#1][#2]#3{{%
\setbox10=\hbox{\raise #2cm \hbox{\hskip #1cm #3}}%
\ht10=0pt%
\dp10=0pt%
\wd10=0pt%
\box10}}%

\def\placefigure#1#2#3{%
\medskip%
\midinsert%
\vbox{\line{\hfil#2\epsfbox{#3}#1\hfil}%
\vskip 0.3cm%
\line{\noindent\hfil\sl Figure \nextfigno\hfil}}%
\medskip%
\endinsert%
}


\newif\ifinproclaim%
\global\inproclaimfalse%
\def\proclaim#1{%
\medskip%
%
%
\bgroup%
\inproclaimtrue%
\setbox10=\vbox\bgroup\leftskip=0.8em\noindent{\bf #1}\sl%
}

\def\endproclaim{%
\egroup%
\setbox11=\vtop{\noindent\vrule height \ht10 depth \dp10 width 0.1em}%
\wd11=0pt%
\setbox12=\hbox{\copy11\kern 0.3em\copy11\kern 0.3em}%
\wd12=0pt%
\setbox13=\hbox{\noindent\box12\box10}%
\noindent\unhbox13%
\egroup%
\medskip\ignorespaces%
}

\def\proclaim#1{%
\medskip%
\bgroup%
\inproclaimtrue%
\noindent{\bf #1}%
\nobreak\medskip%
\sl%
}

\def\endproclaim{%
\mypar\egroup\penalty\proclaimpenalty\medskip\ignorespaces%
}

\def\noskipproclaim#1{%
\medskip%
\bgroup%
\inproclaimtrue%
\noindent{\bf #1}\nobreak\sl%
}

\def\endnoskipproclaim{%
\mypar\egroup\penalty\proclaimpenalty\medskip\ignorespaces%
}


\def\ninn{{n\in\Bbb{N}}}

\def\proof{{\noindent\bf Proof:\ }}

\def\msup{\mathop{{\mathroman Sup}}}
\def\minf{\mathop{{\mathroman Inf}}}

\def\psl{\Bbb{P}SL}
\def\qed{~$\square$}
\def\munion{\mathop{\cup}}
\def\minter{\mathop{\cap}}
\def\myitem#1{%
    \noindent\hbox to .5cm{\hfill#1\hss}
}

\catcode`\@=11
\def\Eqalign#1{\null\,\vcenter{\openup\jot\m@th\ialign{%
\strut\hfil$\displaystyle{##}$&$\displaystyle{{}##}$\hfil%
&&\quad\strut\hfil$\displaystyle{##}$&$\displaystyle{{}##}$%
\hfil\crcr #1\crcr}}\,}
\catcode`\@=12

\def\makeop#1{%
\global\expandafter\def\csname op#1\endcsname{{\mathroman #1}}}%

\def\makeopsmall#1{%
\global\expandafter\def\csname op#1\endcsname{{\mathroman{\lowercase{#1}}}}}%

\makeopsmall{ArcTan}%
\makeopsmall{ArcCos}%
\makeop{Arg}%
\makeop{Det}%
\makeop{Log}%
\makeop{Re}%
\makeop{Im}%
\makeop{Dim}%
\makeopsmall{Tan}%
\makeop{Ker}%
\makeopsmall{Cos}%
\makeopsmall{Sin}%
\makeop{Exp}%
\makeopsmall{Tanh}%
\makeop{Tr}%
\makeop{End}%
\makeop{Long}%
\makeop{Ch}%
\makeop{Exp}%
\makeop{Eval}%
\makeop{Lift}%
\makeop{Int}%
\makeop{Ext}%
\makeop{Aire}%
\makeop{Im}%
\makeop{Conf}%
\makeop{Exp}%
\makeop{Mod}%
\makeop{Log}%
\makeop{Ext}%
\makeop{Int}%
\makeop{Dist}%
\makeop{Aut}%
\makeop{Id}%
\makeop{GL}%
\makeop{SO}%
\makeop{Homeo}%
\makeop{Vol}%
\makeop{Ric}%
\makeop{Hess}%
\makeop{Euc}%
\makeop{Isom}%
\makeop{Max}%
\makeop{Long}%
\makeop{Fix}%
\makeop{Wind}%
\makeop{Mush}%
\makeop{Ad}%
\makeop{loc}%
\makeop{Len}%
\makeop{Area}%
\def\opTildeHomeo{{\mathroman H\widetilde{omeo}}}%
\font\mycirclefont=cmsy7
\def\textcircle{{\raise 0.3ex \hbox{\mycirclefont\char'015}}}
\def\mathcircle{\mathop{\textcircle}}

\let\emph=\bf

\hyphenation{quasi-con-formal}

%
%

\ifmakebiblio%
    \openout\biblio=biblio.tex %
    {%
        \edef\fileoutput{\write\biblio{\bgroup\leftskip=2em}}%
        \fileoutput
    }%
\fi%

\newref{And}{}
\newref{BallGromSch}{}
\newref{Frib}{}
\newref{GallKapMard}{}
\newref{LabA}{}
\newref{LeRoux}{}
\newref{SmiB}{Smith G., Hyperbolic Plateau problems, Preprint, Orsay (2005)}%
\newref{SmiE}{Smith G., Th\`ese de doctorat, Paris (2004)}%
\newref{SmiF}{Smith G., A Homomorphism for Tress Groups in the Sphere, in preparation}%

\ifmakebiblio%
    {\edef\fileoutput{\write\biblio{\egroup}}%
    \fileoutput}%
\fi%

%
%
%
\document
\myfontdefault
\global\chapno=1
\global\showpagenumflag=1
\def\Pagetitle{}
\null
\vfill
\def\centre{\rightskip=0pt plus 1fil \leftskip=0pt plus 1fil \spaceskip=.3333em \xspaceskip=.5em \parfillskip=0em \parindent=0em}%
\def\textmonth#1{\ifcase#1\or January\or February\or March\or April\or May\or June\or July\or August\or September\or October\or November\or December\fi}
\font\abstracttitlefont=cmr10 at 14pt
{\abstracttitlefont\centre Equivariant Plateau Problems\par}
\bigskip
{\centre Graham Smith\par}
\bigskip
{\centre \the\day\ \textmonth\month\ \the\year\par}
\bigskip
{\centre Max Planck Institute for Mathematics in the Sciences,\par
Inselstrasse 22.,\par
D-04103 Leipzig,\par
GERMANY\par}
\bigskip
\noindent{\emph Abstract:\ }Let $(M,Q)$ be a compact, three dimensional manifold of strictly negative
sectional curvature. Let $(\Sigma,P)$ be a compact, orientable surface of hyperbolic type (i.e. of genus at
least two). Let $\theta:\pi_1(\Sigma,P)\rightarrow\pi_1(M,Q)$ be a homomorphism. Generalising a recent 
result of Gallo, Kapovich and Marden concerning necessary and sufficient conditions for the
existence of complex projective structures with specified holonomy to manifolds of non-constant negative
curvature, we obtain necessary conditions on $\theta$ for the existence of a so called $\theta$-equivariant
Plateau problem over $\Sigma$, which is equivalent to the existence of a strictly convex immersion 
$i:\Sigma\rightarrow M$ which realises $\theta$ (i.e. such that $\theta=i_*$).
\bigskip
\noindent{\emph Key Words:\ }Kleinian groups, Fuchsian groups, Plateau problem, complex projective structures, immersions
\bigskip
\noindent{\emph AMS Subject Classification:\ }57M50 (30F10, 30F40, 32G15)
\vfill
\nextoddpage
\def\Pagetitle{\sl Equivariant Plateau Problems}
\newhead{Introduction}
\noindent In this paper, we study means of obtaining constant curvature realisations of homomorphisms
of fundamental groups of surfaces into fundamental groups of compact, three dimensional manifolds of strictly
negative sectional curvature.
\medskip
\noindent Let $(M,Q)$ be a pointed, compact three dimensional Riemannian manifold of strictly negative
sectional curvature. Let $(\Sigma,P)$ be a closed (orientable) surface of hyperbolic type (i.e. of genus
at least $2$). Let $(\tilde{M},\tilde{Q})$ be the universal cover of $(M,Q)$. Let 
$\theta:\pi_1(\Sigma,P)\rightarrow\pi_1(M,Q)$ be a representation of $\pi_1(\Sigma,P)$ in $\pi_1(M,Q)$.
We observe that it is in many ways preferable to work in the framework of pointed manifolds. Firstly, 
$\pi_1(M,Q)$ is canonically defined in terms of $(M,Q)$, and, secondly, the action of $\pi_1(M,Q)$ over
$(\tilde{M},\tilde{Q})$ is well defined.
\medskip
\noindent $\tilde{M}$ is a Hadamard manifold. We denote by $\partial_\infty\tilde{M}$ its ideal boundary,
and we make the following definition:
\proclaim{Definition \nextprocno}
\noindent Let $(\tilde{\Sigma},\tilde{P})$ be the universal cover of $(\Sigma,P)$. Let
$\theta:\pi_1(\Sigma,P)\rightarrow\pi_1(M,Q)$ be a homomorphism. A {\emph $\theta$-equivariant Plateau
problem\/} over $\Sigma$ is a function $\varphi:\tilde{\Sigma}\rightarrow\partial_\infty\tilde{M}$ such
that, for all $\gamma\in\pi_1(\Sigma,P)$:
$$
\varphi\circ\gamma = \theta(\gamma)\circ\varphi.
$$
\endproclaim
\noindent Such structures are introduced and discussed in detail by Labourie in \cite{LabA} and by the author in \cite{SmiB} and \cite{SmiE}.
\medskip
\noindent The group $\pi_1(M,Q)$ acts canonically on $\tilde{M}\munion\partial_\infty\tilde{M}$. 
We underline that $\pi_1(M,Q)$ acts on $\tilde{M}\munion\partial_\infty\tilde{M}$ {\emph from the right\/}. Thus, throughout this entire paper, composition is to be read {\emph from left to right\/}. A subgroup
$\Gamma$ of $\pi_1(M,Q)$ is said to be {\emph non-elementary\/} if and only if it has no fixed points in
$\tilde{M}\munion\partial_\infty\tilde{M}$, and we will say that it is {\emph elementary\/} otherwise.
In our case, the only elementary subgroups of $\pi_1(M,Q)$ are the trivial group, and isomorphic copies
of $\Bbb{Z}$.
\medskip
\noindent Let $\opHomeo_0(\partial_\infty\tilde{M})$ denote the connected component of the group of 
homeomorphisms of $\partial_\infty\tilde{M}$ which preserves the identity. The group $\pi_1(M,Q)$ acts faithfully over $\partial_\infty\tilde{M}$, and may thus be considered as a subgroup of 
$\opHomeo_0(\partial_\infty\tilde{M})$. Let $\opTildeHomeo_0(\partial_\infty\tilde{M})$ be the universal cover of $\opHomeo_0(\partial_\infty\tilde{M})$, and let
$\pi:\opTildeHomeo_0(\partial_\infty\tilde{M})\rightarrow\opHomeo_0(\partial_\infty\tilde{M})$ be the canonical projection. We have the following short exact sequence:
$$
0\rightarrow\Bbb{Z}_2\mathop{\rightarrow}^{i}\opTildeHomeo_0(\partial_\infty\tilde{M})
\mathop{\rightarrow}^{\pi}\opHomeo_0(\partial_\infty\tilde{M})\rightarrow 0.
$$
\noindent For $\Gamma$ an arbitrary group, and for 
$\varphi:\Gamma\rightarrow\opHomeo_0(\partial_\infty\tilde{M})$ a homomorphism, we define a {\emph lifting\/}
of $\varphi$ in $\opTildeHomeo_0(\partial_\infty\tilde{M})$ to be a homomorhism
$\hat{\varphi}:\Gamma\rightarrow\opTildeHomeo_0(\partial_\infty\tilde{M})$ such that 
$\pi\circ\hat{\varphi} = \varphi$. We say that $\varphi$ {\emph lifts\/} if such a lifting may be found. 
\medskip
\noindent Returning to the special case where $M$ is of constant sectional curvature equal to $-1$, we
recall an existence result, obtained by Gallo, Kapovich and Marden in \cite{GallKapMard}, which, when 
translated into our framework, may be expressed as follows:
\proclaim{Theorem \nextprocno\ [Gallo, Kapovich, Marden, 2000]}
\noindent Suppose that $M$ is of constant sectional curvature equal to $-1$. Then, there exists a 
$\theta$-equivariant Plateau problem $\varphi$ over $\Sigma$ if and only if $\theta$ is non-elementary 
and lifts.
\endproclaim
\proclabel{GallKapMarden}
\noindent The requirement that $\theta$ be non-elementary arises from the fact that $\varphi$ defines
a $\psl(2,\Bbb{C})$ structure over the surface $\Sigma$. If $\theta$ were elementary, then it would either be a subset of the rotation group, in which case it would define an $SO(3)$ structure over $\Sigma$, or
it would be a subset of the affine group, in which case it would define an affine structure over $\Sigma$.
Neither of these are possible, since they would induce non-negative curvature metrics over $\Sigma$, which
are excluded by the Gauss-Bonnet Theorem. This results yields the following interesting corollary:
\proclaim{Corollary \nextprocno}
\noindent Suppose that $M$ has constant sectional curvature equal to $-1$. There exists a locally convex, immersion $i:\Sigma\rightarrow M$ such that $\theta = i_*$ if and only if $\theta$ is non-elementary and lifts.
\endproclaim
\noindent The objective of this paper is to obtain an analogue of this result in the more general case
where $M$ is only of strictly negative sectional curvature. The main result of this paper is the following:
\proclaim{Theorem \nextprocno}
\noindent Suppose that $(M,Q)$ is a pointed, compact manifold of strictly negative sectional curvature. 
Let $(\Sigma,P)$ be a pointed, compact surface of hyperbolic type (i.e. of genus at least two). Let
$\theta:\pi_1(\Sigma,P)\rightarrow\pi_1(M,Q)$ be a homomorphism. Suppose that $\theta$ is non-elementary and
may be lifted to a homomorphism $\hat{\theta}$ of $\pi_1(\Sigma,P)$ into the group 
$\opTildeHomeo_0(\partial_\infty\tilde{M})$. Then there exists an equivariant Plateau problem for $\theta$.
\endproclaim
\proclabel{PresentationChIVExistence}
\noindent This result permits us to obtain a locally convex realisation of $\theta$:
\proclaim{Theorem \nextprocno}
\noindent If $\theta$ is non-elementary and lifts, then there exists a convex immersion 
$i:\Sigma\rightarrow M$ such that $\theta = i_*$.
\endproclaim
\proclabel{PresentationChIVConvexRealisation}
\noindent Theorem \procref{PresentationChIVExistence} is, in many respects, a generalisation to non-constant
curvature of Theorem \procref{GallKapMarden}, and will be proven using essentially the same techniques,
following the observation that many of the algebraic results used in \cite{GallKapMard}, which thus depend
on the structure of $\psl(2,\Bbb{C})=\opIsom(\Bbb{H}^3)$, may be expressed in a purely topological manner,
and may thus be applied to $\opIsom(\tilde{M})$. 
\medskip
\noindent Our approach differs from that of Gallo, Kapovich and Marden in various respects. By far the most significant, however, is our treatment of the algebraic obstruction. In \cite{GallKapMard}, the complex structure allows the algebraic obstruction to be rapidly treated using classical techniques associated to index theory. Since the case that we study is purely topological, such powerful tools are no longer available, and a much deeper understanding of the algebraic properties of the construction described in \cite{GallKapMard} is required. This forms the content of section $4$ and Appendix $A$. The more perspicacious reader will also observe that this deeper understanding allows us to simplify the trouser decomposition of \cite{GallKapMard}, in that it is no longer necessary that it take the form of a tree.
\medskip
\noindent The remaining differences between our approaches involve the sort of simplifications that one would expect to arise from any in-depth study of such a significant work. First, in section $3$, we adopt the perspective that the homomorphism, $\theta$, as opposed to the paths that we study, is changed by the Dehn twists. These approaches are trivially logically equivalent, but we believe that our perspective facilitates the reading of this paper by hiding the complexity of the construction inside the algebra. Moreover, we construct the Dehn twists used in \cite{GallKapMard} as compositions of much simpler Dehn twists. We believe that this also facilitates the understanding of our proofs.
\medskip
\noindent Finally, in section $2$, since we only work with mappings of hyperbolic type, we may bypass most of the technically complicated algebraic reasoning employed in \cite{GallKapMard} to obtain the algebraic results necessary for the construction of Schottky groups. The simpler conditions of the case we study allow us to prove these results in terms of elementary arguments concerning the fixed point sets of hyperbolic mappings and a general property of Hadamard manifolds with sectional curvature bounded above by $k<0$. This alternative approach is crucial for our results, since the algebraic approach used in \cite{GallKapMard} is only valid for the M\"obius group. In a way, this part of our work may also be considered as an illustrative example, facilitating the understanding of Gallo, Kapovich and Marden's result without requiring an in-depth understanding of the technical complexity of their own paper.
\medskip
\noindent In the second section of this paper, we obtain algebraic results concerning the construction of Schottky groups. In the third section, we use Dehn twists to construct a trouser decomposition of $\Sigma$ in such a manner that the $\theta$-image of the fundamental group of each trouser is a Schottky group. In the fourth section, we describe invariant domains of Schottky groups and study the algebraic properties of these domains. In the fifth section, we show how invariant domains may be used to construct $\theta$-equivariant Plateau problems over trousers, and we then show how to glue these functions together
to obtain a proof of Theorems \procref{PresentationChIVExistence} and \procref{PresentationChIVConvexRealisation}. Finally, in the appendix, we prove the homeomorphism equivalence of Schottky groups. In particular, this permits us to construct invariant domains for any Schottky group.
\medskip
\noindent I would like to thank Fran\c{c}ois Labourie for having drawn my attention to this problem, Kevin
Costello for help with algebraic topology, and Fr\'ed\'eric Leroux, for having introduced me to the results of Friberg. I would also like to
thank the Max Planck Institute for Mathematics in the Sciences in Leipzig for providing the comfort required
to write the english version of this paper.
\newhead{Schottky Groups}
\newsubhead{Introduction}
\noindent Throughout this section, $M$ will denote a three dimensional Hadamard manifold. The construction
that we will carry out in the sequel makes considerable use of Schottky groups. We define such groups as
follows:
\proclaim{Definition \nextprocno} 
\noindent Let $C_a^-$, $C_a^+$, $C_b^-$ and $C_b^+$ be disjoint Jordan curves in the sphere $S_2$ oriented
such that each one of these curves is situated in the exterior of the three others. We denote the group of homeomorphisms of $S_2$ which preserve orientation by $\opHomeo_0(S_2)$. We will adopt the convention that
this group acts on $S_2$ {\emph from the right\/}. A {\emph Schottky Group\/} is a subgroup $\Gamma$ of 
$\opHomeo_0(S_2)$ generated by two elements $a$ and $b$ such that:
$$
\opExt(C_a^-)\cdot a = \opInt(C_a^+),\qquad \opExt(C_b^-)\cdot b=\opInt(C_b^+).
$$
\noindent We will refer to the curves $(C_a^\pm,C_b^\pm)$ as {\emph generating circles\/} for the Schottky
group $\Gamma$ with respect to the pair of generators $(a,b)$.
\endproclaim
\noindent Let $\alpha:M\rightarrow M$ be an isometry of $M$. The mapping $\alpha$ may be extended to a 
homeomorphism of $M\munion\partial_\infty M$. A mapping is {\emph hyperbolic\/} if and only if it has
precisely two fixed points in $\partial_\infty M$ and no fixed points in $M$ (see \cite{BallGromSch}). In the sequel, we will work with actions of fundamental groups of compact manifolds of strictly negative sectional curvature on the universal covers of these manifolds. In such subgroups 
$\Gamma\subseteq \opIsom(M)$, every element that is different to the identity is necessarily hyperbolic.
\medskip
\noindent In this section, we show how Schottky subgroups of $\opHomeo_0(\partial_\infty M)$ may be constructed. We obtain the following results:
\proclaim{Proposition \nextprocno}
\noindent Let $M$ be a Hadamard manifold. Let $\Gamma$ be a subgroup of the group of isometries of $M$ 
consisting only of hyperbolic elements (and the identity). Let $\alpha$,$\beta$,$\xi$ and $\eta$ be
elements of $\Gamma$ such that the subgroup $\langle\alpha,\beta\rangle$ of $\Gamma$ is non-elementary.
There exists $K\in\Bbb{N}$ and a function
$N:\left\{\left|n\right|\geqslant K\right\}\rightarrow\Bbb{N}$ such that, for all $\left|k\right|\geqslant K$
and for all $\left|n\right|\geqslant N(k)$:
\medskip
\noindent\myitem{(1)} $\beta\alpha^k$ is hyperbolic,
\medskip
\noindent\myitem{(2)} $\delta_k=\eta\beta\alpha^k$ is hyperbolic,
\medskip
\noindent\myitem{(3)} $\delta_k^n\alpha$ is hyperbolic and has no fixed point in common with $\beta\alpha^k$, and
\medskip
\noindent\myitem{(4)} $\delta_k^n\xi$ is hyperbolic.
\endproclaim
\proclabel{NastyAlgebraI}
\proclaim{Proposition \nextprocno}
\noindent Let $M$ be a Hadamard manifold. Let $\Gamma$ be a subgroup of the group of isometries of $M$
consisting only of hyperbolic elements (and the identity). Let $\alpha$,$\beta$,$\gamma$ and $\eta$ be 
elements of $\Gamma$ such that the subgroup $\langle\alpha,\beta\rangle$ generated by $\alpha$ and
$\beta$ is non-elementary. Suppose moreover that $\xi$ and $\eta$ are both hyperbolic.
\medskip
\noindent For $k\in\Bbb{Z}$, we define $\delta_k=\eta\beta\alpha^k$. After exchanging $\xi$ and $\eta$ if
necessary, there exists $K\in\Bbb{N}$ and a function 
$N:\left\{\left|n\right|\geqslant K\right\}\rightarrow\Bbb{N}$ such that, for all 
$\left|k\right|\geqslant K$ and for all $\left|n\right|\geqslant N(k)$, 
$\langle\delta_k^n\xi\delta_k^{-n},\eta\rangle$ is a Schottky group.
\endproclaim
\proclabel{NastyAlgebraII}
\proclaim{Proposition \nextprocno}
\noindent Let $M$ be a Hadamard manifold. Let $\Gamma$ be a subgroup of the group of isometries of $M$
consisting only of hyperbolic elements (and the identity). Let $\alpha$,$\beta$ be elements of $\Gamma$
such that the subgroup $\langle\alpha,\beta\rangle$ is non-elementary. 
\medskip
\noindent For all $\sigma\in\Gamma$, there exists $K\in\Bbb{N}$ such that, for all 
$\left|k\right|\geqslant K$, $J_k = \sigma^k\alpha\sigma^{-k}\beta$ is hyperbolic.
\endproclaim
\proclabel{NastyAlgebraIV}

\proclaim{Proposition \nextprocno}
\noindent Let $M$ be a Hadamard manifold. Let $\Gamma$ be a subgroup of the group of isometries of $M$
consisting only of hyperbolic elements (and the identity). Let $\alpha$,$\beta$,$\xi$ and $\eta$ be elements
of $\Gamma$ such that:
\medskip
\myitem{(i)} the subgroup $\langle\alpha,\beta\rangle$ is non-elementary,
\medskip
\myitem{(ii)} $\alpha\beta\alpha^{-1}\beta^{-1}=\eta\xi$,
\medskip
\myitem{(iii)} the subgroup $\langle\xi,\eta\rangle$ is a Schottky group,
\medskip
\myitem{(iv)} $\xi\beta$ is not equal to the identity.
\medskip
\noindent Then there exists $k$ and $n$ such that, if we define $\gamma_k$ and $\delta_k$ by:
$$
\gamma_k = \beta\alpha^k,\qquad \delta_k = \xi\gamma_k,
$$
\noindent then, $\langle \gamma_k, (\delta_k^n\alpha)^{-1}\gamma_k^{-1}(\delta_k^n\alpha)\rangle$
and $\langle \xi, \delta_k^{-n}\eta\delta_k^n\rangle$ are Schottky groups:
\endproclaim
\proclabel{NastyAlgebraV}
\noindent In the second part of this section, we review properties of isometries of Hadamard manifolds.
In the third part, we study properties of conical open sets about geodesics. In the fourth part, we proof
general results about the construction of Schottky groups, and in the final part, we prove the above propositions.
\newsubhead{Hyperbolic Isometries of Hadamard Manifolds}
\noindent It is well known that a Hadamard manifold is canonically diffeomorphic (under the action of the 
exponential mapping) to its tangent space over any given point. Similarly, its ideal boundary is
canonically homeomorphic to the unit sphere in the tangent space over any given point. We have 
the following analogous result for normal bundles over geodesics:
\proclaim{Lemma \nextprocno}
\noindent Let $M$ be a Hadamard manifold of sectional curvature bounded above by $-\epsilon<0$. 
Let $\gamma$ be a geodesic in $M$. Let $N_\gamma$ be the normal bundle over $\gamma$. In otherwords:
$$
N_\gamma = \left\{ X \in T_{\gamma(t)}M | t\in\Bbb{R}, X_p\perp\partial_t\gamma(t) \right\}.
$$
\noindent The restriction of the exponential mapping to $N_\gamma$ is a diffeomorphism between $N_\gamma$ and $M$.
\medskip
\noindent Likewise, let $\Sigma_\gamma$ be the normal sphere bundle over $\gamma$ contained in $N_\gamma$.
In otherwords:
$$
\Sigma_\gamma = \left\{ X \in N_\gamma | \|X\| = 1 \right\}.
$$
\noindent Let us define $\partial_\infty\gamma$ by:
$$
\partial_\infty\gamma = \left\{\gamma(-\infty),\gamma(\infty)\right\}.
$$
\noindent The restriction of the Gauss-Minkowski mapping defines a homeomorphism between $\Sigma_\gamma$ and 
$\partial_\infty M\setminus\partial_\infty\gamma$.
\endproclaim
\proclabel{EmptyCylinder}
\noindent Let us suppose that the sectional curvature of $M$ is bounded above by
$-\epsilon<0$. Let $p_\pm$ be the fixed points of $\alpha$. There exists a unique oriented geodesic $\gamma$ such that $\gamma(-\infty)=p_-$ and $\gamma(+\infty)=p_+$. This geodesic is 
preserved by $\alpha$. Moreover, there exists $T_0\in\Bbb{R}$ such that, for all $t\in\Bbb{R}$:
$$
(\alpha\circ\gamma)(t) = \gamma(t+T_0).
$$
\noindent We define $\|\alpha\|$, the {\emph minimal displacement\/} of $\alpha$, by 
$\|\alpha\|=\left|T_0\right|$. We may show that:
$$
\|\alpha\| = \minf_{x\in M}d(x,x\cdot\alpha).
$$
\newsubhead{Conical Open Sets About Geodesics}
\noindent Let $P$ be a point in $\partial_\infty M$. For $q\in M\munion\partial_\infty M$, and
$r\in M$, we denote by $\widehat{prq}$ the angle at $r$ between the geodesics joining $r$ to $p$ and 
$r$ to $q$ respectively. For $\theta\in (0,\pi)$, we consider the neighbourhood $\Omega_\theta(p;r)$
of $p$ in $M\munion\partial_\infty M$ given by:
\subheadlabel{ConeControl}
$$
\Omega_\theta(p;r) = \left\{q\in M\munion\partial_\infty M\text{ s.t. }\widehat{prq}<\theta\right\}.
$$
\noindent We will require the following technical but elementary result:
\proclaim{Lemma \nextprocno}
\noindent Let $k>0$ be a positive real number. There exists a continuous function 
$\Delta_k:\Bbb{R}^+\times(0,\pi)\rightarrow(0,\pi)$ such that if:
\medskip
\myitem{(i)} $M$ is a Hadamard manifold with curvature less than or equal to $-k$, and
\medskip
\myitem{(ii)} $\gamma:\Bbb{R}\rightarrow M$ is a geodesic in  $M$,
\medskip
\noindent Then, for all $R\in\Bbb{R}^+$ and for all $t\in\Bbb{R}$:
$$
\Omega_\theta(\gamma(-\infty),\gamma(t))\subseteq
\Omega_{\Delta_k(R,\theta)}(\gamma(-\infty),\gamma(t+R)).
$$
\noindent Moreover:
\medskip
\myitem{(i)} $\Delta$ is strictly increasing in $\theta$ and strictly decreasing in $R$,
\medskip
\myitem{(ii)} $\Delta(R,\theta)$ tends to zero (locally uniformly in $\theta$) as $R$ tends to
infinity, and
\medskip
\myitem{(iii)} for all $R$ and for all $\theta$, $\Delta(R,\theta)<\theta$.
\endproclaim
\proclabel{LemmaConeControl}
\proof Let $M_k$ be the two dimensional Hadamard manifold of constant curvature equal to $-k$. Let
$XYZ$ be a triangle in $M_k\munion \partial_\infty M_k$ such that:
$$
X\in\partial_\infty M_k,\qquad YZ = R,\qquad \widehat{XYZ}=\pi-\theta.
$$
\noindent We define $\Delta_k(R,\theta)$ by:
$$
\Delta_k(R,\theta)=\widehat{YZX}.
$$
\noindent Let $q$ be a point in $\Omega_\theta(\gamma(-\infty);\gamma(t))$. Let us denote $A=\gamma(t)$,
$B=\gamma(t+R)$ and $C=q$. Let $A'B'C'$ be the comparison triangle in $M_k$ such that:
$$
A'B'=AB,\qquad A'C'=AC,\qquad \widehat{C'A'B'}=\widehat{CAB}.
$$
\noindent Let $C''$ be the point in $\partial_\infty M$ obtained by following the geodesic running from
$A'$ to $C'$ to its end point in $\partial_\infty M_k$. We have:
$$
\widehat{CBA}\leqslant \widehat{C'B'A'}\leqslant \widehat{C''B'A'}\leqslant \widehat{YZX} = \Delta_k(R,\theta).
$$
\noindent The result now follows.\qed
\medskip
\noindent We also require the following result
\proclaim{Lemma \nextprocno}
\noindent Let $p$ be a point in $\partial_\infty M$. Let $\Omega$ be a neighbourhood of $p$ in
$M\munion\partial_\infty M$. There exists a neighbourhood $U$ of $p$ in $\partial_\infty\Omega$ such that,
for $\gamma$ in $M$ a geodesic:
$$
\gamma(-\infty),\gamma(+\infty)\in U\Rightarrow\gamma\subseteq\Omega.
$$
\endproclaim
\proclabel{LemmaWeakConvexity}
\proof Indeed, otherwise, by constructing a sequence and taking limits, we obtain a geodesic $\gamma\in M$
such that:
$$
\gamma(0)\in\partial\Omega\minter M,\qquad \gamma(-\infty),\gamma(+\infty)=p.
$$
\noindent This is absurd.\qed
\newsubhead{Constructing Schottky Groups}
\noindent The now obtain the following result which shows how Schottky groups may be constructed in general:
\proclaim{Lemma \nextprocno}
\noindent Let $k>0$ be a real number and let $M$ be a Hadamard manifold with curvature less than
$-k$. Let $\alpha$ be a hyperbolic isometry of $M$ having $p_+$ and $p_-$ as attractive and 
repulsive fixed points respectively. Let $p_0$ be a point distinct from $p_\pm$. For every $B>0$, there exists a neighbourhood $\Omega$ of $p_0$ in $\partial_\infty M$ such that if $q_\pm\in\Omega$ and 
if $\beta$ is a hyperbolic isometry of $M$ having $q_\pm$ as fixed points such that $\|\beta\|>B$,
then the subgroup $\langle\alpha,\beta\rangle$ of $\opIsom(M)$ generated by $\alpha$ and $\beta$ is a Schottky group.
\endproclaim
\proclabel{LemmaSchottkyI}
\proof Let us also denote by $\alpha$ the unique geodesic in $M$
such that $\alpha(\pm\infty)=p_\pm$. Let $\Sigma_\alpha$ be the normal sphere bundle over $\alpha$. Let
$\overrightarrow{n}:\Sigma_\alpha\rightarrow\partial_\infty$ be the Gauss-Minkowski mapping. By Lemma
\procref{EmptyCylinder}, there exists a unique vector $X\in\Sigma_\alpha$ such that 
$\overrightarrow{n}(X)=p_0$. We suppose that $\alpha$ is normalised such that $X$ lies in the fibre
over $\alpha(0)$. We define $D=\|\alpha\|/2$. Let us denote by $\Sigma_{\alpha,D}^\pm$ the restriction
of $\Sigma_\alpha$ to $\alpha((\pm D,\pm\infty))$. Let us now define:
$$
U_\pm = \overrightarrow{n}(\Sigma_{\alpha,D}^\pm)\munion{p_\pm}.
$$
\noindent The mapping $\alpha$ sends the interior of the complement of $U_-$ onto $U_+$. By Lemma
\procref{LemmaConeControl}:
$$
U_\pm\subseteq \Omega_{\Delta_k(D,\pi/2)}(p_\pm;\alpha(0)).
$$
\noindent Since $\Delta_k(D,\pi/2)<\pi/2$, by continuity, there exists $\epsilon\in\Bbb{R}^+$ and
$\theta\in(0,\pi)$ such that, for all $t\in(-\epsilon,\epsilon)$:
$$
\Omega_{2\theta}(p_0;\alpha(t))\minter U_\pm = \emptyset.
$$
\noindent Let $U$ be a neighbourhood of $p_0$ in $M\munion\partial_\infty M$ such that:
$$
\Delta_k(d(\gamma,U),\pi-\theta)<\Delta_k(B/2,\pi/2).
$$
\noindent Let $V$ be a neighbourhood of $p_0$ in $\partial_\infty U$ as in Lemma 
\procref{LemmaWeakConvexity}. Let $q_\pm$ be points in $V$ and let $\beta$ be the unique geodesic in $M$
such that $\beta(\pm\infty)=q_\pm$.
\medskip
\noindent Let $\gamma:\Bbb{R}\rightarrow M$ be the unique geodesic meeting both $\alpha$ and $\beta$ at
right angles. By restricting $U$ if necessary, we may suppose that 
$\gamma(0)\in\alpha((-\epsilon,\epsilon))$ and that:
$$
\gamma(+\infty)\in\Omega_\theta(p_0;\gamma(0)).
$$
\noindent We may suppose that $\beta$ meets $\gamma$ at $\beta(0)$. Let us define $\Sigma^\pm_{\beta,B/2}$
as for $\Sigma^\pm_{\alpha,D}$ and:
$$
V_\pm=\overrightarrow{n}(\Sigma^\pm_{\beta,B/2})\munion\left\{q_\pm\right\}.
$$
\noindent The mapping $\beta$ sends the interior of the complement of $V_-$ onto $V_+$. Moreover:
$$
\Omega_{\Delta_k(B/2,\pi/2)}(\gamma(-\infty);\beta(0))\minter V_\pm =\emptyset.
$$
\noindent However:
$$
\Omega_\theta(\gamma(+\infty);\gamma(0))\minter U_\pm=\emptyset.
$$
\noindent Consequently:
$$\matrix
&U_\pm\hfill&\subseteq\Omega_{\pi-\theta}(\gamma(-\infty);\gamma(0))\hfill\cr
\Rightarrow\hfill&U_\pm\hfill&\subseteq\Omega_{\Delta_k(B/2,\pi/2)}(\gamma(-\infty);\beta(0)).\hfill\cr
\endmatrix$$
\noindent The sets $U_\pm$ and $V_\pm$ are thus disjoint, and the result now follows by the
classical ``Ping-Pong'' Lemma.\qed
\medskip
\noindent This yields the following useful corollary:
\proclaim{Corollary \nextprocno}
\noindent Let $k>0$ be a real number and let $M$ be a Hadamard manifold with curvature less than 
$-k$. Let $\Gamma$ be a subgroup of $\opIsom(M)$ which only contains hyperbolic elements
(and the identity). Let $\gamma$ be a hyperbolic element of $\Gamma$ and let $p_\pm$ be the attractive and
repulsive fixed points of $\gamma$ respectively. If $\alpha,\beta$ are elements of $\Gamma$ having no fixed
points in common with $\gamma$, then there exists $N>0$ such that, for $\left|n\right|\geqslant N$, the
subgroup of $\opIsom(M)$ generated by $\gamma^n\alpha\gamma^{-n}$ and by $\beta$ is a Schottky group.
\endproclaim
\proclabel{SchottkyGroups}
\proof Since $\alpha$ does not fix $p_-$, we find that the fixed points of $\gamma^{-n}\alpha\gamma^n$ tend
towards $p_+$ as $n$ tends towards $+\infty$. By hypothesis, $p_+$ is not a fixed point of $\beta$. Moreover, $\alpha$ is not equal to the identity and is thus hyperbolic.
Next, we see that:
$$
\|\gamma^{-n}\alpha\gamma^n\| = \|\alpha\| > 0.
$$
\noindent By the preceeding lemma, there exists $N\geqslant 0$ such that, for $n\geqslant N$, the subgroup
of $\opIsom(M)$ generated by $\gamma^{-n}\alpha\gamma^n$ and $\beta$ is a Schottky group. A similar reasoning
permits us to prove the result for $n$ negative, and the result now follows.\qed
\newsubhead{Proofs of Main Results of This Section}
\noindent We are now in a position to prove the propositions stated in the introduction to this section:
\medskip
{\noindent\bf Proof of Proposition \procref{NastyAlgebraI}:}
\medskip
\noindent\myitem{(1)} This holds for all $k$.
\medskip
\noindent\myitem{(2)} This holds for all but one values of $k$. There thus exists $K$ such that, for 
$\left|k\right|\geqslant K$, the application $\delta_k$ is hyperbolic.
\medskip
\noindent\myitem{(3)} For $\left|k\right|\geqslant K$, the mapping $\delta_k^n\alpha$ is hyperbolic for all but one values of $n$. Suppose that there exists $n\neq n'$ such that $\delta_k^n\alpha$ and 
$\delta_k^{n'}\alpha$ have the same fixed points as $\beta\alpha^k$. The mapping $\delta_k$ thus has the same
fixed points as $\beta\alpha^k$. Consequently, so does $\alpha$, and thus also $\beta$. This is absurd.
\medskip
\noindent\myitem{(4)} For $\left|k\right|\geqslant K$, this holds for all but one values of $n$, and the result now follows.\qed
\medskip
{\noindent\bf Proof of Proposition \procref{NastyAlgebraII}:\ }We begin by showing that, after exchanging 
$\xi$ and $\eta$ if necessary, there exists 
$K\in\Bbb{N}$ such that for $\left|k\right|\geqslant K$, the mapping $\delta_k=\eta\beta\alpha^k$ has no
fixed point in common with $\xi$. We suppose the contrary in order to obtain a contradiction. Suppose that there exists $k\neq k'$ such that $\eta\beta\alpha^k$ and $\eta\beta\alpha^{k'}$ both share the same fixed points
as $\xi$. Suppose, moreover, that $\xi\beta\alpha^k$ and $\xi\beta\alpha^{k'}$ both share the same fixed points as $\eta$. From the first condition, we conclude that $\alpha$ has the same fixed points as $\xi$, and
from the second we conclude that $\alpha$ has the same fixed points as $\eta$. It therefore follows that
both $\alpha$ and $\beta$ have the same fixed points, which is absurd. The result now follows.
\medskip
\noindent We now suppose that there exists $k\neq k'$ such that $\delta_k$ and $\delta_{k'}$
have the same fixed points as $\eta$. It follows that $\alpha$ does too, and thus so does $\beta$, which
is absurd. Consequently, there exists $K$ such that for $\left|k\right|\geqslant K$:
$$
\opFix(\delta_k)\minter \opFix(\eta)=\emptyset.
$$
\noindent In particular, $\delta_k$ is hyperbolic. Moreover, we may assume that, for all 
$\left|k\right|\geqslant K$:
$$
\opFix(\delta_k)\minter \opFix(\xi) = \emptyset.
$$
\noindent It follows by corollary \procref{SchottkyGroups} that there exists $N_1(k)$ such that, for 
$\left|n\right|\geqslant N_1(k)$ the group generated by $\langle \delta_k^n\xi\delta_k^{-n}, \eta \rangle$
is a Schottky group. The result now follows.\qed
\medskip
{\noindent\bf Proof of Proposition \procref{NastyAlgebraIV}:\ }If $\sigma=\opId$, then $J_k=\alpha\beta$ for
all $k$. This is necessarily hyperbolic, and the 
result follows in this case. Now suppose that $\sigma$ is hyperbolic. Suppose that there exists 
$k\neq k'$ such that $J_k=J_{k'}=\opId$. Then:
$$
\opId = J_kJ_{k'}^{-1}.
$$
\noindent If we denote $n=k-k'$, we obtain:
$$\matrix
& \sigma^n\alpha\sigma^{-n} \hfill&= \alpha \hfill\cr
\Rightarrow\hfill & \opFix(\alpha)\cdot\sigma^{-n} \hfill&= \opFix(\alpha) \hfill\cr
\Rightarrow\hfill &\opFix(\alpha)\hfill&=\opFix(\sigma).\hfill\cr
\endmatrix$$
\noindent Likewise, $\opFix(\beta) = \opFix(\sigma) = \opFix(\alpha)$, which is absurd. The result now
follows.\qed
\medskip
{\noindent\bf Proof of Proposition \procref{NastyAlgebraV}:\ }
\medskip
\noindent \myitem{(1)} Since $\langle\xi,\eta\rangle$ is a Schottky group:
$$\matrix
&\alpha\beta\alpha^{-1}\beta^{-1} \hfill&= \eta\xi \hfill&\neq\xi\hfill\cr
\Rightarrow\hfill &\alpha\beta\alpha^{-1}\beta^{-1}\xi^{-1} \hfill&\neq \opId\hfill\cr
\Rightarrow\hfill& \beta\alpha^{-1}\beta^{-1}\xi^{-1}\alpha\hfill& \neq \opId.\hfill\cr
\endmatrix$$
\noindent\myitem{(2)} We suppose that there exists $k\neq k'$ such that $\beta\alpha^k$ and 
$\beta\alpha^{k'}$ share the same fixed points. It follows that $\alpha$, and thus also $\beta$, share the
same fixed points as these mappings, which is absurd. Consequently, there exists $K_1$ such that, for 
$\left|k\right|\geqslant K_1$, the application $\gamma_k$ has no fixed points in common with $\xi$ or 
$\beta\alpha^{-1}\beta^{-1}\xi^{-1}\alpha$. By conjugation, for $\left|k\right|\geqslant K_1$, 
$\alpha\gamma_k\alpha^{-1}$ has no fixed point in common with $\alpha\beta\alpha^{-1}\beta^{-1}\xi^{-1}$. Moreover, for $\left|k\right|\geqslant K_1$, the application $\delta_k = \xi\gamma_k$ does not have any fixed
points in common with $\gamma_k$.
\medskip
\noindent\myitem{(3)} If $\opFix(\delta_k) = \opFix(\eta)$, then $\delta_k$ commutes with $\eta$. Thus, $\langle \xi, \delta_k^{-n}\eta\delta_k^n \rangle = \langle \xi, \eta\rangle$
is a Schottky group for all $n$. Otherwise, $\opFix(\delta_k)\minter \opFix(\eta)=\emptyset$. Since 
$\opFix(\delta_k)\minter \opFix(\xi)\neq\emptyset$, it follows by corollary \procref{SchottkyGroups} that there exists $N_1(k)$ such that, for $\left|n\right|\geqslant N_1(k)$, the group 
$\langle \xi, \delta_k^{-n}\eta\delta_k^n \rangle$ is a Schottky group.
\medskip
\noindent\myitem{(4)} If $\delta_k$ shares the same fixed points as $\alpha\gamma_k\alpha^{-1}$, then so 
does $\alpha\beta\alpha^{-1}\beta^{-1}\xi^{-1}$, which is absurd by $(2)$. If $\delta_k$ shares the same
fixed points as $\gamma_k$, then so does $\xi$, which is also absurd by $(2)$. It thus follows by
corollary \procref{SchottkyGroups} that there exists $N_2(k)\geqslant N_1(k)$ such that, for 
$\left|n\right|\geqslant N_2(k)$, the group 
$\langle \alpha\gamma_k\alpha^{-1}, {\delta_k^n}^{-1}\gamma_k^{-1}\delta_k^n \rangle$ is a Schottky group, and the result now follows by conjugation with $\alpha$.\qed
\newhead{The Trouser Decomposition}
\newsubhead{Introduction}
\noindent Throughout this section $(M,Q)$ denotes a pointed, compact, Riemannian manifold of strictly
negative sectional curvature, $(\tilde{M},\tilde{Q})$ denotes its universal cover, $(\Sigma,P)$ denotes
a pointed, compact surface of genus at least $2$ and $\theta:\pi_1(\Sigma,P)\rightarrow\pi_1(M,Q)$
denotes a representation of $\pi_1(\Sigma,P)$ in $\pi_1(M,Q)$. In \cite{GallKapMard}, Gallo, Kapovich and Marden obtain a trouser decomposition
of $\Sigma$ such that the $\theta$-image of the fundamental group of each trouser is a Schottky group.
In this section, we will prove the following lemma, which is the generalisation to our framework of this part of Gallo, Kapovich and Marden's result:
\proclaim{Lemma \nextprocno}
\noindent Let $(M,Q)$ be a pointed, compact, three dimensional Riemannian manifold of strictly negative
sectional curvature. Let $(\Sigma,\theta)$ be a pointed, orientable, compact surface of hyperbolic type 
(i.e. of genus greater than or equal to $2$) with holonomy in $\pi_1(M,Q)$. Then there exist
bound, marked trousers $(T_i,\theta_i,\beta_i)_{1\leqslant i\leqslant 2g-2}$ with holonomy in 
$\pi_1(M,Q)$ such that:
\medskip
\myitem{(i)} for all $i$, the image of $\theta_i$ is a Schottky group, and
\medskip
\myitem{(ii)}
$$
(\Sigma,\theta)\cong\mathcircle_{i=1}^g(\munion_{i=1}^{2g-2}(T_i,\theta_i,\beta_i)).
$$
\endproclaim
\proclabel{SchottkyTrouserDecomposition}
\noindent The terms and terminology used in the statement of this lemma will be explained in the second and third parts of this section. The proof of this lemma is an immediate consequence of Propositions
\procref{DecompositionStepI}, \procref{DecompositionStepII} and \procref{DecompositionStepIII}.
\medskip
\noindent In the second part of this section, we define the notion of a marked surface with holonomy in a
group $G$. In the third part, we define the notion of a binding, which is required as a bookkeeping measure
in order to recover correctly the fundamental group after the surface has been cut up into trousers. In the
remaining parts, we successively decompose the surface $\Sigma$, proving Propositions
\procref{DecompositionStepI}, \procref{DecompositionStepII} and \procref{DecompositionStepIII}.
\newsubhead{Marked Surfaces and Handles}
\noindent The group $\pi_1(M,Q)$ may be considered as a subgroup of the group of isometries of $\tilde{M}$.
Since $M$ is compact, every element of $\pi_1(M,Q)\setminus\left\{\opId\right\}$ has a hyperbolic
action over $(\tilde{M},\tilde{Q})$. If $a$ and $b$ are two points in $\partial_\infty\tilde{M}$, and
if $\gamma\in\pi_1(M,Q)$, then $\gamma$ cannot send $a$ to $b$ and $b$ to $a$, since, otherwise, it would be
elliptic. Moreover, we recall that the fixed point sets of two elements of $\pi_1(M,Q)$ are either
disjoint or equal.
\medskip
\noindent Let $G$ be a group. We define a {\emph marked surface\/} with {\emph holonomy\/} in $G$ to be a triplet of:
\medskip
\myitem{(i)} a pointed compact surface with boundary $(\Sigma,\partial\Sigma,P_0)$, with $\partial\Sigma$
oriented such that the interior of $\Sigma$ lies to its left,
\medskip
\myitem{(ii)} for every connected component $C_\alpha$ of $\partial\Sigma$, a point $Q_\alpha\in C_\alpha$,
which is the base point of $C_\alpha$, and
\medskip
\myitem{(iii)} a homomorphism $\theta:\pi_1(\Sigma,P_0)\rightarrow G$.
\medskip
\noindent By slight abuse of notation, we will refer to the marked surface with holonomy simply by
$(\Sigma,\theta)$. Every closed curve in $\Sigma$ defines canonically a conjugacy class in 
$\pi_1(\Sigma,P_0)$. For such a closed curve, $a$, we denote the corresponding conjugacy class by $[a]$.
The $\theta$-image of $[a]$ defines a conjugacy class in $G$ which we denote by $\theta[a]$. If 
$G$ is a subgroup of $\pi_1(M,Q)$, and if one element in a given conjugacy class of $G$ is hyperbolic,
then all elements in that conjugacy class are hyperbolic, and the class is said to be {\emph hyperbolic\/}. A closed curve $a$ in $(\Sigma,\theta)$ is then said to be {\emph hyperbolic\/} if and only if $\theta[a]$ is hyperbolic.
\medskip
\noindent Every pair $(a,b)$ of simple closed curves in $\Sigma$ with the same base point
generates a conjugacy class of pairs in $\pi_1(\Sigma,P_0)$. Indeed, with $x$ a curve joining $P_0$ to
the common base point of $a$ and $b$, we define:
$$
(a,b)_x = (x^{-1}ax,x^{-1}bx) \in\pi_1(\Sigma,P_0)\times\pi_1(\Sigma,P_0).
$$
\noindent The conjugacy class of this pair is independant of $x_0$, and we denote it by $[a,b]$.
\medskip
\noindent The $\theta$-image of the class $[a,b]$ defines a conjugacy class of pairs in $G$ which we
denote by $\theta[a,b]$. Every element of $\theta[a,b]$ generates a subgroup of $G$, and $\theta[a,b]$ thus defines a conjugacy class of subgroups of $G$ which we denote by $\langle\theta[a,b]\rangle$. If $G$ is a subgroup of $\pi_1(M,Q)$, a subgroup in a given conjugacy class in $G$ is a Schottky group if and only
if every group in that class is too, and the class is said to be {\emph Schottky\/}. The pair $(a,b)$ is
said to be {\emph Schottky\/} if and only if $\langle\theta[a,b]\rangle$ is Schottky. We define a 
{\emph non-elementary\/} pair and a non-elementary class in the same manner.
\medskip
\noindent Let $(\Sigma,\theta)$ be a marked surface with holonomy in $G$. A {\emph marked handle\/} in
$(\Sigma,\theta)$ is a pair $(a,b)$ of non-dividing simple closed curves in $(\Sigma,\theta)$ such that $a$ and $b$ only intersect at their common base point, and $b$ crosses $a$ from right to left. 
When $G$ is a subgroup of the group $\pi_1(M,Q)$, we say that $(a,b)$ is
a {\emph hyperbolic handle\/} if and only if both $a$ and $b$ are hyperbolic. We say that the handle 
$(a,b)$ is {\emph non-elementary\/} or {\emph Schottky\/} if and only if, viewed as a pair, it is
non-elementary or Schottky respectively.
\medskip
\noindent Finally, for $(\Sigma,\theta)$ a marked surface with boundary with holonomy in $G$, and for a given family, $(Q_i)_{i\in F}$ of points in $(\Sigma,\theta)$, we define a {\emph sash\/} to be a
set of simple curves $(\gamma_i)_{i\in F}$ in $(\Sigma,\theta)$ indexed by the same set as 
$(Q_i)_{i\in F}$ such that:
\medskip
\myitem{(i)} for every $i$, $\gamma_i$ goes from $P_0$ to $Q_i$,
\medskip
\myitem{(ii)} for every $i\neq j$, $\gamma_i$ only intersects $\gamma_j$ at $P_0$,
\medskip
\myitem{(iii)} for every $i$, $\gamma_i$ can only intersect $\partial\Sigma$ at $Q_i$.
\medskip
\noindent Sashes permit us to identify objects in $\Sigma$ (curves, marked handles, etc.)
explicitly with elements or pairs of elements of $\pi_1(\Sigma,P_0)$.
\medskip
\noindent When we define many such geometric structures over $(\Sigma,\theta)$, we will assume them to be disjoint, except, possibly, at their end or base points.
\newsubhead{Bindings}
\noindent Let $G$ be a group. Let $(\Sigma,\theta)$ be a marked surface with holonomy in $G$. Suppose that
$\Sigma$ has $n$ boundary components, and let $(C_i, Q_i)_{1\leqslant i\leqslant n}$ be the pointed oriented
boundary components of $\Sigma$. For every $i$, denote the set of all homotopy classes of curves in $\Sigma$
going from $P_0$ to $Q_i$ by $\pi_1(\Sigma,P_0,Q_i)$. A {\emph binding\/} of $(\Sigma,\theta)$ is a set of
mappings $(\beta_i)_{1\leqslant i\leqslant n}$ such that:
\medskip
\myitem{(i)} for every $i$, $\beta_i$ maps $\pi_1(\Sigma,P_0,Q_i)$ into $G$, and
\medskip
\myitem{(ii)} for every $i$, $\beta_i$ is equivariant under $\theta$ with respect to the canonical right 
action of $\pi_1(\Sigma,P_0)$ on $\pi_1(\Sigma,P_0,Q_i)$. In otherwords, for all 
$\xi\in\pi_1(\Sigma,P_0,Q_i)$, and for all $\eta\in\pi_1(\Sigma,P_0)$:
$$
\beta_i(\xi\eta) = \beta_i(\xi)\theta(\eta).
$$
\noindent We will denote the binding by $\beta$, and, for every $i$, we will write $\beta$ instead of $\beta_i$ when there is no danger of ambiguity. A marked surface $(\Sigma,\theta)$ carrying a binding, $\beta$, will be said to be {\emph bound\/}. We will denote the bound marked surface with holonomy
in $G$ by $(\Sigma,\theta,\beta)$. For any oriented, pointed boundary component $(C,Q)$ of $\Sigma$, we 
define the element $\beta(C)\in G$ such that, for all $\xi\in\pi_1(\Sigma,P_0,Q)$:
$$
\beta(C) = \beta(\xi)\theta(\xi^{-1}C\xi)\beta(\xi)^{-1}.
$$
\noindent Let $(\Sigma,\theta,\beta)$ and $(\Sigma',\theta',\beta')$ be bound marked surfaces with holonomy
in $G$. Let $(C,Q)$ and $(C',Q')$ be oriented, pointed boundary components of $\Sigma$ and $\Sigma'$ 
respectively. Let $\varphi:(C,Q)\rightarrow({C'}^{-1},Q')$ be an orientation reversing homeomorphism. 
Let $\Sigma\cup_\varphi\Sigma'$ be the surface obtained by joining $\Sigma$ to $\Sigma'$ along $\varphi$.
We identify all objects (points, curves, etc.) in $\Sigma$ and $\Sigma'$ with the corresponding objects
in $\Sigma\munion_\varphi\Sigma'$ and vice-versa. We take the base point of $\Sigma\cup_\varphi\Sigma'$
to be that of $\Sigma$.
\medskip
\noindent We say that $(\theta,\beta)$ and $(\theta',\beta')$ may be {\emph joined\/} along $\varphi$ if and
only if $\beta(C)^{-1} = \beta'(C')$. In this case, we define $\theta\cup_\varphi\theta'$, the {\emph joined union\/} of $\theta$ and $\theta'$ {\emph along\/} $\varphi$,
such that:
\medskip
\myitem{(i)} for all $x\in\pi_1(\Sigma,P_0)$:
$$
(\theta\cup_\varphi\theta')(x) = \theta(x),
$$
\noindent and
\medskip
\myitem{(ii)} for all $x\in\pi_1(\Sigma',P'_0)$, $\xi\in\pi_1(\Sigma,P_0,Q)$ and 
$\eta\in\pi_1(\Sigma',P_0',Q')$:
$$
(\theta\cup_\varphi\theta')(\xi^{-1}\eta x\eta^{-1}\xi) = 
\beta(\xi)^{-1}\beta'(\eta)\theta'(x)\beta'(\eta)^{-1}\beta(\xi).
$$
\noindent The join condition on $(\theta,\beta)$ and $(\theta',\beta')$ ensures that 
$\theta\cup_\varphi\theta'$ is a well defined homomorphism from $\pi_1(\Sigma\cup_\varphi\Sigma',P_0)$ into
$G$. We may also define $\beta\cup_\varphi\beta'$, the {\emph joined union\/} of $\beta$ and $\beta'$ {\emph along\/} $\varphi$ in the same manner, with the join condition again ensuring that 
$\beta\munion_\varphi\beta'$ is well defined. In the sequel, the homeomorphism $\varphi$ will be suppressed, and we thus denote:
$$
(\Sigma,\theta,\beta)\cup(\Sigma',\theta',\beta')=
(\Sigma\cup\Sigma',\theta\cup\theta',\beta\cup\beta').
$$
\noindent In a similar manner, we may join a bound marked surface to itself along an orientation
reversing homeomorphism between two boundary components. We denote the join of $(\Sigma,\theta,\beta)$
to itself by:
$$
{\mathcircle}(\Sigma,\theta,\beta) = ({\mathcircle}\Sigma,{\mathcircle}\theta,{\mathcircle}\beta).
$$
\noindent For surfaces with many boundary components, this process may be iterated, and we will refer to the
$n$-fold join of $(\Sigma,\theta,\beta)$ to itself by:
$$
{\mathcircle}^n(\Sigma,\theta,\beta).
$$
\newsubhead{Slicing Open The Surface}
\noindent In this subsection, we will prove the following result:
\proclaim{Proposition \nextprocno}
\noindent Let $(M,Q_0)$ be a pointed, compact, three dimensional manifold of strictly negative sectional curvature. let $(\Sigma,\theta)$ be a compact surface without boundary of genus $g\geqslant 2$ with holonomy in $\pi_1(M,Q_0)$.
\medskip
\noindent If the image of $\theta$ is non-elementary, then there exists 
$(\Sigma_{2g-2},\theta_{2g-2},\beta_{2g-2})$, a bound marked surface with holonomy in $\pi_1(M,Q_0)$ of genus
$1$ with a non-elementary marked handle $(a,b)$ and having $2g-2$ boundary components, such that:
\medskip
\myitem{(i)} every connected component of $\partial\Sigma_{2g-2}$ is hyperbolic, and
\medskip
\myitem{(ii)} $(\Sigma,\theta)$ is homeomorphic to 
${\mathcircle}^{g-1}(\Sigma_{2g-2},\theta_{2g-2},\beta_{2g-2})$.
\endproclaim
\proclabel{DecompositionStepI}
\noindent If $\Sigma$ is a closed surface of genus $g\geqslant 2$, then we define a
{\emph canonical basis\/} of $\Sigma$ to be an n-tuple $[(a_1,b_1), ..., (a_n, b_n)]$ of pairs of simple closed curves in $\Sigma$ which correspond to handles in $\Sigma$, such that:
\medskip
\noindent\myitem{(i)} for all $i\neq j$ the pairs $(a_i,b_i)$ and $(a_j,b_j)$ are disjoint,
\medskip
\noindent\myitem{(ii)} for all $i$, the curves $a_i$ and $b_i$ have a common base point, which is their only shared point and which we denote by $Q_i$, and
\medskip
\noindent\myitem{(iii)} for all $i$, $b_i$ crosses $a_i$ from right to left.
\medskip
\noindent We have the following result:
\proclaim{Proposition \nextprocno}
\noindent Let $(M,Q_0)$ be a pointed compact three dimensional Hadamard manifold of strictly negative 
sectional curvature. Let $(\Sigma,\theta)$ be a compact surface without boundary of genus $g\geqslant 2$
with holonomy in $\pi_1(M,Q_0)$. Let $[(a_1,b_1),...,(a_n,b_n)]$ be a canonical basis of $\Sigma$.
\medskip
\noindent If $\theta$ is non-elementary, then there exists a homeomorphism $\Phi$ of $\Sigma$ such that
$\Phi_*(a_1,b_1)$ is non-elementary.
\endproclaim
\proclabel{ObtainingANonElemHandle}
\proof Since the image of $\theta$ is non-elementary, after composing with a homeomorphism of $\Sigma$ 
which permutes the generators of $\pi_1(\Sigma)$, we may suppose that the class $\theta(a_1)$ is hyperbolic. 
For all $i$, let $Q_i$ be the common base point of $a_i$ and $b_i$. Let $(c_i)_{1\leqslant i\leqslant g}$ be
a sash of $(Q_i)_{1\leqslant i\leqslant g}$. Let us define $a=c_1^{-1}a_1c_1$ and 
$b=c_1^{-1}b_1c_1$. Let us next define $\alpha = \theta(a)$ and $\beta = \theta(b)$.
There are two cases:
\medskip
\noindent\myitem{(1)} $\beta$ does not have any fixed points in common with $\alpha$. In this case, we take
$\Phi=\opId$ and we obtain the desired result.
\medskip
\noindent\myitem{(2)} $\beta$ has two fixed points in common with $\alpha$ (in particular, $\beta$ could be
the identity). Since $\Gamma$ is non-elementary, there exists 
$x\in\left\{c_i^{-1}a_ic_i,c_i^{-1}b_ic_i | 2 \leqslant i \leqslant n\right\}$ such that
$\xi = \theta(x)$ does not have any fixed points in common with $\alpha$.
\medskip
\noindent We may orient $x$ such that $C=ax$ is freely homotopic to a simple closed curve. Let $T$ be the 
Dehn twist about $C$. In order to simplify the discussion, we identify $Q_1$ with $P_0$. Figure 
\figref{FigureDehnTwistI} shows the handle of $\Sigma$ defined by $(a,b)$ opened up by cutting along $a$ 
and $b$. 
\placefigure{}{%
\placelabel[3.6][0.3]{$P_0$}%
\placelabel[0.1][4]{$P_0$}%
\placelabel[3.6][4]{$P_0$}%
\placelabel[0.1][0.3]{$P_0$}%
\placelabel[0.63][0]{Rest of $\Sigma$}%
\placelabel[2.55][3.95]{$b$}%
\placelabel[2.55][0.3]{$b$}%
\placelabel[3.7][2.6]{$a$}%
\placelabel[0.1][2.6]{$a$}%
\placelabel[3.1][2.6]{$a'$}%
\placelabel[2.25][3.1]{$x_A$}%
\placelabel[2.25][1.3]{$x_B$}%
\placelabel[1.3][1.3]{$x_C$}%
\placelabel[1.3][3.1]{$x_D$}%
}{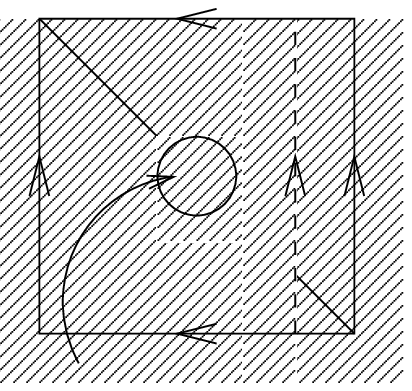}
\figlabel{FigureDehnTwistI}
\noindent The curve $x$ terminates at $P_0$ without crossing either $a$ or $b$. This corresponds to one of
the four topologically distinct configurations $x_A$,$x_B$,$x_C$ or $x_D$ shown in figure 
\figref{FigureDehnTwistI}. By defining $C$
using a slightly displaced copy of $a$, we may ensure that $C$ never crosses $a$. If $x$ is in one of the
configurations $A$ or $B$, this would correspond to taking a copy of $a$ displaced to its left, and if $x$
is in one of the configurations $C$ or $D$, then this would correspond to taking a copy of $a$ displaced to
its right, as shown in figure \figref{FigureDehnTwistI}. If we choose $x$ in configuration $A$, then the Dehn twist, $T$, about $C$ satisfies:
$$
T_*a = a, \qquad T_*b = bax.
$$
\noindent Now:
$$
\opFix(\beta\alpha\xi) = \opFix(\alpha) \Leftrightarrow \opFix(\xi) = \opFix(\alpha).
$$
\noindent It thus follows that $\opFix(\beta\alpha\xi)\minter \opFix(\alpha) =\emptyset$. Defining $\Psi=T$,
we obtain the desired homeomorphism.\qed
\medskip
\noindent We next have:
\proclaim{Proposition \nextprocno}
\noindent Let $(M,Q_0)$ be a pointed compact three dimensional Hadamard manifold of strictly negative
sectional curvature. Let $(\Sigma,\theta)$ be a compact surface without boundary of genus $g\geqslant 2$
with holonomy in $\pi_1(M,Q_0)$. Let $[(a_1,b_1),...,(a_g,b_g)]$ be a canonical basis of $\Sigma$.
\medskip
\noindent If $(a_1,b_1)$ is non-elementary, then, for all $i$ different to $1$, there exists a homeomorphism
$\Psi$ of $\Sigma$ such that:
\medskip
\noindent\myitem{(1)} $\Psi_*(a_1,b_1)$ is non-elementary,
\medskip
\noindent\myitem{(2)} $\Psi_*a_i$ and $\Psi_*b_i$ are hyperbolic, and
\medskip
\noindent\myitem{(3)} for all $j$ different from $1$ and $i$, the application $\Psi_*$ 
acts on $a_j$ and $b_j$ by conjugation.
\endproclaim
\proclabel{MakingOtherHandlesHyperbolic}
\noindent Using induction, we thus obtain the following result:
\proclaim{Corollary \nextprocno}
\noindent With the same hypotheses as in Proposition \procref{MakingOtherHandlesHyperbolic}, there exists a
homeomorphism $\Psi$ such that:
\medskip
\noindent\myitem{(1)} $\Psi_*(a_1,b_1)$ is non-elementary, and
\medskip
\noindent\myitem{(2)} for all $i$ different from $1$, $\Psi_*a_i$ and $\Psi_*b_i$ are
hyperbolic.
\endproclaim
\proclabel{CorMakingOtherHandlesHyperbolic}
{\noindent\bf Proof of Proposition \procref{MakingOtherHandlesHyperbolic}:\ } For all $i$, let $Q_i$ be the
common base point of $a_i$ and $b_i$. Let $(c_i)_{1\leqslant i\leqslant g}$ be a sash of 
$(Q_i)_{1\leqslant i\leqslant g}$. As
before, we denote $a=c_1^{-1}a_1c_1$ and $b=c_1^{-1}b_1c_1$. We now define $\alpha = \theta(a)$, and 
$\beta = \theta(b)$. For a given $i$ different to $1$, we denote $x = c_i^{-1}a_ic_i$, and 
$y = c_i^{-1}b_ic_i$, and we define $\xi = \theta(x)$ and $\eta = \theta(y)$. 
\medskip
\noindent We orient $y$ such that $C=yb$ is freely homotopic to a simple closed curve. We denote
$d_0=yb$. Let $T_0$ be a Dehn twist about $C$. The curves $c_1$ and $c_i$ may be chosen such that:
$$\matrix
(T_0)_*a = d_0a, \hfill& (T_0)_*b = b, \hfill\cr
(T_0)_*x = d_0x, \hfill& (T_0)_*y = y. \hfill\cr
\endmatrix$$
\noindent Since $(T_0)_*$ leaves $d_0$ invariant, we obtain, for all $n$:
$$\matrix
(T_0^n)_*a = d_0^na, \hfill& (T_0^n)_*b = b, \hfill\cr
(T_0^n)_*x = d_0^nx, \hfill& (T_0^n)_*y = y. \hfill\cr
\endmatrix$$
\noindent Moreover, for all $j$ different from $1$ and $i$, since $d_0$ only intersects $c_j$ near $P_0$, it
follows that $(T_0)_*$ acts on $(a_j,b_j)$ by conjugation. Let $T_a$ be the Dehn twist about $a$ such that,
for all $k$:
$$\matrix
(T_a^k)_*a = a, \hfill& (T_a^k)_*b = ba^k, \hfill\cr
(T_a^k)_*x = x, \hfill& (T_a^k)_*y = y. \hfill\cr
\endmatrix$$
\noindent For all $j$ different to $1$, since $a$ stays away from $(a_j,b_j)$, $(T_a)_*$ has no effect on
this pair. We define $\Psi_{k,n}:\Sigma\rightarrow\Sigma$ by:
$$
\Psi_{k,n} = T_a^kT_0^nT_a^{-k}.
$$
\noindent If, for all $k$, we define $d_k$ by $d_k=yba^k$, then we obtain:
$$\matrix
(\Psi_{k,n})_*a = d_k^na, \hfill& (\Psi_{k,n})_*ba^k = ba^k, \hfill\cr
(\Psi_{k,n})_*x = d_k^nx, \hfill& (\Psi_{k,n})_*y = y. \hfill\cr
\endmatrix$$
\noindent Where $d_k = yba^k$. Moreover, for all $j$ different to $1$ and $i$, the mapping $\Psi_{k,n}$
acts on $(a_j,b_j)$ by conjugation. We now choose $k$ and $n$ as in Lemma \procref{NastyAlgebraI}. If $\eta$
is hyperbolic, then, taking $\Psi=\Psi_{k,n}$, we obtain the desired result. Otherwise, $\eta=\opId$, and
we may take a Dehn twist $T$ about $x$ such that $Ty=xy$ and which leaves all the other curves invariant. We see that:
$$
\theta(\Psi_{k,n})_*T_*y = \theta{\Psi_{k,n}}_*xy = \theta d_k^nxy=\delta^n_k\xi,
$$
\noindent and this is hyperbolic. We thus denote $\Psi = \Psi_{k,n}\circ T$ and we once again obtain the
desired result.\qed
\medskip
\noindent The proof of Proposition \procref{DecompositionStepI} is now elementary:
\medskip
{\noindent\bf Proof of Proposition \procref{DecompositionStepI}:\ } Let
$[(a_1,b_1),...,(a_n,b_n)]$ be a canonical basis of $\Sigma$. By Proposition 
\procref{ObtainingANonElemHandle} and corollary \procref{CorMakingOtherHandlesHyperbolic}, we may suppose 
that $(a_1,b_1)$ is a non-elementary marked handle and that, for all $2\leqslant i\leqslant n$, 
$a_i$ and $b_i$ are hyperbolic.
\medskip
\noindent We obtain $\Sigma_{2g-2}$ by cutting $\Sigma$ along each $a_i$ for $2\leqslant i\leqslant n$. We
identify objects (points, curves etc.) in $\Sigma$ with the corresponding objects in $\Sigma_{2g-2}$, and
vice-versa. For each $i$, we denote by $C_{i,-}$ the copy of $a_i$ in $\Sigma_{2g-2}$ from which $b_i$
leaves and by $C_{i,+}$ the copy of $a^{-1}_i$ in $\Sigma_{2g-2}$ at which $b_i$ arrives, and we define 
$Q_{i,\pm}$ to be the copy of $Q_i$ lying in $C_{i,\pm}$. Let $(\gamma_i)_{2\leqslant i\leqslant n}$ be 
a sash of $(Q_{i,-})_{2\leqslant i\leqslant n}$. We define the the binding $\beta_{2g-2}$ over
$\Sigma_{2g-2}$ such that, for all $i$:
$$\matrix
\beta_{2g-2}(\gamma_i) \hfill&= \opId,\hfill\cr
\beta_{2g-2}(b_i\gamma_i) \hfill&= \theta(\gamma_i^{-1}b_i\gamma_i).\hfill\cr
\endmatrix$$
\noindent We may verify that this binding satisfies the appropriate join conditions and that:
$$
(\Sigma,\theta) \cong {\mathcircle}^{g-1}(\Sigma_{2g-2},\theta_{2g-2},\beta_{2g-2}).
$$
\noindent Finally, we define the marked handle $(a,b)$ in $\Sigma_{2g-2}$ by:
$$
(a,b) = (a_1,b_1).
$$
\noindent Since $(a_1,b_1)$ is hyperbolic, so is $(a,b)$, and the result now follows.\qed
\newsubhead{Pruning}
\noindent In this subsection, we aim to prove the following result:
\proclaim{Proposition \nextprocno}
\noindent Let $(M,Q)$ be a pointed, compact, three dimensional manifold of strictly negative sectional
curvature. Let $(\Sigma_{2g-2},\theta_{2g-2},\beta_{2g-2})$ be a bound marked surface with
holonomy in $\pi_1(M,Q)$ of genus $1$ having $2g-2$ boundary components, all of which are hyperbolic, and a non-elementary marked handle $(a,b)$. There exists:
\medskip
\myitem{(i)} a bound marked surface $(\Sigma_2,\theta_2,\beta_2)$ with holonomy in $\pi_1(M,Q)$, and
\medskip
\myitem{(ii)} $2g-4$ bound, marked trousers $(T_i,\varphi_i,\beta_i)_{1\leqslant i\leqslant 2g-4}$ with
holonomy in $\pi_1(M,Q)$,
\medskip
\noindent such that:
$$
(\Sigma_{2g-2},\theta_{2g-2},\beta_{2g-2})\cong (\Sigma_2,\varphi_2,\beta_2)\munion_{i=1}^{2g-4}
(T_i,\varphi_i,\beta_i),
$$
\noindent and:
\medskip
\myitem{(i)} $(\Sigma_2,\theta_2,\beta_2)$ has genus $1$, two hyperbolic boundary components and a non-elementary marked handle $(a',b')$, and
\medskip
\myitem{(ii)} for all $i$, the image of $\varphi_i$ in $\pi_1(M,Q)$ is a Schottky group.
\endproclaim
\proclabel{DecompositionStepII}
\noindent We begin with the following proposition:
\proclaim{Proposition \nextprocno}
\noindent Let $(M,Q)$ be a pointed compact three dimensional manifold of strictly negative sectional
curvature. Let $(\Sigma_n,\theta_n,\beta_n)$ be a bound marked surface with holonomy in $\pi_1(M,Q)$ of
genus $1$ having $n$ hyperbolic boundary components and a non-elementary marked handle $(a,b)$.
\medskip
\noindent Let $x_1$ and $y_1$ be two boundary components of $\Sigma_n$. Let $c$ be a simple curve joining 
$x$ to $y$ which is disjoint from $x_1$, $y_1$, $a_1$ and $b_1$ except possibly at its end points.
\medskip
\noindent There exists a homeomorphism 
$\Psi:(\Sigma_n,\partial \Sigma_n,P_0)\rightarrow (\Sigma_n,\partial \Sigma_n,P_0)$ of $\Sigma_n$ such that:
\medskip
\noindent\myitem{(1)} $\Psi_*(a_1,b_1)$ is non-elementary,
\medskip
\noindent\myitem{(2)} if $z$ is a boundary component of $\Sigma$ different to $x_1$ and $y_1$, then 
$\Phi_*$ acts on $z$ by conjugation (thus, if $\theta(z)$ is hyperbolic, then $\theta\Phi_*z$ is
hyperbolic), and
\medskip
\noindent\myitem{(3)} $\Psi_*(x_1,c^{-1}y_1c)$ is Schottky (and thus, in particular, $x_1c^{-1}y_1c$ is hyperbolic).
\endproclaim
\proclabel{MakingTrousersSchottky}
\proof Let $Q_0$ be the common base point of $a_1$ and $b_1$. Let $Q_x$ and $Q_y$ be the base points of 
$x$ and $y$ resp. Let $(c_0,c_x)$ be a sash of $(Q_0,Q_x)$. We define:
$$\matrix
a = c_0^{-1}a_1c_0, \hfill& b = c_0^{-1}b_1c_0, \hfill\cr
x = c_x^{-1}x_1c_x, \hfill& y = (cc_x)^{-1}y_1(cc_x). \hfill\cr
\endmatrix$$
\noindent We denote $\alpha = \theta_n(a)$, $\beta = \theta_n(b)$, $\xi = \theta_n(x)$ and
$\eta = \theta_n(y)$. 
\medskip
\noindent We may suppose that $b$ is oriented such that $C=yb$ is homotopic to a simple closed
curve. we denote $d_0=yb$. Let $T_0$ be a Dehn Twist about $C$. We may suppose that $c_a$ and $c_x$ are
chosen such that:
$$\matrix
(T_0)_*a = d_0a, \hfill& (T_0)_*b = b, \hfill\cr
(T_0)_*x = d_0xd_0^{-1}, \hfill& (T_0)_*y = y. \hfill\cr
\endmatrix$$
\noindent Since $(T_0)_*$ is a homomorphism which leaves $d_0$ invariant, we obtain, for all $n$:
$$\matrix
(T_0^n)_*a = d_0^na, \hfill& (T_0^n)_*b = b, \hfill\cr
(T_0^n)_*x = d_0^nxd_0^{-n}, \hfill& (T_0^n)_*y = y. \hfill\cr
\endmatrix$$
\noindent We choose a Dehn twist $T_a$ about $a$ such that, for all $k$:
$$\matrix
(T_a^k)_*a = a, \hfill& (T_a^k)_*b = ba^k, \hfill\cr
(T_a^k)_*x = x, \hfill& (T_a^k)_*y = y. \hfill\cr
\endmatrix$$
\noindent We define $\Psi_{k,n}:\Sigma_{2n-2}\rightarrow\Sigma_{2n-2}$ by:
$$
\Psi_{k,n} = T_a^kT_0^nT_a^{-k}
$$
\noindent Denoting $d_k=yba^k$, We have:
$$\matrix
(\Psi_{k,n})_*a = d_k^na, \hfill& (\Psi_{k,n})_*ba^k = ba^k, \hfill\cr
(\Psi_{k,n})_*x = d_k^nxd_k^{-n}, \hfill& (\Psi_{k,n})_*y = y. \hfill\cr
\endmatrix$$
\noindent We choose $k$ and $n$ as in Lemmata \procref{NastyAlgebraI} and \procref{NastyAlgebraII}, and we 
denote $\Psi = \Psi_{k,n}$. The first result follows from the first and third conclusions of Lemma
\procref{NastyAlgebraI}. Since $\Psi$ is a product of Dehn twists, it acts on paths corresponding to
boundary components by conjugation. The second result now follows. The third result follows from Lemma
\procref{NastyAlgebraII}.\qed
\medskip
\noindent This yields:
\proclaim{Corollary \nextprocno}
\noindent With the same hypotheses as in Proposition \procref{MakingTrousersSchottky}, there exists:
\medskip
\myitem{(i)} a bound marked trouser $(T,\theta,\beta)$ with holonomy in $\pi_1(M,Q)$, and
\medskip
\myitem{(ii)} a bound marked surface $(\Sigma_{n-1},\theta_{n-1},\beta_{n-1})$ with holonomy in $\pi_1(M,Q)$,
\medskip
\noindent such that:
$$
(\Sigma_n,\theta_n,\beta_n) \cong (\Sigma_{n-1},\theta_{n-1},\beta_{n-1})\cup(T,\theta,\beta),
$$
\noindent and,
\medskip
\myitem{(i)} $(\Sigma_{n-1},\theta_{n-1},\beta_{n-1})$ has genus $1$, $n-1$ hyperbolic boundary components, 
and a non-elementary marked handle $(a',b')$,
\medskip
\myitem{(ii)} the image of $\theta$ in $\pi_1(M,Q)$ is a Schottky group.
\endproclaim
\proclabel{CorMakingTrousersSchottky}
\proof We may suppose that $(\Sigma,\theta)$ satisfies the conclusions of Proposition 
\procref{MakingTrousersSchottky}. Cutting $\Sigma_n$ along a curve freely homotopic to $x_1c^{-1}y_1c$
which does not intersect either $a$ or $b$, we obtain the desired result as in the proof of
Proposition \procref{DecompositionStepI}.\qed
\medskip
\noindent We may now prove Proposition \procref{DecompositionStepII}:
\medskip
{\noindent\bf Proof of Proposition \procref{DecompositionStepII}:\ } This follows directly from
induction and corollary \procref{CorMakingTrousersSchottky}.\qed
\newsubhead{Untying The Root}
%
%
%
%
\noindent In this subsection, we aim to prove the following result:
\proclaim{Proposition \nextprocno}
\noindent Let $(M,Q)$ be a pointed compact three dimensional manifold of strictly negative sectional
curvature. Let $(\Sigma_2,\theta_2,\beta_2)$ be a bound marked surface with holonomy in $\pi_1(M,Q)$ of
genus $1$ having $2$ hyperbolic boundary components and a non-elementary marked handle $(a,b)$.
\medskip
\noindent There exist bound marked trousers $(T_i,\theta_i,\beta_i)_{1\leqslant i\leqslant 2}$ with
holonomy in $\pi_1(M)$ such that:
\medskip
\myitem{(i)} $(\Sigma_2,\theta_2,\beta_2)\cong [\mathcircle (T_1,\theta_1,\beta_1)]\cup(T_2,\theta_2,\beta_2)$, and
\medskip
\myitem{(ii)} the images of $\theta_1$ and $\theta_2$ are Schottky groups.
\endproclaim
\proclabel{DecompositionStepIII}
\noindent We begin by proving the following proposition:
\proclaim{Proposition \nextprocno}
\noindent Let $(M,Q)$ be a pointed compact three dimensional manifold of strictly negative sectional
curvature. Let $(\Sigma_2,\theta_2,\beta_2)$ be a bound marked surface with holonomy in $\pi_1(M,Q)$ of
genus $1$ having $2$ boundary components and a non-elementary marked handle $(a,b)$.
\medskip
\noindent Let $x_1$ and $y_1$ be the two oriented boundary components of $\Sigma_2$, and let $c$ be a curve 
in $\Sigma_2$ joining the base point of $y$ to that of $x$, but otherwise disjoint from $a_1$,$b_1$,$x_1$ and
$y_1$. There exists a homeomorphism $\Psi$ of $\Sigma_2$ such that $\Psi_*(c^{-1}x_1c,y_1)$ and 
$\Psi_*(a_1^{-1}b_1a_1,b_1)$ are Schottky.
\endproclaim
\proclabel{TheFinalCut}
\proof Let $Q_0$ be the common base point of $a$ and $b$. Let $Q_x$ and $Q_y$ be the base points of $x$ and
$y$ respectively. Let $(c_0,c_x,c_y)$ be a sash of $(Q_0,Q_x,Q_y)$. We denote:
$$\matrix
a = c_a^{-1}a_1c_a, \hfill& b = c_a^{-1}b_1c_a, \hfill\cr
x = c_x^{-1}x_1c_x, \hfill& y = c_y^{-1}y_1c_y. \hfill\cr
\endmatrix$$
\noindent We may assume that $c = a^{-1}b^{-1}ab$ is freely homotopic to $(yx)^{-1}$. We denote:
$$
\alpha = \theta_2(a), \beta = \theta_2(b), \xi = \theta_2(x), \eta = \theta_2(y), \sigma = \theta_2(c).
$$
\noindent By Proposition \procref{MakingTrousersSchottky}, we may suppose that $\langle \xi, \eta \rangle$
is a Schottky group and that $\langle\alpha,\beta\rangle$ is non-elementary. Let $T_c$ be the Dehn twist
about $c$ such that, for all $m$:
%
%
%
%
%
%
$$\matrix
T_c^ma = a, \hfill& T_c^mb = b, \hfill\cr
T_c^mx = c^mxc^{-m}, \hfill& T_c^my = c^myc^{-m}. \hfill\cr
\endmatrix$$
\noindent By Proposition \procref{NastyAlgebraIV}, after replacing $\theta$ with $\theta\circ(T_c)_*$ if necessary, we may suppose that $J=\xi\beta$ is not equal to the identity. For $k\in\Bbb{N}$, we define $d_k=xba^k$. We may assume that $xb$ is freely homotopic to a simple
closed curve, and, as in the proofs of Lemmata \procref{DecompositionStepI} and 
\procref{DecompositionStepII}, using Dehn twists, for all $k$ and for all $n$, we may construct a
homeomorphism $\Psi_{k,n}$ of $\Sigma_2$ such that:
$$\matrix
(\Psi_{k,n})_*a \hfill&= d_k^na, \hfill& (\Psi_{k,n})_*ba^k \hfill&= ba^k, \hfill\cr
(\Psi_{k,n})_*x \hfill&= x, \hfill& (\Psi_{k,n})_*y \hfill&= d_k^{-n}yd_k^n. \hfill\cr
\endmatrix$$
\noindent Where $d_k = xba^k$. Choosing $k$ and $n$ as in Lemma \procref{NastyAlgebraV}, we obtain the
desired result.\qed
\medskip
\noindent We now obtain Proposition \procref{DecompositionStepIII}:
\medskip
{\noindent\bf Proof of Proposition \procref{DecompositionStepIII}:\ }This follows directly by applying
Proposition \procref{TheFinalCut} and then cutting along curves freely homotopic to $yx$ and to $b$.\qed
\newhead{Invariant Domains of Schottky Groups}
\newsubhead{Introduction}
\noindent Throughout this section, $(M,Q)$ will be a pointed three dimensional Hadamard manifold
and $\Gamma=\langle\alpha,\beta\rangle$ will be a Schottky subgroup of $\opIsom(M,Q)$. 
In this section, we study the algebraic properties of Schottky groups.
\medskip
\noindent We define an {\emph invariant domain\/} of $\Gamma$ to be a Jordan domain $\Omega$ contained in 
$\partial_\infty M$ which is invariant under the action of $\Gamma$. Invariant domains are easy to
construct. Indeed, let $\Gamma'$ be any Schottky subgroup of $\psl(2,\Bbb{C})$ which preserves the real
line. Let $\phi:\Gamma\rightarrow\Gamma'$ be an isomorphism, and let 
$\Phi:\partial_\infty M\rightarrow\hat{\Bbb{C}}$ be a homeomorphism which intertwines with
$\phi$ (which we construct in appendix \headref{MakingHomeomorphism}). $\Phi^{-1}(\Bbb{\hat{R}})$ is a
Jordan curve in $\partial_\infty M$, and both connected
components of its complement are invariant domains. 
\medskip
\noindent Trivially, for any invariant domain, $\Omega$:
$$
\opFix(\Gamma)\subseteq\partial\Omega.
$$
\noindent We define $\hat{\Omega}$ by:
$$
\hat{\Omega} = \overline{\Omega}\setminus\opFix(\Gamma).
$$
\noindent For every element $\gamma$ of $\Gamma$ we define $\gamma_-$ and $\gamma_+$ to be the 
repulsive and attractive fixed points of $\gamma$ respectively. We say that $\gamma$ is simple when
there exists no other element $\eta$ in $\Gamma$ and no $n\geqslant 2$ such that $\gamma=\eta^n$.
For $\gamma$ a simple element, we say that a connected component $I$ of 
$\partial\Omega\setminus\opFix(\Gamma)$ is {\emph adapted to\/} $\gamma$ if and only if, viewed as an oriented subarc of $\partial\Omega$, it runs from $\gamma_-$ to $\gamma_+$ (and thus has $\gamma_\pm$ as 
its endpoints).
\medskip
\noindent Such a connected component, when it exists, is unique. We say that $\Omega$ is {\emph adapted\/}
to the generators $(\alpha,\beta)$ when both of $\alpha$ and $\beta$ have adapted components in
$\partial\Omega$. Trivially, the invariant domain $\Omega$, constructed as above is adapted to the
generators $(\alpha,\beta)$.
\medskip
\noindent We show (Proposition \procref{QuotientOfInvDomainExists}) that $\Gamma$ acts properly
discontinuously over $\hat{\Omega}$ and that its quotient is a compact topological surface with boundary. When $\Omega$ is adapted to the generators $(\alpha,\beta)$, the adapted components of $\alpha$ and $\beta$ project down to boundary components of $\hat{\Omega}/\Gamma$, and since $\alpha$ and $\beta$ cannot be conjugate in $\Gamma$, these components are distinct. Thus, since its fundamental group is the free group on two generators, $\hat{\Omega}/\Gamma$ is a trouser.
\medskip
\noindent This would already allow us to construct equivariant Plateau problems over a given bound, marked
trouser $(T,\theta,\beta)$ with holonomy in $\pi_1(M,Q)$ such that $\theta$ is a Schottky group.
The solution is by no means unique, and in order to solve the more general problem, we are required to study
in more depth the algebraic properties of adapted invariant domains.
\medskip
\noindent If $\gamma$ is a hyperbolic element of $\Gamma$, we define the torus $\Bbb{T}_\gamma$ by:
$$
\Bbb{T}_\gamma = (\partial_\infty\tilde{M}\setminus\left\{\gamma_\pm\right\})/\langle\gamma\rangle.
$$
\noindent For any hyperbolic element $\gamma$ of $\Gamma$, we canonically define a preferred one dimensional subspace $L_\gamma$ of $H_1(\Bbb{T}_\gamma)$. We may also associate to any invariant domain 
$\Omega$ of $\Gamma$ a unique element $[\Omega]_\gamma$ in $H_1(\Bbb{T}_\gamma)$. This element necessarily
lies in $L_\gamma$ (see section \subheadref{HomologicalClasses}).
\medskip
\noindent Let $\opHomeo_0(\partial_\infty M)$ denote the connected component of the space of
homeomorphisms of $\partial_\infty M$ which contains the identity. There exists a canonical 
homomorphism of $\pi_1(M,Q)$ into $\opHomeo_0(\partial_\infty M)$. We may show (see
corollary \procref{FundamentalGroupOfHomeo}) that if $\opTildeHomeo_0(\partial_\infty M)$ 
denotes the universal cover of $\opHomeo_0(\partial_\infty M)$, then it is a two-fold covering.
For every hyperbolic element $\gamma$ of $\Gamma$, we may canonically define a mapping
$\opLift_\gamma:L_\gamma\rightarrow\opTildeHomeo_0(\partial_\infty M)$.
\medskip
\noindent The main result of this section, which allows us to establish the obstruction to constructing
a $\pi_1(M,Q)$ structure over a closed surface $\Sigma$ with holonomy in $\pi_1(M,Q)$ may now be
expressed as follows:
\proclaim{Lemma \nextprocno}
\noindent Let $\Gamma=\langle\alpha,\beta\rangle$ be a Schottky group. Let us denote $\gamma=\alpha\beta$.
Let $([a]_\alpha,[b]_\beta,[c]_\gamma)$ be a triplet in $L_\alpha\times L_\beta\times L_\gamma$. There exists
an invariant domain, $\Omega$, of $\Gamma$ in $\partial_\infty M$ adapted to the generators 
$(\alpha,\beta)$ such that:
$$
[\Omega]_\alpha = [a]_\alpha,\qquad [\Omega]_\beta = [b]_\beta,\qquad [\Omega]_\gamma = [c]_\gamma,
$$
\noindent if and only if:
$$
\opLift_\gamma([c]_\gamma)^{-1}\opLift_\alpha([a]_\alpha)\opLift_\beta([b]_\beta) \neq \opId.
$$
\endproclaim
\proclabel{PrincipalObstructionToJoin}
\noindent The proof of this result follows immediately from Propositions \procref{LemmaBraguette} and
\procref{InverseZipper}.
\medskip
\noindent In the second part of this section, we review the geometric properties of invariant domains, 
justifying the assertions made in this introduction. In the third part, we review the topology of the group
of homeomorphisms of the sphere and the properties of braid groups of order three in the sphere. In the
fourth part, we define the element $[\Omega]_\gamma$, for any hyperbolic element, $\gamma$, of 
$\Gamma$ and any invariant domain, $\Omega$, of $\Gamma$. Finally, in the fifth part, we prove Propositions \procref{LemmaBraguette} and \procref{InverseZipper}.
\newsubhead{Invariant Domains}
\noindent In this section, we study the geometric properties of invariant domains. The following result is fairly trivial:
\proclaim{Propositon \nextprocno}
\noindent $\Gamma$ acts properly discontinuously on $\hat{\Omega}$ and $\hat{\Omega}/\Gamma$ is compact.
\endproclaim
\proclabel{QuotientOfInvDomainExists}
\noindent Since the fundamental group of $\hat{\Omega}/\Gamma$ is isomorphic to the free group on two
generators, $\hat{\Omega}/\Gamma$ is either a trouser or a punctured torus. We consequently restrict our
attention to invariant domains which are adapted to these generators, since, in this case, the adapted
components of $\alpha$ and $\beta$ project down to distinct boundary components of $\hat{\Omega}/\Gamma$ and
$\hat{\Omega}/\Gamma$ is consequently a trouser.
\medskip
\noindent In the sequel, we will require the following technical result concerning adapted invariant 
domains:
\proclaim{Proposition \nextprocno}
\noindent Let $\Omega$ be an invariant domain of $\Gamma$ adapted to the generators $(\alpha,\beta)$. Let 
$\gamma$ and $\delta$ be the elements of $\Gamma$ defined by:
$$
\gamma=\alpha\beta,\qquad \delta=\beta\alpha.
$$
\noindent If $\alpha_\pm$,$\beta_\pm$,$\gamma_\pm$ and $\delta_\pm$ are the fixed points of these four elements, then they are distributed round $\partial\Omega$ in the following order:
$$
\alpha_-,\alpha_+,\delta_+,\delta_-,\beta_-,\beta_+,\gamma_+,\gamma_-.
$$
\endproclaim
\proclabel{OrderOfFixedPoints}
\proof Let $p$ be one of the fixed points of $\gamma$. We have:
$$
(p\cdot\alpha)\cdot\delta = p\cdot(\alpha\beta\alpha) = (p\cdot\gamma)\cdot\alpha = p\cdot\alpha.
$$
\noindent Consequently, $\alpha$ sends the fixed points of $\gamma$ to those of $\delta$, and so
$\gamma_\pm\cdot\alpha = \delta_\pm$. Likewise $\delta_\pm\cdot\beta = \gamma_\pm$. Let $I_\alpha$ and
$I_\beta$ be the adapted components of $\alpha$ and $\beta$ respectively. 
\medskip
\noindent Since $\alpha$ preserves orientation, it also preserves $I_\alpha$ and $I_\alpha^C$, Moreover, it
shifts all points in $I_\alpha$ and $I_\alpha^C$ towards $\alpha^+$. Consequently, the attractive and
repulsive fixed points of $\alpha$,$\gamma$ and $\delta$ appear in the following order as one moves around
$\partial\Omega$ in the positive direction:
$$
\alpha_-,\alpha_+,\delta_\pm,\gamma_\pm.
$$
\noindent Similarly, the mapping $\beta$ shifts all points in $I_\beta$ and $I_\beta^C$ towards $\beta^+$.
Thus, the attractive and repulsive fixed points of $\beta$,$\gamma$ and $\delta$ appear in the following order as one moves around $\partial\Omega$ in the positive direction:
$$
\beta_-,\beta_+,\gamma_\pm,\delta_\pm.
$$
\noindent Since the mapping $\alpha$, which sends the fixed points of $\gamma$ to those of $\delta$, is
an orientation preserving mapping, it follows that, as one moves round $\partial\Omega$ in the positive direction, $\gamma_-$ and $\gamma_+$ appear in the same order as $\delta_-$ and $\delta_+$.
Moreover, since $I_\alpha$ and $I_\beta$ are adapted components, they have no fixed points of $\Gamma$ in
their interiors, and so, combining all this information, we find that these eight fixed points appear round
$\partial\Omega$ in the following order:
$$
\alpha_-,\alpha_+,\delta_\pm,\delta_\mp,\beta_-,\beta_+,\gamma_\pm,\gamma_\mp.
$$
\noindent Let $p$ be the fixed point of $\gamma$ lying closest to $\beta_+$ in $I_\beta^C$. Let $q$ be the
fixed point of $\delta$ lying closest to $\alpha_+$ in $I_\alpha^C$. The point $p$ is also the fixed
point of $\gamma$ lying closest to $\alpha_+$ in $I_\alpha^C$. Since $\alpha$ preserves orientation, we
have $p\cdot\alpha=q$. The points $\alpha_\pm$ lie between $q$ and $\beta_+$ in $I_\beta^C$. Consequently,
the points $\alpha_\pm\cdot\gamma=\alpha_\pm\cdot\beta$ lie between $p=p\cdot\gamma=q\cdot\beta$ and 
$\beta_+$ in $I_\beta^C$.
\medskip
\noindent Let $I_1$ and $I_2$ be the two connected components of 
$\partial\Omega\setminus\left\{\gamma_\pm\right\}$. The mapping $\gamma$ shifts all points in the intervals
$I_1$ and $I_2$ towards $\gamma_+$. We may suppose that $I_1$ is the component containing $\alpha_\pm$.
Since $\alpha_+\cdot\gamma$ lies between $p$ and $\beta_+$ in $I_\beta^C$, it follows that it also lies
between $\alpha_+$ and $p$ in $I_1$. The point $p$ is consequently the attractive fixed point of $\gamma$
and the result now follows.\qed
\newsubhead{Braids and the Topology of $\opHomeo_0(\partial_\infty M)$}
\noindent For any topological surface $\Sigma$, let $\opHomeo_0(\Sigma)$ be the connected component of
$\opHomeo(\Sigma)$ that contains the identity. There exists a canical embedding of the group $\pi_1(M,Q)$
into $\opHomeo_0(\partial_\infty M)$. We may thus consider the mapping $\theta$ as a homomorphism of 
$\pi_1(\Sigma,P_0)$ taking values in the group $\opHomeo_0(\partial_\infty M)$.
\medskip
\noindent Let $X\subseteq Y$ be topological spaces. We define a {\emph strong deformation retraction\/} of 
$Y$ onto $X$ to be a mapping $\psi:I\times Y\rightarrow Y$ such that:
\medskip
\myitem{(i)} $\psi_0:Y\rightarrow Y$ is the identity,
\medskip
\myitem{(ii)} $\psi_1(Y)\subseteq X$, and
\medskip
\myitem{(iii)} for all $t\in I$, the restriction of $\psi_t$ to $X$ is the identity.
\medskip
\noindent Let us denote by $S^2\subseteq\Bbb{R}^3$ the sphere of radius $1$ in $\Bbb{R}^3$. We recall the
following result concerning the homotopy type of $\opHomeo_0(S^2)$ (see \cite{Frib}, \cite{LeRoux}):
\proclaim{Theorem \nextprocno\ {\bf [Friberg, 1973]}}
\noindent The space $\opHomeo_0(S^2)$ retracts by strong deformations onto $\opSO(3,\Bbb{R})$.
\endproclaim
\proclabel{ThmFriberg}
\noindent In particular, we obtain:
\noskipproclaim{Corollary \nextprocno}
$$
\pi_1(\opHomeo_0(\partial_\infty M),\opId) \cong \pi_1(\opSO(3,\Bbb{R}),\opId) \cong \Bbb{Z}_2.
$$
\endnoskipproclaim
\proclabel{FundamentalGroupOfHomeo}
\noindent A {\emph braid of order $3$\/} in $\partial_\infty M$ is a triple 
$\pmb{\gamma}=(\gamma_1,\gamma_2,\gamma_3)$ where $\gamma_1,\gamma_2,\gamma_3:I\rightarrow S^2$ are
curves such that, for all $t\in I$, the points $\gamma_1(t)$, $\gamma_2(t)$ and $\gamma_3(t)$ are all 
distinct. The interested reader may find a more detailed treatment of braids in general in \cite{SmiF} or
in appendix $D$ of \cite{SmiE}. For all $t\in I$, we denote:
$$
\pmb{\gamma}(t) = (\gamma_1(t),\gamma_2(t),\gamma_3(t)).
$$
\noindent We call the point $\pmb{\gamma}(0)$ the {\emph base point\/} of the braid $\pmb{\gamma}$. We
say that the braid is closed if and only if $\pmb{\gamma}(0) = \pmb{\gamma}(1)$. Let $\pmb{\gamma}_0$ and 
$\pmb{\gamma}_1$ be two braids having the same endpoints. A 
{\emph homotopy\/} between $\pmb{\gamma}_0$ and $\pmb{\gamma}_1$ is a continuous family 
$(\pmb{\eta}_t)_{t\in I}$ of braids having the same extremities as $\pmb{\gamma}_0$ and 
$\pmb{\gamma}_1$ such that:
$$
\pmb{\eta}_0 = \pmb{\gamma}_0,\qquad \pmb{\eta}_1 = \pmb{\gamma}_1.
$$
\noindent For ${\mathbf{p}}=(p_1,p_2,p_3)$ a triplet of distinct points in $\partial_\infty M$, we
denote by $T_{{\mathbf{p}}}$ the family of braids of order $3$ in $\partial_\infty M$ having 
${\mathbf{p}}$ as a base point. Likewise, we denote by $T_{{\mathbf{p}}}^0$ the subfamily of $T_{{\mathbf{p}}}$
consisting of all the closed braids in $T_{{\mathbf{p}}}$. Let $\sim$ be the homotopy equivalence relation
over $T_{{\mathbf{p}}}$. The law of composition of curves yields the law of composition of braids, and the set $T_{{\mathbf{p}}}/\sim$ thus forms a semigroup. Likewise, $T_{{\mathbf{p}}}^0/\sim$ forms a group. For 
$\pmb{\gamma}=(\gamma_1,\gamma_2,\gamma_3)$ a braid in $T_{{\mathbf{p}}}$, we denote by
$[\pmb{\gamma}]=[\gamma_1,\gamma_2,\gamma_3]$ its projection in $T_{{\mathbf{p}}}/\sim$. 
\medskip
\noindent Let $C_0(I,\opHomeo_0(\partial_\infty M))$ denote the family of continuous curves in
$\opHomeo_0(\partial_\infty M)$ leaving from the identity. For any triplet 
${\mathbf{p}}=(p_1,p_2,p_3)$ of distinct points in $\partial_\infty M$, we define the mapping 
$\Cal{T}_{{\mathbf{p}}}$ from $C_0(I,\opHomeo_0(\partial_\infty M))$ to $T_{{\mathbf{p}}}$ such that,
for all $c\in C_0(I,\opHomeo_0(\partial_\infty M))$, and for all $t\in I$:
$$
(\Cal{T}_{{\mathbf{p}}}c)(t) = (p_1\cdot c(t),p_2\cdot c(t),p_3\cdot c(t)).
$$
\noindent In the case where $\partial_\infty M=\hat{\Bbb{C}}$, for all ${\mathbf{p}}$, the mapping
$\Cal{T}_{{\mathbf{p}}}$ restricts to a homeomorphism from $C_0(I,\psl(2,\Bbb{C}))$ onto $T_{{\mathbf{p}}}$.
Consequently, in the general case, the mapping $\Cal{T}_{{\mathbf{p}}}$ is surjective. Moreover, for all
${\mathbf{p}}$, the mapping $\Cal{T}_{{\mathbf{p}}}$ quotients down to a surjective mapping from 
$\opTildeHomeo_0(\partial_\infty M)$ to $T_{{\mathbf{p}}}/\sim$ which we also denote by
$\Cal{T}_{{\mathbf{p}}}$. Theorem \procref{ThmFriberg} permits us to show that this mapping is bijective.
\medskip
\noindent Let $\pmb{\gamma}$ be a braid in $T_{{\mathbf{p}}}$. We may suppose that there
exists a point $p_\infty$ in $\partial_\infty M$ which does not lie in the image of 
$\pmb{\gamma}$. Let $\alpha:\partial_\infty M\setminus\left\{p_\infty\right\}\rightarrow\Bbb{R}^2$
be a homeomorphism. For any closed curve $\eta$ in $\Bbb{R}^2\setminus\left\{0\right\}$, let $\opWind(\eta)$
be the winding number of $\eta$ about $0$. The quantity:
$$
\sum_{i<j}\opWind(\alpha\circ\gamma_i-\alpha\circ\gamma_j)\text{ Mod }2
$$
\noindent is well defined and independant of $\alpha$ and $p_\infty$ (see, for example, \cite{SmiF}, or appendix $D$ of \cite{SmiE}). We thus define $\opWind_r(\pmb{\gamma})$, the {\emph relative winding 
number\/} of $\pmb{\gamma}$ to be equal to this quantity.
\medskip
\noindent $\opWind_r$ defines an isomorphism from $T^0_{{\mathbf{p}}}/\sim$ to $\Bbb{Z}_2$. We define the
mapping $W_{r,\pmb{p}}$ which sends $\pi_1(\opHomeo_0(\partial_\infty M))$ into $\Bbb{Z}_2$ by:
$$
W_{r,\pmb{p}} = \opWind_r\circ\Cal{T}_{{\mathbf{p}}}.
$$
\noindent The mapping $W_{r,\pmb{p}}$ defines an isomorphism from $\pi_1(\opHomeo_0(\partial_\infty M))$ to 
$\Bbb{Z}_2$. Moreover, it does not depend on the choice of ${\mathbf{p}}$, and we thus denote 
$W_r=W_{r,\pmb{p}}$.
\newsubhead{Homological Classes}
\noindent For $\gamma\in\Gamma$, we define the torus $\Bbb{T}_\gamma$ by:
\subheadlabel{HomologicalClasses}
$$
\Bbb{T}_\gamma = (\partial_\infty M\setminus\left\{\gamma_-,\gamma_+\right\})/\langle\gamma\rangle.
$$
\medskip
\noindent Let $C$ be a closed curve which turns once around the cylinder 
$\partial_\infty M\setminus\left\{\gamma_\pm\right\}$ in such a manner that $\gamma_+$ lies to its
left (in otherwords, in its interior). $C$ defines a unique homology class $[C]_\gamma$ in 
$H_1(\Bbb{T}_\gamma)$. Using Poincare duality, we define $L_\gamma\subseteq H_1(\Bbb{T}_\gamma)$ by:
$$
L_\gamma = \left\{[a]\in H_1(\Bbb{T}_\gamma)\ |\ \langle [C_\gamma],[a] \rangle = 1 \right\}.
$$
\noindent Heuristically, $L_\gamma$ is a straight line which contains all the curves in $\Bbb{T}_\gamma$
which cross $C_\gamma$ exactly once, going from right to left. For $[a],[b]\in L_\gamma$, we define
$[a]-[b]$ to be the unique integer $n$ such that:
$$
[a] = [b] + n[C]_\gamma.
$$
\noindent Conversely, for $[a]\in L_\gamma$ and $n\in\Bbb{Z}$, we define $[a]+n\in L_\gamma$ such that $([a]+n)-[a]=n$.
\medskip
\noindent Let us denote by $\pi_\gamma$ the canonical projection from 
$\partial_\infty M\setminus\left\{\gamma_-,\gamma_+\right\}$ onto $\Bbb{T}_\gamma$. Let $\Omega$ be an
invariant domain of $\Gamma$ in $\partial_\infty M$. If $c:I\rightarrow\Omega$ is a curve such that
$c(0)=\gamma_-$ and $c(1)=\gamma_+$, then the homology class of $\pi_\gamma\circ c$ in $\Bbb{T}_\gamma$
lies in $L_\gamma$. Moreover, it does not depend on $c\subseteq\Omega$, and we thus define:
$$
[\Omega]_\gamma = [\pi_\gamma\circ c] \in L_\gamma.
$$
\newsubhead{Liftings of Applications}
\noindent Let $\gamma$ be an element in $\Gamma$. Let $\gamma_\pm$ be the fixed points of $\gamma$. Let 
$c:I\rightarrow\partial_\infty M\setminus\left\{\gamma_\pm\right\}$ be a curve such that:
$$
c(0)\cdot\gamma = c(1).
$$
\noindent The curve $c$ projects down to an element of $L_\gamma$. We consider $\gamma_\pm$ as being
constant curves, and we define $\hat{\gamma}_c$ to be the unique lifting of $\gamma$ in 
$\opTildeHomeo_0(\partial_\infty M)$ such that:
$$
\Cal{T}_{(\gamma_-,c(0),\gamma_+)}\hat{\gamma}_c = [\gamma_-,c,\gamma_+].
$$
\noindent Trivially, $\hat{\gamma}_c$ only depends on the homotopy class of $[\pi_\gamma\circ c]$ in
$L_\gamma\subseteq H_1(\Bbb{T}_\gamma)$. We thus obtain a mapping 
$\opLift_\gamma:L_\gamma\rightarrow\opTildeHomeo_0(\partial_\infty M)$ which is defined such
that, for all $c$ with $c(0)\cdot\gamma=c(1)$:
$$
\opLift_\gamma([\pi_\gamma\circ c]) = \hat{\gamma}_c.
$$
\noindent We have:
\proclaim{Proposition \nextprocno}
\noindent Let $c,c':I\rightarrow M\setminus\left\{\gamma_\pm\right\}$ be curves such that 
$c(0)\cdot\gamma=c(1)$ and $c'(0)\cdot\gamma=c'(1)$. Then:
$$
\opLift_\gamma([\pi_\gamma\circ c]) = \opLift_\gamma([\pi_\gamma\circ c']) 
\Leftrightarrow [\pi_\gamma\circ c] - [\pi_\gamma\circ c'] = 0
\text{ Mod }2.
$$
\endproclaim
\proclabel{ConditionForEquivalenceOfLifts}
\proof Let us define:
$$
\alpha = \opLift_\gamma([\pi_\gamma\circ c])^{-1}\cdot\opLift_\gamma([\pi_\gamma\circ c']).
$$
\noindent The element $\alpha$ is a lifting of the identity in $\opTildeHomeo_0(\partial_\infty M)$, and:
$$
\Cal{T}_{(\gamma_-,c'(0),\gamma_+)}\alpha=[\gamma_-,c^{-1}c',\gamma_+].
$$
\noindent The element $\alpha$ is equal to the identity if and only if:
$$
W_r[\gamma_-,c^{-1}c',\gamma_+]=0.
$$
\noindent However, since $\gamma_\pm$ are trivial paths, we obtain:
$$
W_r[\gamma_-,c^{-1}c',\gamma_+]=[\pi_\gamma\circ c']-[\pi_\gamma\circ c]\text{ Mod }2.
$$
\noindent The result now follows.\qed
\medskip
\noindent If $\Omega$ is an invariant domain of $\Gamma$ in $\partial_\infty M$, then, for every
element $\gamma\in\Gamma$, we define $\hat{\gamma}_\Omega$ by:
$$
\hat{\gamma}_\Omega = \opLift_\gamma[\Omega]_\gamma.
$$
\noindent We obtain the following result:
\proclaim{Proposition \nextprocno}
\noindent Let $\Gamma=\langle\alpha,\beta\rangle$ be a Schottky group. Let $\Omega$ be an invariant domain of $\Gamma$ in $\partial_\infty M$ which is adapted to the generators $(a,b)$. Let us denote 
$\gamma=\alpha\beta$, and let us define $\Delta$ by:
$$
\Delta = \hat{\gamma}_\Omega^{-1}\hat{\alpha}_\Omega\hat{\beta}_\Omega.
$$
\noindent Then $\Delta$, which is a closed curve in $\opHomeo_0(\partial_\infty M)$, is homotopically
non-trivial.
\endproclaim
\proclabel{LemmaBraguette}
\proof By Lemma \procref{HomeoEquivalence}, it suffices to prove this result when
$\partial_\infty\tilde{M}=\hat{\Bbb{C}}$ and $\Gamma$ is a subgroup of $\psl(2,\Bbb{C})$. In this case, for
any triplet ${\mathbf{p}}=(p_1,p_2,p_3)$ of distinct points in $\hat{\Bbb{C}}$, the restriction of the
mapping $\Cal{T}_{{\mathbf{p}}}$ to $C_0(I,\psl(2,\Bbb{C}))$ is bijective. Let $\Cal{S}_{{\mathbf{p}}}$ denote
its inverse. Since, for any braid, $\pmb{\eta}$ in $T_{{\mathbf{p}}}$, the point ${\mathbf{p}}$ is the
base point of $\pmb{\eta}$, and is thus determined by $\pmb{\eta}$, we may write $\Cal{S}$ instead of
$\Cal{S}_{{\mathbf{p}}}$.
\medskip
\noindent $\partial\Omega$ is homeomorphic to the circle $S^1$ and is invariant under the action of $\Gamma$
on $\Bbb{C}$. By Proposition \procref{OrderOfFixedPoints}, since $\Omega$ is adapted to the generators $(\alpha,\beta)$, the $8$ points $\alpha_\pm$,$\beta_\pm$,$\gamma_\pm$ and 
$\gamma_\pm\cdot\alpha$ are distributed around $\partial\Omega$ according to the diagram in figure 
\figref{FigureDistributionOfFixedPoints}.
%
\placefigure{}{%
\placelabel[0.75][3.6]{$\alpha_-$}%
\placelabel[3.5][0.75]{$\beta_+$}%
\placelabel[-0.65][0.75]{$\gamma_+\cdot\alpha$}%
\placelabel[-0.05][0.1]{$\gamma_-\cdot\alpha$}%
\placelabel[2.8][0.1]{$\beta_-$}%
\placelabel[2.8][3.6]{$\gamma_-$}%
\placelabel[3.5][2.85]{$\gamma_+$}%
\placelabel[2][2]{$\Omega$}%
\placelabel[0][2.85]{$\alpha_+$}%
}{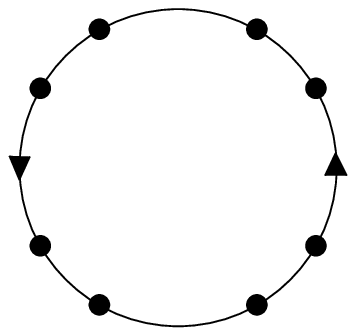}
\figlabel{FigureDistributionOfFixedPoints}
\noindent For every point, $p$ in $\partial_\infty M$, let us denote also by $p$ the constant curve
which sends the unit interval onto $p$. Let $p_0$ be a point in $\Omega$. Let $a$ be a curve in $\Omega$
joining $p_0$ to $p_0\cdot\alpha$. Let $b$ be a curve in $\Omega$ joining $p_0\cdot\alpha$ to 
$p_0\cdot(\alpha\beta)$. By definition:
$$\matrix
\hat{\alpha}_\Omega \hfill&= [\Cal{S}(\alpha_-,a,\alpha_+)], \hfill\cr
\hat{\beta}_\Omega \hfill&= [\Cal{S}(\beta_-,b,\beta_+)], \hfill\cr
\hat{\gamma}_\Omega^{-1} \hfill&= [\Cal{S}(\gamma_-,(a\cdot b)^{-1},\gamma_+)]. \hfill\cr
\endmatrix$$
\noindent Let $\xi_\pm:I\rightarrow\partial\Omega$ be such that $\xi_\pm$ avoids $\alpha_\mp$ and that:
$$
\xi_\pm(0) = \alpha_\pm, \qquad \xi_\pm(1) = \gamma_\pm.
$$
\noindent Since $\partial\Omega$ is invariant under the action of $\alpha$, the curves $\xi_\pm\cdot\alpha$
lie in $\partial\Omega$ and join $\alpha_\pm$ to $\gamma_\pm\cdot\alpha$ whilst avoiding $\alpha_\mp$. Let 
$(\aleph_t^\pm)_{t\in I}$ be a continuous family of curves in $\partial\Omega$ such that:
\medskip
\myitem{(i)} for all $s$:
$$
\aleph_0^\pm(s) = \alpha_\pm,
$$
\noindent and,
\medskip
\myitem{(ii)} for all $t$:
$$\matrix
\aleph_t^\pm(0) \hfill&= \xi_\pm(t), \hfill\cr
\aleph_t^\pm(1) \hfill&= \xi_\pm(t)\cdot\alpha.\hfill\cr
\endmatrix$$
\noindent To be more precise, we deform $\aleph^-$ slightly towards the interior of $\partial\Omega$ and 
$\aleph^+$ slightly towards the exterior of $\partial\Omega$ so that they do not intersect each other. Heuristically, $\aleph_1^\pm$ is a curve which goes from $\gamma_\pm$ to $\gamma_\pm\cdot\alpha$ in the
clockwise direction. For all $t$ and for all $s$, we may suppose that the three points 
$(\aleph_t^-(0),a(0),\aleph_t^+(0))$ are distinct. Moreover, for all $t$:
$$
(\aleph_t^-(0),a(0),\aleph_t^+(0))\cdot\alpha = (\aleph_t^-(1),a(1),\aleph_t^+(1)).
$$
\noindent Consequently, the family $(\aleph_t^-,a,\aleph_t^+)_{t\in I}$ defines a homotopy between the
braids $(\alpha_-,a,\alpha_+)$ and $(\aleph_1^-,a,\aleph_1^+)$, respecting the fact that $\alpha$ sends 
the start point of each braid onto its end point. Since the set of liftings of $\alpha$ in $\widetilde{\psl}(2,\Bbb{C})$ is discrete, it follows by continuity that:
$$
[\Cal{S}(\aleph_1^-,a,\aleph_1^+)] = [\Cal{S}(\alpha_-,a,\alpha_+)].
$$
\noindent We define $(\mybeth_t^\pm)_{t\in I}$ in a similar fashion using $\beta$ instead of $\alpha$, and
we obtain:
$$\matrix
\hat{\gamma}_T^{-1}\hat{\alpha}_T\hat{\beta}_T\hfill&=[\Cal{S}(\gamma_-,(a\cdot b)^{-1},\gamma_+)
\Cal{S}(\alpha_-,a,\alpha_+)\Cal{S}(\beta_-,b,\beta_+)]\hfill\cr
&=[\Cal{S}(\gamma_-,(a\cdot b)^{-1},\gamma_+)\Cal{S}(\aleph_1^-,a,\aleph_1^+)
\Cal{S}(\mybeth_1^-,b,\mybeth_1^+)]\hfill\cr
&=[\Cal{S}(\gamma_-\cdot\aleph_1^-\cdot\mybeth_1^-,(a\cdot b)^{-1}\cdot a\cdot b,\gamma_+\cdot\aleph_1^+\cdot\mybeth_1^+)].\hfill\cr
\endmatrix$$
\noindent We thus define the closed braid $\Cal{T}$ by:
$$
\Cal{T} = (\gamma_-\cdot\aleph_1^-\cdot\mybeth_1^-,(a\cdot b)^{-1}\cdot a\cdot b,\gamma_+\cdot\aleph_1^+\cdot\mybeth_1^+).
$$
\noindent Since $\aleph_1^\pm$ and $\mybeth_1^\pm$ stay close to $\partial\Omega$ and 
$(a\cdot b)^{-1}\cdot a\cdot b$ lies in the interior of $\Omega$, there exists a homotopy between 
$(a\cdot b)^{-1}\cdot a\cdot b$ and the constant curve $p_0$ which stays away from $\aleph_1^\pm$ and 
$\mybeth_1^\pm$. Consequently:
$$\matrix
\Cal{T} \hfill&\sim (\gamma_-\cdot\aleph_1^-\cdot\mybeth_1^-,p_0, \gamma_+\cdot\aleph_1^+\cdot\mybeth_1^+) \hfill\cr
&\sim (\aleph_1^-\cdot\mybeth_1^-,p_0, \aleph_1^+\cdot\mybeth_1^+). \hfill\cr
\endmatrix$$
\noindent Heuristically, $\aleph_1^+\cdot\mybeth_1^+$ is a curve which turns once about $\partial\Omega$ in
the clockwise direction. Moreover, this curve lies to the exterior of $\partial\Omega$. There thus exists a
homotopy which stays away from $p_0$ and $\aleph_1^-\cdot\mybeth_1^-$ between this curve and the constant
curve $q_+$. Thus:
$$
\Cal{T} \sim (\aleph_1^-\cdot\mybeth_1^-,p_0, q_+).
$$
\noindent Finally, heuristically, $\aleph_1^-\cdot\mybeth_1^-$ is a simple closed curve which separates 
$p_0$ from $q_+$. Consequently:
$$\matrix
&\opWind_r(\Cal{T}) \hfill&= 1\hfill\cr
\Rightarrow\hfill&W_r(\hat{\gamma}_\Omega^{-1}\hat{\alpha}_\Omega\hat{\beta}_\Omega) \hfill&=1.\hfill\cr
\endmatrix$$
\noindent The result now follows.\qed
\medskip
\noindent We also obtain the converse to this result:
\proclaim{Proposition \nextprocno}
\noindent Let $\Gamma=\langle\alpha,\beta\rangle$ be a Schottky group. If $[a]_\alpha\in L_\alpha$, 
$[b]_\beta\in L_\beta$ and $[c]_\gamma\in L_\gamma$ are such that:
$$
(\opLift_\gamma[c]_\gamma)^{-1}(\opLift_\alpha[a]_\alpha)(\opLift_\beta[b]_\beta)\neq\opId,
$$
\noindent then there exists an invariant domain $\Omega$ of $\partial_\infty M$ adapted to the
generators $(\alpha,\beta)$ such that:
$$
[\Omega]_\alpha = [a]_\alpha, \qquad [\Omega]_\beta = [b]_\beta, \qquad [\Omega]_\gamma = [c]_\gamma.
$$
\endproclaim
\proclabel{InverseZipper}
\proof By Lemma \procref{HomeoEquivalence}, it suffices to prove this result in the case where 
$\partial_\infty M=\hat{\Bbb{C}}$ and $\Gamma$ is a Schottky subgroup of $\psl(2,\Bbb{C})$.
\medskip
\noindent Let $(C_\alpha^{\pm},C_\beta^\pm)$ be generating circles of $\Gamma$ with respect to the generators
$(\alpha,\beta)$. We define the set $X\subseteq\hat{\Bbb{C}}$ by:
$$
X = \hat{\Bbb{C}}\setminus\left(
\opInt(C_\alpha^+)\munion\opInt(C_\alpha^-)\munion\opInt(C_\beta^+)\munion\opInt(C_\beta^-)\right).
$$
\noindent We now define the curves $(a,b,c_1,c_2)$ to be non-intersecting, simple curves, lying in the
interior of $X$ except at their end points, such that:
$$\matrix
a(1)\hfill&=a(0)\cdot\alpha,\hfill&\qquad&b(1)\hfill&=b(0)\cdot\beta,\hfill\cr
c_1(0)\hfill&=c_2(1)\cdot\alpha,\hfill&\qquad&c_2(0)\hfill&=c_2(1)\cdot\beta.\hfill\cr 
\endmatrix$$
\noindent We refer to the quadruplet $(a,b,c_1,c_2)$ as {\emph generating curves\/} for $\Gamma$ with
respect to the generating circles $(C_\alpha^\pm,C_\beta^\pm)$. By taking the images of these four curves
under the actions of elements of $\gamma$, and by then adjoining $\opFix(\gamma)$, we obtain uniquely from
these four curves a Jordan curve which is invariant under the action of $\Gamma$ and which we denote by 
$\Gamma(a,b,c_1,c_2)$. The interior of $\Gamma(a,b,c_1,c_2)$ is an invariant domain of $\Gamma$ which is
adapted to the generators $(\alpha,\beta)$. Let us denote this domain by $\Omega(a,b,c_1,c_2)$. By definition:
$$
[a]_\alpha = [\Omega(a,b,c_1,c_2)]_\alpha,\qquad [b]_\beta = [\Omega(a,b,c_1,c_2)]_\beta.
$$
\noindent We define the curve $c$ by $c=c_2^{-1}(c_1\cdot\beta^{-1})^{-1}$, and we observe that:
$$
c(1) = c(0)\cdot\gamma.
$$
\noindent Consequently:
$$
[c]_\gamma = [\Omega(a,b,c_1,c_2)]_\gamma.
$$
\noindent By Proposition \procref{LemmaBraguette}:
$$
\opLift_\gamma([c]_\gamma)^{-1}\opLift_\beta([b]_\beta)\opLift_\alpha([a]_\alpha) \neq \opId.
$$
\noindent Let $T_1$ be the Dehn twist about $C_\alpha^+$ such that:
$$
T_1[a]_\alpha = [a]_\alpha + 1, \qquad T_1[b]_\beta = [b]_\beta, \qquad T_1[c]_\gamma = [c]_\gamma + 1.
$$
\noindent Let $T_2$ be the Dehn twist about $C_\beta^+$ such that:
$$
T_2[a]_\alpha = [a]_\alpha, \qquad T_2[b]_\beta = [b]_\beta + 1, \qquad T_2[c]_\gamma = [c]_\gamma + 1.
$$
\noindent Finally, let $T_3$ be a Dehn twist about a curve that separates $C_\alpha^\pm$ from $C_\beta^\pm$.
We may choose $T_3$ such that:
$$
T_3[a]_\alpha = [a]_\alpha, \qquad T_3[b]_\beta = [b]_\beta, \qquad T_3[c]_\gamma = [c]_\gamma + 2.
$$
\noindent Using combinations of these three Dehn twists, for any triple 
$(x,y,z)\in L_\alpha\times L_\beta\times L_\gamma$ satisfying:
$$
([a]_\alpha - x) + ([b]_\beta - y) + ([c]_\gamma - z) = 0 \text{ Mod }2,
$$
\noindent we can construct a quadruplet $(a',b',c_1',c_2')$ of generating curves of $\Gamma$ such that:
$$
[\Omega(a',b',c_1',c_2')]_\alpha = x,\qquad 
[\Omega(a',b',c_1',c_2')]_\beta = y,\qquad 
[\Omega(a',b',c_1',c_2')]_\gamma = z.
$$
\noindent By Proposition \procref{ConditionForEquivalenceOfLifts}, the triple $(x,y,z)$ satisfies this condition precisely when:
$$
\opLift_\gamma(z)^{-1}\opLift_\alpha(x)\opLift_\beta(y) \neq \opId.
$$
\noindent The result now follows.\qed\goodbreak
\newhead{Constructing The Local Homeomorphism}
\newsubhead{Introduction}
\noindent Throughout this section, $(M,Q)$ will be a pointed, compact, three dimensional manifold of
strictly negative sectional curvature, $(\tilde{M},\tilde{Q})$ its universal cover, and 
$\Gamma=\langle\alpha,\beta\rangle$ a Schottky subgroup of $\pi_1(M,Q)$.
\medskip
\noindent In this section, we will prove the main results of this paper. First we have:
\proclaim{Theorem \procref{PresentationChIVExistence}}
\noindent Suppose that $(M,Q)$ is a pointed, compact manifold of strictly negative sectional curvature. 
Let $(\Sigma,P)$ be a pointed, compact surface of hyperbolic type (i.e. of genus at least two). Let
$\theta:\pi_1(\Sigma,P)\rightarrow\pi_1(M,Q)$ be a homomorphism. Suppose that $\theta$ is non-elementary and
may be lifted to a homomorphism $\hat{\theta}$ of $\pi_1(\Sigma,P)$ into the group 
$\opTildeHomeo_0(\partial_\infty\tilde{M})$. Then there exists an equivariant Plateau problem for $\theta$.
\endproclaim
\noindent We then prove:
\proclaim{Theorem \procref{PresentationChIVConvexRealisation}}
\noindent If $\theta$ is non-elementary and lifts, then there exists a convex immersion 
$i:\Sigma\rightarrow M$ such that:
$$
\theta = i_*.
$$
\endproclaim
\noindent In the second part of this section, we define the notion of a $\pi_1(M,Q)$ structure, and
prove an existence result for such structures over bound, marked trousers. In the third part of this
section, we state results which summarise the glueing technique using by Gallo, Kapovich and Marden in
section $8$ of \cite{GallKapMard}. In the fourth section we provide a proof of Theorem
\procref{PresentationChIVExistence}. Finally, in the fifth section we show how to obtain convex solutions to
the Plateau problem in Hadamard manifolds, and this permits us to prove Theorem 
\procref{PresentationChIVConvexRealisation}.
\newsubhead{Constructing the Solution Over Trousers}
\noindent Let $(\Sigma,\theta,\beta)$ be a bound marked surface. Let $\tilde{\Sigma}$ be the universal
cover of $\Sigma$. We define a $\pi_1(M,Q)$ {\emph structure\/} over $(\Sigma,\theta,\beta)$ to be a 
local homeomorphism $\varphi:\tilde{\Sigma}\rightarrow\partial_\infty\tilde{M}$ such that:
\medskip
\myitem{(i)} $\varphi$ is equivariant under the action of $\theta$, and
\medskip
\myitem{(ii)} for every boundary component $C$ of $\tilde{\Sigma}$, the restriction of $\varphi$ to $C$
is a homeomorphism onto its image. In otherwords $\varphi(C)$ is a non self-intersecting curve.
\medskip
\noindent Let $P_0$ be the base point of $\Sigma$, and let $\tilde{P}_0$ be the corresponding base point
in the universal cover $\tilde{\Sigma}$ of $\Sigma$. Let $(C,P)$ be a pointed boundary component of $\Sigma$.
Let $\eta$ be an element of $\pi_1(\Sigma,P_0,P)$, and let $\tilde{\eta}$ be the lift of $\eta$ such that:
$$
\tilde{\eta}(0) = \tilde{P}_0.
$$
\noindent Let $\hat{P}_\eta$ be the endpoint of $\tilde{\eta}$. $\hat{P}_\eta$ is thus a lift of $P$.
Viewing $C$ as a parametrised simple closed curve in $\Sigma$ leaving $P$, we define $\hat{C}_\eta$ to be
the lift of $C$ starting from $\hat{P}_\eta$. The curve $C_\eta$ is thus a segment of one of the boundary
components of $\tilde{\Sigma}$. Moreover, by definition:
$$
\hat{C}_\eta(0)\cdot(\eta^{-1}C\eta) = \hat{C}_\eta(1).
$$
\noindent Since $\varphi\circ\hat{C}_\eta$ is non self-intersecting, it must avoid the fixed points of 
$\theta(\eta^{-1}C\eta)$ in $\partial_\infty\tilde{M}$. Consequently $\varphi\circ\hat{C}_\eta$ projects
down to a closed curve in $\Bbb{T}_{\theta(\eta^{-1}C\eta)}$. We will denote this element by 
$[C]_{\eta,\varphi}$.
\medskip
\noindent Let $\gamma$ and $\xi$ be two elements of $\pi_1(M,Q_0)$. Let $\xi$ be another element. The mapping
$\xi$ sends $\gamma_\pm$ onto the fixed points $\delta_\pm$ of $\delta=\xi^{-1}\gamma\xi$. It follows
that $\xi$ defines a homeomorphism from $\Bbb{T}_\gamma$ to $\Bbb{T}_\delta$, which we will also denote
by $\xi$. Moreover for any $\alpha\in\pi_1(\Sigma)$:
$$\matrix
&\hat{C}_\eta\cdot\alpha\hfill&= \hat{C}_{\eta\alpha}\hfill\cr
\Rightarrow\hfill&\varphi(\hat{C}_\eta\cdot\alpha) \hfill&=\varphi(\hat{C}_{\eta\alpha})\hfill\cr
\Rightarrow\hfill&(\varphi\circ\hat{C}_\eta)\cdot\theta(\alpha)\hfill&=(\varphi\circ\hat{C}_{\eta\alpha})\hfill\cr
\Rightarrow\hfill&[C]_{\eta,\varphi}\cdot\theta(\alpha)\hfill&= [C]_{\eta\alpha,\varphi}\hfill.\cr
\endmatrix$$
\noindent We may thus define the element $\beta_\varphi(C)$ lying in $\Bbb{T}_{\beta(C)}$ such that, for all
$\eta\in\pi_1(\Sigma,P_0,P)$:
$$
\beta_\varphi(C) = [C]_{\eta,\varphi}\cdot\beta(\eta)^{-1}.
$$
%
%
%
%
%
%
%
\noindent We are now in a position to construct $\pi_1(M,Q_0)$ structures over bound, marked trousers:
\proclaim{Proposition \nextprocno}
\noindent Let $(M,Q_0)$ be a pointed three dimensional Hadamard manifold of strictly negative sectional
curvature. Let $(T,\theta,\beta)$ be a bound marked trouser with holonomy in $\pi_1(M,Q_0)$. Suppose that
the image of $\theta$ is a Schottky group.
\medskip
\noindent Let $P_0$ be the base point of $T$ and let $(C_j,Q_j)_{1\leqslant j\leqslant 3}$ be the three
oriented boundary components of $T$. Let $(\xi_j)_{1\leqslant j\leqslant 3}$ be a sash of $T$
with respect to the points $(Q_j)_{1\leqslant j\leqslant 3}$ such that:
$$
\xi_3^{-1}C_3^{-1}\xi_3=(\xi_1^{-1}C_1^{-1}\xi_1)(\xi_2^{-1}C_2^{-1}\xi_2).
$$
\noindent Let $(x_1,x_2,x_3)\in L_{\beta(C_1)}\times L_{\beta(C_2)}\times L_{\beta(C_3)}$ be a triplet such
that:
$$
\opLift_{\theta(\xi_3^{-1}C_3\xi_3)}(x_3\cdot\beta(\xi_3))^{-1}
\opLift_{\theta(\xi_1^{-1}C_1\xi_1)}(x_1\cdot\beta(\xi_1))
\opLift_{\theta(\xi_2^{-1}C_2\xi_2)}(x_2\cdot\beta(\xi_2))\neq\opId.
$$
\noindent Then, there exists a $\pi_1(M,Q_0)$ structure, $\varphi$, over $(T,\theta,\beta)$ such that, for 
each $i$:
$$
\beta_\varphi(C_i) = x_i.
$$
\endproclaim
\proclabel{ExistenceForTrousers}
\proof Let $P_T$ be the base point of $T$. Let us denote $a=\xi_1^{-1}C_1\xi_1$ and $b=\xi_2^{-1}C_2\xi_2$.
Let $\Gamma=\opIm(\theta)$. Denote $\alpha=\theta(a)$ and $\beta=\theta(b)$. By Lemma 
\procref{PrincipalObstructionToJoin}, there exists an invariant domain $\Omega$ of $\Gamma$
adapted to the generators $(\alpha,\beta)$ such that:
$$
[\Omega]_\alpha = x_1\cdot\beta(\xi_1),\qquad
[\Omega]_\beta = x_2\cdot\beta(\xi_2),\qquad
[\Omega]_\gamma = x_3\cdot\beta(\xi_3).
$$
\noindent Let $P_0$ be a point in the interior of $\hat{\Omega}$. Let $\pi:\hat{\Omega}\rightarrow\hat{\Omega}/\Gamma$ be the canonical projection. Let $\phi:\pi_1(\hat{\Omega}/\Gamma,\pi(P_0))\rightarrow\Gamma$ be the unique isomorphism defined such that, for all $\gamma\in\Gamma$ and for any curve $c$ joining $P_0$ to $\gamma(P_0)$:
$$
\phi([\pi\circ c]) = \gamma.
$$
\noindent In particular, since $\Omega$ is adapted to the generators $(\alpha,\beta)$, we see that there
exist boundary components $a'$ and $b'$ of $\hat{\Omega}/\Gamma$ and a sash $(c_{a'},c_{b'})$ of
$(\hat{\Omega}/\Gamma,\pi(P_0))$ with respect to the base points of $a'$ and $b'$ respectively such that:
$$
\phi(c_{a'}^{-1}a'c_{a'}) = \alpha,\qquad\phi(c_{b'}^{-1}b'c_{b'}) = \beta.
$$
\noindent There exists a homeomorphism, $\varphi:(T,P_T)\rightarrow(\hat{\Omega}/\Gamma,\pi(P_0))$ such that:
$$
\varphi_*a\sim c_{a'}^{-1}a'c_{a'},\qquad \varphi_*b\sim c_{b'}^{-1}b'c_{b'}.
$$
\noindent Consequently:
$$
\theta = \phi\circ\varphi_*.
$$
\noindent Let $(\tilde{T},\tilde{P}_T)$ denote the pointed universal cover of $(T,P_T)$. The mapping
$\varphi$ lifts to a unique pointed homeomorphism, $\tilde{\varphi}$, from $(\tilde{T},\tilde{P}_T)$ to 
$(\hat{\Omega},P_0)$ which is equivariant under the action of $\theta$. The curve 
$\tilde{\varphi}\circ C_1$ lies entirely in $\hat{\Omega}$. Thus, by definition of $[\Omega]_\alpha$:
$$\matrix
&[C_1]_{\xi_1,\tilde{\varphi}}\hfill
&=[\Omega]_\alpha\hfill
&=x_1\cdot\beta(\xi_1)\hfill\cr
\Rightarrow\hfill
&\beta_{\tilde{\varphi}}(C_1)\hfill
&=[C_1]_{\xi_1,\tilde{\varphi}}\cdot\beta(\xi_1)^{-1}\hfill
&=x_1.\hfill\cr
\endmatrix$$
\noindent Similarly, $\beta_{\tilde{\varphi}}(C_2)=x_2$ and $\beta_{\tilde{\varphi}}(C_3)=x_3$. The
mapping $\tilde{\varphi}$ is consequently the desired $\pi_1(M,Q_0)$ structure, and the result now 
follows.\qed
\newsubhead{Joining the Trousers}
\noindent Let $(\Sigma_i,\theta_i,\beta_i)_{i\in\left\{1,2\right\}}$ be two bound marked surfaces with
holonomy in $\pi_1(M,Q)$. For each $i$, let $(C_i,Q_i)$ be a pointed boundary component of $\Sigma_i$.
Let $\psi:(C_1,Q_1)\rightarrow (C_2^{-1},Q_2)$ be a homeomorphism and suppose that 
$(\Sigma_1,\theta_1,\beta_1)$ may be joined to $(\Sigma_2,\theta_2,\beta_2)$ along $\psi$. We recall that
this means that $\beta_1(C_1)^{-1}=\beta_2(C_2)$ and consequently that:
$$
\Bbb{T}_{\beta_1(C_1)^{-1}} = \Bbb{T}_{\beta_2(C_2)}.
$$
\noindent We observe that:
$$
L_{\gamma^{-1}} = -L_\gamma.
$$
\noindent For each $i$, let $\varphi_i$ be a $\pi_1(M,Q)$ structure over $(\Sigma_i,\theta_i,\beta_i)$. We
say that $\varphi_1$ may be {\emph joined\/} to $\varphi_2$ along $\psi$ if and only if:
$$
-\beta_{1,\varphi_1}(C_1) = \beta_{2,\varphi_2}(C_2).
$$
\noindent The glueing procedure described by Gallo, Kapovich and Marden in section $8$ of 
\cite{GallKapMard} yields the following result:
\proclaim{Proposition \nextprocno}
\noindent Let $(M,Q)$ be a compact, pointed, three dimensional manifold of strictly negative sectional 
curvature.
\medskip
\noindent Let $(\Sigma_i,\theta_i,\beta_i)_{i\in\left\{1,2\right\}}$ be two bound marked surfaces with
holonomy in $\pi_1(M,Q)$. For each $i$, let $(C_i,Q_i)$ be pointed boundary components of $\Sigma_i$ and let
$\psi:(C_1,Q_1)\rightarrow(C_2^{-1},Q_2)$ be a homeomorphism such that $(\Sigma_1,\theta_1,\beta_1)$ may be
joined to $(\Sigma_2,\theta_2,\beta_2)$ along $\psi$.
\medskip
\noindent Suppose that there exists, for each $i$, a $\pi_1(M,Q)$ structure $\varphi_i$ over
$(\Sigma_i,\theta_i,\beta_i)$. If moreover, $\varphi_1$ may be joined to $\varphi_2$ along $\psi$, then 
there exists a $\pi_1(M,Q)$ structure $\varphi$ over
$(\Sigma_1,\theta_1,\beta_1)\cup(\Sigma_2,\theta_2,\beta_2)$ such that, after identifying objects
in each of $\Sigma_1$ and $\Sigma_2$ with the corresponding objects in $\Sigma_1\cup\Sigma_2$, for 
each $i$ and for every pointed boundary component $(C',Q')$ of $\Sigma_i$:
$$
(\beta_1\munion\beta_2)_\varphi(C') = \beta_{i,\varphi_i}(C').
$$
\endproclaim
\proclabel{JoiningDistinctTrousers}
\noindent Let $(\Sigma,\theta,\beta)$ be a bound marked surface with holonomy in $\pi_1(M,Q)$. Let
$(C_i,Q_i)_{i\in\left\{1,2\right\}}$ be pointed boundary components of $\Sigma$. Let
$\psi:(C_1,Q_1)\rightarrow(C_2^{-1},Q_2)$ be an orientation reversing homeomorphism and suppose that $(\Sigma,\theta,\beta)$ may be joined to itself along $\psi$. Let $\varphi$ be a $\pi_1(M,Q)$ structure over $(\Sigma,\theta,\beta)$. We say that $\varphi$ may be {\emph joined to itself\/} along $\psi$ if and only if:
$$
-\beta_\varphi(C_1) = \beta_\varphi(C_2).
$$
\noindent Once again, the glueing procedure described by Gallo, Kapovich and Marden in section $8$
of \cite{GallKapMard} permits us to obtain the following analogue of Lemma \procref{JoiningDistinctTrousers}:
\proclaim{Proposition \nextprocno}
\noindent Let $(M,Q)$ be a compact, pointed, three dimensional manifold of strictly negative sectional
curvature.
\medskip
\noindent Let $(\Sigma,\theta,\beta)$ be a bound, marked surface with holonomy in $\pi_1(M,Q)$. Let
$(C_i,Q_i)_{i\in\left\{1,2\right\}}$ be pointed boundary components of $\Sigma$. Let 
$\psi:(C_1,Q_1)\rightarrow(C_2^{-1},Q_2)$ be a homeomorphism such that $(\Sigma,\theta,\beta)$ may be 
joined to itself along $\psi$.
\medskip
\noindent Suppose that there exists, for each $i$, a $\pi_1(M,Q)$ structure, $\varphi$, over 
$(\Sigma,\theta,\beta)$. If, moreover, $\varphi$ may be joined to itself along $\psi$, then there exists a
$\pi_1(M,Q)$ structure $\varphi'$ over $\mathcircle(\Sigma,\theta,\beta)$ such that, after identifying
objects in $\Sigma$ with the corresponding objects in $\mathcircle\Sigma$, for every pointed
boundary component $(C',Q')$ of $\Sigma$:
$$
(\mathcircle\beta)_{\varphi'}(C') = (\mathcircle\beta)_\varphi(C).
$$
\endproclaim
\proclabel{JoiningSameTrouser}
\newsubhead{The Construction of a Local Homeomorphism}
\noindent We are now in a position to prove Theorem \procref{PresentationChIVExistence}:
\medskip
{\bf\noindent Proof of Theorem \procref{PresentationChIVExistence}:\ }By Lemma 
\procref{SchottkyTrouserDecomposition}, there exist $2g-2$ bound, marked trousers
$(T_i,\theta_i,\beta_i)_{1\leqslant i\leqslant 2g-2}$ with holonomy in $\pi_1(M,Q)$ such that the image
of every $\theta_i$ is a Schottky group, and:
$$
(\Sigma,\theta)\cong \mathcircle_{i=1}^g(\munion_{i=1}^{2g-2}(T_i,\theta_i,\beta_i)).
$$
\noindent Let $\pi:\opTildeHomeo_0(\partial_\infty\tilde{M})\rightarrow\opHomeo_0(\partial_\infty\tilde{M})$
be the canonical projection. For every $i$, using the lifting, $\hat{\theta}$, of $\theta$, we may construct 
liftings $\hat{\theta}_i$ and $\hat{\beta}_i$ of $\theta_i$ and $\beta_i$ such that:
\medskip
\myitem{(i)} for all $i$, $\pi\circ\hat{\theta}_i=\theta_i$ and $\pi\circ\hat{\beta}_i=\beta_i$, and
\medskip
\myitem{(ii)} the $(T_i,\hat{\theta}_i,\hat{\beta}_i)_{1\leqslant i\leqslant 2g-2}$ may be joined along the
same edges as the $(T_i,\theta_i,\beta_1)_{1\leqslant i\leqslant 2g-2}$.
\medskip
\noindent For each $i$, we label the boundary components of $T_i$ with distinct numbers from $1$ to $3$. For all $1\leqslant i\leqslant 2g-2$,$1\leqslant j\leqslant 3$, we denote by $(C_{i,j},Q_{i,j})$ the $j$'th
pointed boundary component of $T_i$. Let $\alpha$ be the involution (having no fixed points) of 
$A=\left\{1,...,2g-2\right\}\times\left\{1,2,3\right\}$ such that, for all $(i,j)\in A$, the boundary component $(C_{i,j},Q_{i,j})$ is joined to $(C_{\alpha(i,j)},Q_{\alpha(i,j)})$ in the above decomposition. For all $i$,$j$, we define $E_{i,j}\in L_{\beta(C_{i,j})}$ such that:
\medskip
\myitem{(i)} for all $(i,j)$, $\opLift_{\beta_i(C_{i,j})}E_{i,j}=\neq\hat{\beta}_i(C_{i,j})$, and
\medskip
\myitem{(ii)} for all $(i,j)$, $E_{i,j}=-E_{\alpha(i,j)}$.
\medskip 
\noindent Choose $1\leqslant i\leqslant 2g-2$. Let $\Cal{P}_i$ be the base point of $T_i$. Let
$(\xi_j)_{1\leqslant j\leqslant 3}$ be a binding sash of $T_i$ with respect to the points
$(Q_{i,j})_{1\leqslant j\leqslant 3}$ such that:
$$
\xi_3^{-1}C_{i,3}\xi_3 = (\xi_1^{-1}C_{i,1}\xi_1)(\xi_2^{-1}C_{i,2}\xi_2).
$$
\noindent For each $j$, let us denote $\gamma_j=\xi_j^{-1}C_{i,j}\xi_j$ and $a_j=E_{i,j}\cdot\beta(\xi_j)$.
Let $\Delta\in\opTildeHomeo_0(\partial_\infty\tilde{M})$ be the element of $\pi^{-1}(\opId)$ that is different to the identity. We recall that $\Delta$ commutes with every element of 
$\opTildeHomeo(\partial_\infty\tilde{M})$. For each $j$, we have:
$$\matrix
\opLift_{\theta(\gamma_j)}(a_j)\hfill
&=\beta(\xi_j)^{-1}[\opLift_{\beta(\xi_j)\theta(\gamma_j)\beta(\xi_j)^{-1}}(E_{i,j})]\beta(\xi_j)\hfill\cr
&=\hat{\beta}(\xi_j)^{-1}\opLift_{\beta(C_{i,j})}(E_{i,j})\hat{\beta}(\xi_j)\hfill\cr
&=\hat{\beta}(\xi_j)^{-1}\hat{\beta}(C_{i,j})\Delta\hat{\beta}(\xi_j)\hfill\cr
&=\hat{\theta}(\gamma_j)\Delta.\hfill\cr
\endmatrix$$
\noindent Thus:
$$\matrix
\opLift_{\theta(\gamma_3)}(a_3)^{-1}
\opLift_{\theta(\gamma_1)}(a_1)\opLift_{\theta(\gamma_2)}(a_2)\hfill
&=\hat{\theta}(\gamma_3)^{-1}\hat{\theta}(\gamma_1)\hat{\theta}(\gamma_2)\Delta^3\hfill\cr
&=\hat{\theta}(\gamma_3^{-1}\gamma_1\gamma_2)\Delta^3\hfill\cr
&=\Delta^3\hfill\cr
&=\Delta.\hfill\cr
\endmatrix$$
\noindent Thus, by Proposition \procref{ExistenceForTrousers}, there exists a $\pi_1(M,Q)$ structure
$\varphi_i$ over $(T_i,\theta_i,\beta_i)$ such that, for each $j$:
$$
\beta_{i,\varphi_i}(C_{i,j}) = E_{i,j}.
$$
\noindent By definition of $(E_{i,j})_{(i,j)\in A}$, these $\pi_1(M,Q)$ structures may all be joined to each
other, and the existence of a $\pi_1(M,Q)$ structure over $(\Sigma,\theta)$ now follows using a process of
induction and Propositions \procref{JoiningDistinctTrousers} and \procref{JoiningSameTrouser}. The result now follows.\qed
\newsubhead{Convex Subsets of Hadamard Manifolds}
\noindent We now require the following results concerning convex sets:
\proclaim{Proposition \nextprocno}
\noindent Let $p$ be a point in $\partial_\infty\tilde{M}$ and let $\Omega$ be a neighbourhood of $p$ in 
$\partial_\infty\tilde{M}$. There exists a complete convex subset $X$ of $\tilde{M}$ 
such that:
$$
p\notin\partial_\infty X,\qquad\Omega^c\subseteq\partial_\infty X.
$$
\endproclaim
\proclabel{ExistenceOfConvexSets}
\proof Let $K\geqslant k>0$ be such that the sectional curvature of $\tilde{M}$ is contained in
$[-K,-k]$. For all $q$ in $\tilde{M}$, for $r\in\partial_\infty\tilde{M}$ and for $\theta\in(0,\pi)$,
we define $\Omega_\theta(r;q)$ as in section \subheadref{ConeControl}. By Theorem $3.1$ of \cite{And}, there exists $\pi/2>\psi>0$ such that, for all $q\in\tilde{M}$ and for all $r\in\partial_\infty\tilde{M}$, there exists a complete convex subset $X$ of $\tilde{M}$ such that:
$$
\partial_\infty\Omega_\psi(r;q)\subseteq \partial_\infty X.
$$
\noindent Let $q_0$ be any point in $\tilde{M}$. Let $\gamma$ be the unique geodesic running from $q_0$
to $p$, normalised such that $\gamma(0)=q_0$ and $\gamma(\infty)=p$. Let $\delta\in\Bbb{R}^+$ be such that:
$$
\partial_\infty\Omega_\delta(p;q_0)\subseteq\Omega.
$$
\noindent Let us define $B\in\Bbb{R}^+$ such that:
$$
\Delta_k(B,\pi-\delta)<\psi.
$$
\noindent Let $q_1=\gamma(B)$. Let $X$ be a complete convex subset of $\tilde{M}$ such that:
$$
\partial_\infty\Omega_\psi(q_1;\gamma(-\infty))\subseteq \partial_\infty X.
$$
\noindent Thus, by Lemma \procref{LemmaConeControl}:
$$
\Omega^c\subseteq\partial_\infty\Omega_{\pi-\delta}(p;\gamma(-\infty))
\subseteq\partial_\infty\Omega_\psi(q_1;\gamma(-\infty))\subseteq\partial_\infty X.
$$
\noindent The result now follows.\qed
\medskip
\noindent This result now permits us to show:
\proclaim{Lemma \nextprocno}
\noindent Let $D$ be a closed subset of $\partial_\infty\tilde{M}$. There exists a complete strictly convex
subset $X$ with $C^1$ boundary in $\tilde{M}$ such that:
$$
D = \partial_\infty X.
$$
\endproclaim
\proclabel{SolutionToSimpleConvexPlateauProblem}
\proof For all $p\in D^c$, by Proposition \procref{ExistenceOfConvexSets}, there exists a complete convex
subset $X_p$ of $\partial_\infty\tilde{M}$ such that $p\notin\partial_\infty X_p$ and 
$D\subseteq\partial_\infty X_p$. We define $\hat{X}$ by:
$$
\hat{X} = \minter_{p\in C^c} X_p.
$$
\noindent $\hat{X}$ is a complete convex subset of $\tilde{M}$ and $\partial_\infty\hat{X}=D$. In particular,
$\hat{X}$ is non-empty. For $\epsilon\in\Bbb{R}^+$, let $\Sigma_\epsilon$ be the surface obtained by 
moving $\partial\hat{X}$ a distance $\epsilon$ normally along geodesics. $\Sigma_\epsilon$ is strictly
convex, is of type $C^1$ and bounds a complete convex subset $X$ of $\tilde{M}$. $X$ is thus the desired subset, and the result now follows.\qed
\medskip
\noindent We are now in a position to prove Theorem \procref{PresentationChIVConvexRealisation}:
\medskip
{\bf\noindent Proof of Theorem \procref{PresentationChIVConvexRealisation}:\ }We give an elementary account of a technique developed by Labourie in \cite{LabA}. Since $\Sigma$ is compact, and since $\varphi$ is $\theta$-equivariant, there exists a finite open cover 
$\Cal{A}=(\Omega_i)_{1\leqslant i\leqslant k}$ of $\Sigma$ such that, if
$\tilde{\Cal{A}}=(\Omega_{i,\gamma})_{1\leqslant i\leqslant k,\gamma\in\pi_1(\Sigma,P)}$ is the lifting of 
$\Cal{A}$ to $\tilde{\Sigma}$, then, for every $\Omega\in\tilde{\Cal{A}}$:
\medskip
\myitem{(i)} the closure of $\Omega$ in $\tilde{\Sigma}$ is homeomorphic to a closed disc, and
\medskip
\myitem{(ii)} for every $\Omega\in\tilde{\Cal{A}}$, the restriction of $\varphi$ to $\overline{\Omega}$ is a homeomorphism onto its image.
\medskip
\noindent By Lemma \procref{SolutionToSimpleConvexPlateauProblem}, for every $i$, we may find a family 
$(X_{i,\gamma})_{\gamma\in\pi_1(\Sigma,P)}$ of complete strictly convex subsets with $C^1$ 
boundary in $\tilde{M}$ such that:
\medskip
\myitem{(i)} $(X_{i,\gamma})_{\gamma\in\pi_1(\Sigma,P)}$ is $\theta$-equivariant. In other words, for all
$\gamma,\delta\in\pi_1(\Sigma,P)$:
$$
X_{i,\gamma\delta} = X_{i,\gamma}\cdot\theta(\delta),
$$
\noindent and
\medskip
\myitem{(ii)} for all $\gamma\in\pi_1(\Sigma,P)$:
$$
\partial_\infty X_{1,\gamma} = \varphi(\Omega_{i,\gamma})^c.
$$
\noindent We may assume moreover that the boundaries of these sets intersect transversally in the region that
will be of interest to us. Using these convex sets, we construct a family 
$(P_{i,\gamma})_{1\leqslant i\leqslant k,\gamma\in\pi_1(\Sigma,P)}$ of polygonal subsets of
$\tilde{\Sigma}$ and a family $(\psi_{i,\gamma})_{1\leqslant i\leqslant k,\gamma\in\pi_1(\Sigma,P)}$ of
immersions such that:
\medskip
\myitem{(i)} $(P_{i,\gamma})_{1\leqslant i\leqslant k,\gamma\in\pi_1(\Sigma,P)}$ provides a polygonal 
decomposition of $\tilde{\Sigma}$. In otherwords:
$$
\tilde{\Sigma} = 
\munion_{\matrix
{\scriptstyle 1\leqslant i\leqslant k}\cr
{\scriptstyle\gamma\in\pi_{(\Sigma,P)}}\cr
\endmatrix}P_{i,\gamma},
$$
\noindent and two polygons in this family only intersect if they share a common boundary component;
\medskip
\myitem{(ii)} the family $(P_{i,\gamma})_{1\leqslant i\leqslant k,\gamma\in\pi_1(\Sigma,P)}$ is equivariant.
In otherwords, for each $i$, and for all $\gamma,\delta\in\pi_1(\Sigma,P)$:
$$
P_{i,\gamma\delta} = P_{i,\gamma}\cdot\delta;
$$
\myitem{(iii)} for each $i$ and for all $\gamma$, $\psi_{i,\gamma}$ is a homeomorphism
from $P_{i,\gamma}$ onto a subset of the boundary of $X_{i,\gamma}$ in $\tilde{M}$, and 
\medskip
\myitem{(iv)} the family $(\psi_{i,\gamma})_{1\leqslant i\leqslant k,\gamma\in\pi_1(\Sigma,P)}$ is 
$\theta$-equivariant. In otherwords, for each $i$, and for all $\gamma,\delta\in\pi_1(\Sigma,P)$, and 
for all $p\in P_{i,\gamma}$:
$$
\psi_{i,\gamma\delta}(p\cdot\delta) = \psi_{i,\gamma}(p)\cdot\theta(\delta).
$$
\noindent Joining together the elements of 
$(\psi_{i,\gamma})_{1\leqslant i\leqslant k,\gamma\in\pi_1(\Sigma,P)}$, we obtain a $\theta$-equivariant
locally strictly convex immersion $\hat{\psi}$ of $\tilde{\Sigma}$ into $\tilde{M}$ such that if
$\overrightarrow{n}$ is the Gauss-Minkowski mapping sending $U\tilde{M}$ into $\partial_\infty\tilde{M}$,
then:
$$
\overrightarrow{n}\circ\hat{\psi} = \varphi.
$$
\noindent Taking quotients, we obtain a locally strictly convex immersion $\psi$ of $\Sigma$ into $M$ which realises $\theta$. Finally, since $\Sigma$ is compact, by deforming $\psi$ slightly, we may suppose that
it is also smooth and the result now follows.\qed
\global\headno=0
\inappendicestrue
\newhead{Homeomorphism Equivalence of Schottky Groups}
\newsubhead{Introduction}
\noindent Throughout this appendix, $(M,Q)$ will be a pointed three dimensional Hadamard manifold
and $\Gamma=\langle\alpha,\beta\rangle$ a Schottky subgroup of $\opIsom(M,Q)$. 
\headlabel{MakingHomeomorphism}
\medskip
\noindent In this appendix, we provide a proof of the equivalence up to homeomorphisms of the Schottky
groups that we will be using.
\medskip
\noindent We define $\opFix(\Gamma)$, the {\emph fixed point set\/} of $\Gamma$ by:
$$
\opFix(\Gamma) = \overline{\munion_{\gamma\in\Gamma\setminus\left\{\opId\right\}}\opFix(\gamma)}.
$$
\noindent We define a {\emph reduced word\/} over $(\alpha,\beta)$ to be a sequence 
$\pmb{\gamma}=(\gamma_k)_{1\leqslant k\leqslant n}$ of elements of 
$\left\{\alpha^{\pm 1},\beta^{\pm 1}\right\}$ which does not contain any of the subwords $\alpha\alpha^{-1}$,
$\alpha^{-1}\alpha$, $\beta\beta^{-1}$ or $\beta^{-1}\beta$. In otherwords, a reduced word is the shortest
length word expressing the corresponding element of $\Gamma$. Let $W_{\alpha,\beta}$ denote the set of all
reduced words over $(\alpha,\beta)$ of finite length. Let $W^n_{\alpha,\beta}$ denote the subset of 
$W_{\alpha,\beta}$ consisting of words of length $n$. Let $W^\infty_{\alpha,\beta}$ denote the set of all
reduced words over $(\alpha,\beta)$ of infinite length. For all $m\geqslant n$, we define the truncation map
$T_n:W^m_{\alpha,\beta}\rightarrow W^n_{\alpha,\beta}$ to be the map which sends a word of length $m$ to the
word given by its first $n$ letters. This definition is trivially also valid when $m=\infty$. We give
$W_{\alpha,\beta}$ the discrete topology, and we then give $W^\infty_{\alpha,\beta}$ the coarsest 
topology with respect to which every $T_n$ is continuous. In otherwords, a basis of open sets of 
$W^\infty_{\alpha,\beta}$ is given by:
$$
\Cal{B} = \munion_{n\in\Bbb{N}}\left\{T_n^{-1}(\left\{\pmb{\gamma}\right\})\text{ s.t. }
\pmb{\gamma}\text{ is of length }n\right\}.
$$
\noindent This topology is trivially Hausdorff. Moreover, since the alphabet is finite, it is not difficult
to show that $W^\infty_{\alpha,\beta}$ is compact with respect to this topology.
\medskip
\noindent The first result of this appendix is the construction of a canonical homeomorphism from 
$W^\infty_{\alpha,\beta}$ to $\opFix(\Gamma)$. In order to explicitely describe the homeomorphism, we
are required to construct a few more objects. Let $(C_\alpha^\pm,C_\beta^\pm)$ be generating circles for 
$\Gamma$ in $\partial_\infty\tilde{M}$ with respect to the generators $(\alpha,\beta)$. We define the set $S_1$ by $S_1=\left\{C_\alpha^\pm,C_\beta^\pm\right\}$ and we define $S_n$ for $n\geqslant 2$ inductively by:
$$
S_n = \left(\munion_{\gamma\in\left\{\alpha^{\pm 1},\beta^{\pm 1}\right\}} \gamma S_{n-1}\right)
\setminus\munion_{k=0}^{n-1}S_{n-1}.
$$
\noindent For all $n\geqslant 2$, and for all $C\in S_n$, there exists a unique $C'\in S_1$ such that
$C\subseteq\opInt(C')$. We thus orient $C$ such that $\opInt(C)\subseteq\opInt(C')$. We observe that, for
all $n$, and for any distinct $C,C'\in S_n$, the interiors of $C$ and $C'$ are also distinct.
\medskip
\noindent We define the mapping $D_1:W^1_{\alpha,\beta}\rightarrow S_1$ by:
$$
D_1(\alpha^{\pm 1}) = C_\alpha^{\pm},\qquad D_1(\beta^{\pm 1}) = C_\beta^{\pm}.
$$
\noindent We define $D_n$ for $n\geqslant 2$ such that, for all 
$\pmb{\gamma}=(\gamma_k)_{1\leqslant k\leqslant n}$ in $W^n_{\alpha,\beta}$:
$$
D_n(\pmb{\gamma}) = D_1(\gamma_n)\cdot \gamma_{n-1}\cdot...\cdot\gamma_1.
$$
\noindent We then define $D$ over $W_{\alpha,\beta}$ such that it restricts to $D_n$ over each 
$W^n_{\alpha,\beta}$. We observe that, for all $n$, $D_n$ defines a bijection between $W^n_{\alpha,\beta}$
and $S_n$. We now obtain the following result:
\proclaim{Lemma \nextprocno}
\noindent Suppose that the sectional curvature of $M$ is bounded above by $-k<0$.
For all $\pmb{\gamma}=(\gamma_n)_{n\in\Bbb{N}}$ in $W^\infty_{\alpha,\beta}$, there exists a
unique point $\Cal{P}(\pmb{\gamma})\in\opFix(\Gamma)$ such that the sequence 
$(\opInt(D(T_n(\pmb{\gamma}))))_\ninn$ converges towards $\left\{\Cal{P}(\pmb{\gamma})\right\}$ in the
Hausdorff topology. Moreover, $\Cal{P}(\pmb{\gamma})$ is contained in $\opInt(D(T_n(\pmb{\gamma})))$
for every $n$, and $\Cal{P}$ defines a homeomorphism between $W^\infty_{\alpha,\beta}$ and $\opFix(\Gamma)$.
\endproclaim
\proclabel{HomeoOfFixedPointSet}
\noindent Extending this to a homeomorphism between $\partial_\infty M$ and an abstract space
over which $\Gamma$ acts in a trivial manner, we obtain the following result which tells us that all the
Schottky groups that we will be studying are essentially equivalent to Schottky subgroups of
$\psl(2,\Bbb{C})$:
\proclaim{Lemma \nextprocno}
\noindent Suppose that the sectional curvature of $M$ is bounded above by $-k<0$.
Let $\Gamma\subseteq\opIsom_0(M)$ and $\Gamma'\subseteq\opIsom_0(M')$ be Schottky subgroups. There
exists an isomorphism $\phi:\Gamma\rightarrow\Gamma'$ and a homeomorphism 
$\Phi:\partial_\infty M\rightarrow\partial_\infty M'$ such that, for all $\gamma\in\Gamma$:
$$
\Phi\circ\gamma = \varphi(\gamma) \circ \Phi.
$$
\endproclaim
\proclabel{HomeoEquivalence}
\noindent In the second part of this appendix, we prove Lemma \procref{HomeoOfFixedPointSet}, and then
in the third part we show how this result may be used to prove Lemma \procref{HomeoEquivalence}.
\newsubhead{The Fixed Point Set of a Schottky Group}
\noindent In this subsection, we prove Lemma \procref{HomeoOfFixedPointSet}. We begin with the following
more elementary result concerning Hadamard manifolds of strictly negative sectional curvature:
\proclaim{Proposition \nextprocno}
\noindent Let $M$ be a Hadamard manifold of sectional curvature bounded above by $-k<0$. Let $UM$ be the 
unitary bundle over $M$ and let $\pi:UM\rightarrow M$ be the canonical projection. Let $(v_n)_\ninn$
be a sequence of vectors in $UM$ such that $(\pi\circ v_n)_\ninn$ converges to a point $p_0$ in
$\partial_\infty M$. For all $n$, let $\gamma_n$ be the geodesic in $M$ leaving $\pi\circ v_n$ with
velocity $v_n$. After extraction of a subsequence, at least one of $(\gamma_n(-\infty))_\ninn$ and
$(\gamma_n(+\infty))_\ninn$ converges also to $p_0$.
\endproclaim
\proclabel{LimitsOfEndPtsOfGeos}
\proof Let $q_0$ be a point in $M$. For all $n$, let us denote $p_n=\pi\circ v_n$, $D_n=d(q_0,p_n)$, and let 
$q_n\in\partial_\infty M$ be the unique point such that $p_n$ lies on the geodesic segment joining $q_0$
to $q_n$. Trivially $(q_n)_\ninn$ converges to $p_0$ and $(D_n)_\ninn$ tends to infinity. By Lemma \procref{LemmaConeControl}, without loss of generality, for all $n$:
$$
\gamma_n(+\infty)\in\partial_\infty\Omega_{\Delta_k}(D_n,\pi/2)(q_n;q_0).
$$
\noindent Since $(\Delta_k(D_n,\pi/2))_\ninn$ tends to zero, the result follows.\qed
\medskip
\noindent We now define the subset $X$ of $M$ by:
$$
X = M\setminus
\left(\opInt(D_\alpha^+)\munion\opInt(D_\alpha^-)\munion\opInt(D_\beta^+)\munion\opInt(D_\beta^-)\right).
$$
\noindent Trivially:
$$
\partial_\infty X = \partial_\infty M\setminus
\left(\opInt(C_\alpha^+)\munion\opInt(C_\alpha^-)\munion\opInt(C_\beta^+)\munion\opInt(C_\beta^-)\right).
$$
\noindent For all $\gamma\in\Gamma\setminus\left\{\opId\right\}$, let $\gamma^-$ and $\gamma^+$ be the
repulsive and attractive fixed points of $\gamma$ respectively, and let $g_\gamma$ be the unique geodesic
joining $\gamma^-$ to $\gamma^+$. We define $G_\Gamma$ by:
$$
G_\Gamma = \munion_{\gamma\in\Gamma}g_\gamma(\Bbb{R}).
$$
\noindent We obtain the following result:
\proclaim{Proposition \nextprocno}
\noindent The intersection $X\minter G_\Gamma$ is bounded in $M$.
\endproclaim
\proof Suppose the contrary. There exists a sequence $(p_n)_\ninn\in X\minter G_\Gamma$ and 
$p_0\in\partial_\infty M$ such that $(p_n)_\ninn$ converges to $p_0$. For all $n$, let 
$\gamma_n\in\Gamma$
be such that that $p_n$ lies in $g_{\gamma_n}$. For all $n$, let $\gamma_{n,\pm}$ be the fixed points
of $\gamma_n$. By Lemma \procref{LimitsOfEndPtsOfGeos}, we may assume that $(\gamma_{n,+})_\ninn$ tends
to $p_0$, in which case $p_0\in\opFix(\Gamma)$. However, $p_0\in\partial_\infty X$, which is absurd, since
$\opFix(\Gamma)$ and $\partial_\infty X$ are disjoint. The result now follows.\qed
\medskip
\noindent We define the evaluation map $\opEval:W_{\alpha,\beta}\rightarrow\Gamma$ such that, for all
$\pmb{\gamma}=(\gamma_k)_{1\leqslant k\leqslant n}$:
$$
\opEval(\pmb{\gamma}) = \gamma_n\cdot...\cdot\gamma_1.
$$
\noindent We have the following elementary result:
\proclaim{Lemma \nextprocno}
\noindent For every $\pmb{\gamma}=(\gamma_k)_{1\leqslant k\leqslant n}$ in $W_{\alpha,\beta}$, we have:
$$
\opExt(D_1(\gamma_n^{-1}))\cdot\opEval(\pmb{\gamma}) \subseteq \opInt(D_1(\gamma_1)).
$$
\endproclaim
\proclabel{AdvancedPingPong}
\proof We prove this result by induction on the length of $\pmb{\gamma}$. The result follows immediately
from the definition of generating circles when $\pmb{\gamma}$ is of length $1$. We suppose now that the
result is true when $\pmb{\gamma}$ is of length $n$. Let 
$\pmb{\gamma}=(\gamma_k)_{1\leqslant k\leqslant n}$ be a reduced word of length $n+1$ over 
$(\alpha,\beta)$. By the induction hypothesis, we have:
$$
\opExt(D_1(\gamma_{n+1}^{-1}))\cdot\opEval(\pmb{\gamma})\subseteq \opInt(D_1(\gamma_2))\cdot\gamma_1.
$$
\noindent Since $\pmb{\gamma}$ is a reduced word, $\gamma_2\neq\gamma_1^{-1}$. Consequently:
$$\matrix
&\opInt(D_1(\gamma_2)) \hfill&\subseteq \opExt(D_1(\gamma_1^{-1}))\hfill\cr
\Rightarrow\hfill&\opInt(D_1(\gamma_2))\cdot\gamma_1 \hfill& \subseteq \opInt(D_1(\gamma_1)).\hfill\cr
\endmatrix$$
\noindent The result now follows.\qed
\medskip
\noindent This yields the following two corollaries concerning $D(\pmb{\gamma})$ for $\pmb{\gamma}$ in
$W_{\alpha,\beta}$. First, we have:
\proclaim{Corollary \nextprocno}
\noindent For every $\pmb{\gamma}=(\gamma_k)_{1\leqslant k\leqslant n}$ in $W_{\alpha,\beta}$:
$$
\opInt(D(\pmb{\gamma})) = \opExt(D_1(\gamma_n^{-1}))\cdot\opEval(\pmb{\gamma}).
$$
\endproclaim
\proclabel{FormulaForInteriorOfCircles}
\proof This follows from Lemma \procref{AdvancedPingPong} and the fact that $D(\pmb{\gamma})$ is oriented
such that its interior is contained in the interior of the unique circle in $S_1$ in which it lies.\qed
\medskip
\noindent Next we have:
\proclaim{Corollary \nextprocno}
\noindent For every $\pmb{\gamma}=(\gamma_k)_{k\in\Bbb{N}}$ in $W^\infty_{\alpha,\beta}$ and for all 
$n\in\Bbb{N}$:
$$
\opInt(D(T_{n+1}(\pmb{\gamma}))) \subseteq \opInt(D(T_n(\pmb{\gamma}))).
$$
\endproclaim
\proclabel{CirclesAreNested}
\proof Indeed, by Corollary \procref{FormulaForInteriorOfCircles}:
$$
\opInt(D(T_{n+1}(\pmb{\gamma})))=\opInt(D_1(\gamma_{n+1}))\cdot\opEval(T_n(\pmb{\gamma})).
$$
\noindent Since $\pmb{\gamma}$ is a reduced word, $\gamma_{n+1}\neq\gamma_n^{-1}$, and so:
$$
\opInt(D_1(\gamma_{n+1}))\subseteq\opExt(D_1(\gamma_n)).
$$
\noindent The result follows by applying $\opEval(T_n(\pmb{\gamma}))$ to both sides and using 
Corollary \procref{CirclesAreNested}.\qed
\medskip
\noindent We define $W^0_{\alpha,\beta}$ to be the set of all reduced words 
$\pmb{\gamma}=(\gamma_k)_{1\leqslant k\leqslant n}$ in $W_{\alpha,\beta}$ such that 
$\gamma_n\neq\gamma_1^{-1}$. We have the following result:
\proclaim{Propositon \nextprocno}
\noindent $\|\opEval(\pmb{\gamma})\|$ tends to infinity as the length of $\pmb{\gamma}$ tends to 
infinity in $W^0_{\alpha,\beta}$.
\endproclaim
\proclabel{LengthOfEvalTendsToInfinity}
\proof Suppose the contrary. There exists $K>0$ and infinitely many distinct elements $(\gamma_n)_\ninn=(\opEval(\pmb{\gamma}_n))_\ninn$ in $\Gamma$ such that, for all $n$, $\|\gamma_n\| \leqslant K$. Since, for all $n$, $\gamma_n\in W^0_{\alpha,\beta}$, by Lemma \procref{AdvancedPingPong}, we may
suppose that there exist two distinct circles $C_1$ and $C_2$ in $S_1$ such that, for all $n$:
$$
\opExt(C_1)\cdot\gamma_n\subseteq\opInt(C_2).
$$
\noindent Consequently, if $\eta_n$ is the geodesic fixed by $\gamma_n$, then $\eta_n$ intersects
$X$ non-trivially. Thus, for all $n$, there exists $p_n\in X\minter G_\Gamma$ such that:
$$
d(p_n,p_n\cdot\gamma_n)\leqslant K.
$$
\noindent Let $p$ be an arbitrary point in $M$. By compactness, there exists $B>0$ such that,
for all $q\in X\minter G_\Gamma$:
$$
d(p,q) \leqslant B.
$$
\noindent Thus, for all $n$:
$$
d(p,p\cdot\gamma_n))\leqslant 2B + K.
$$
\noindent Consequently, by compactness, the set $\left\{p\right\}\cdot\Gamma$ has a concentration point
in $\tilde{M}$, which is absurd. The result now follows.\qed
\medskip
\noindent For $\pmb{\gamma}$ an element of $W^\infty_{\alpha,\beta}$, we define the sequence 
$(l_n(\pmb{\gamma}))_\ninn$ such that, for all $n$:
$$
l_n(\pmb{\gamma}) = \msup\left\{k\text{ s.t. }1\leqslant k\leqslant n\text{\ \&\ }\gamma_k\neq\gamma_1^{-1}\right\}.
$$
\noindent We now obtain the following partial proof of Lemma \procref{HomeoOfFixedPointSet}:
\proclaim{Proposition \nextprocno}
\noindent Suppose that the sectional curvature of $M$ is bounded above by $-k<0$.
\noindent For all $\pmb{\gamma}=(\gamma_n)_\ninn$ in $W^\infty_{\alpha,\beta}$ such that
$(l_n(\pmb{\gamma}))_\ninn$ tends to infinity, there exists a unique point 
$\Cal{P}(\pmb{\gamma})\in\opFix(\Gamma)$ such that the sequence $(\opInt(D(T_n(\pmb{\gamma}))))_\ninn$
converges towards $\left\{\Cal{P}(\pmb{\gamma})\right\}$ in the Hausdorff topology.
\endproclaim
\proclabel{CtyOfFPSPartI}
\proof For all $n$, we define $\pmb{\mu}_n\in W^{l_n(\pmb{\gamma})}_{\alpha,\beta}$ and 
$\pmb{\nu}_n\in W^{n-l_n(\pmb{\gamma})}_{\alpha,\beta}$ such that:
$$
\pmb{\mu}_n = T_{l_n(\pmb{\gamma})}(\pmb{\gamma}),\qquad T_n(\pmb{\gamma}) = \pmb{\nu}_n\pmb{\mu}_n.
$$
\noindent For all $n$, and for all $1\leqslant k\leqslant l_n(\pmb{\gamma})$, we denote by
$\mu_{n,k}$ the $k$'th letter of $\pmb{\mu}_n$. By corollary \procref{CirclesAreNested}, for all $n$:
$$
\opInt(D(T_n(\pmb{\gamma}))),\qquad\opInt(D(\pmb{\mu}_{n+1}))\subseteq\opInt(D(\pmb{\mu}_n)).
$$
\noindent We define $\Omega$ and $Delta$ as in section \subheadref{ConeControl}. By compactness, there exists $\epsilon>0$ such that, for all $p'\in X\minter G_\Gamma$, for all $q\in\opFix(\Gamma)$ and for all $C\in S_1$:
$$
C\minter\Omega_\epsilon(q;p')=\emptyset.
$$
\noindent Let $\theta>0$ be an angle. By compactness, there exists $\varphi>0$ such that, for all $p,p'\in X\minter G_\Gamma$:
$$
\partial_\infty\Omega_\varphi(q;p)\subseteq\partial_\infty\Omega_\theta(q;p').
$$
\noindent By Proposition \procref{LengthOfEvalTendsToInfinity}, there exists $N\in\Bbb{N}$ such that, for $n\geqslant N$:
$$
\Delta_k(\|\opEval(\pmb{\mu}_n)\|,\pi-\epsilon)<\varphi.
$$
\noindent Let us denote by $\eta_n$ the unique geodesique preserved by $\opEval(\pmb{\mu}_n)$. It
follows that $\eta_n$ intersects $X$ non-trivially. Let $p_n$ be any point in $X\minter\eta_n$. Let
$p_n^-$ and $p_n^+$ be the repulsive and attractive fixed points respectively of 
$\opEval(\pmb{\mu}_n)$. We have:
$$
\opExt(D_1(\mu_{n,l_n(\pmb{\gamma})}^{-1}))\subseteq
\partial_\infty\Omega_{\pi-\epsilon}(p_n^+;p_n).
$$
\noindent Thus, by Lemma \procref{LemmaConeControl} and Lemma \procref{AdvancedPingPong}:
$$\matrix
\opInt(D(\pmb{\mu}_n)) \hfill&=
\opExt(D_1(\mu_{n,l_n(\pmb{\gamma})}^{-1}))\cdot\opEval(\pmb{\mu})\hfill\cr
&\subseteq\partial_\infty\Omega_{\varphi}(p_n^+;p_n)\hfill\cr
&\subseteq\partial_\infty\Omega_\theta(p_n^+;p).\hfill\cr
\endmatrix$$
\noindent Thus, since $\theta$ is arbitrary, the diameter of $(\opInt(D(\pmb{\mu}_n)))_\ninn$ with respect to any given metric on $\partial_\infty M$ tends to zero.
\medskip
\noindent Since $(\overline{\opInt(D(\pmb{\mu}_n))})_\ninn$ is a nested sequence of compact sets, its
intersection is non-empty. Since its diameter is zero, the intersection contains at most $1$ point. 
Since, for all $n$, $\opInt(D(T_n(\pmb{\gamma}_n)))\subseteq\opInt(D(\pmb{\mu}_n))$, and the result follows.\qed
\medskip
\noindent We may then use this result to obtain:
\proclaim{Proposition \nextprocno}
\noindent Suppose that the sectional curvature of $M$ is bounded above by $-k<0$.
For all $\pmb{\gamma}$ in $W^\infty_{\alpha,\beta}$, there exists a unique point 
$\Cal{P}(\pmb{\gamma})\in\opFix(\Gamma)$ such that the sequence $(\opInt(D(T_n(\pmb{\gamma}))))_\ninn$
converges towards $\left\{\Cal{P}(\pmb{\gamma})\right\}$ in the Hausdorff topology.
\endproclaim
\proclabel{CtyOfFPSPartII}
\proof By Proposition \procref{CtyOfFPSPartI}, it suffices to prove the result when 
$(l_n(\pmb{\gamma}))_\ninn$ is bounded. Let $k$ be such that, for all $n$:
$$
l_n(\pmb{\gamma})\leqslant k.
$$
\noindent By definition of $(l_n(\pmb{\gamma}))_\ninn$, for all $m\geqslant k$:
$$
\gamma_m = \gamma_1^{-1}.
$$
\noindent We thus define $\pmb{\mu}=T_k(\pmb{\gamma})$, and we define 
$\pmb{\gamma}'\in W^\infty_{\alpha,\beta}$ such that:
$$
\pmb{\gamma} = \pmb{\gamma}'\pmb{\mu}.
$$
\noindent Trivially, for all $n$, $l_n(\pmb{\gamma}')=n$. By Corollary \procref{FormulaForInteriorOfCircles}:
$$
\opInt(D(T_n(\pmb{\gamma}')))\cdot\opEval(\pmb{\mu})=\opInt(D(T_{n+k}(\pmb{\gamma}))).
$$
\noindent We define $\Cal{P}(\pmb{\gamma})=\Cal{P}(\pmb{\gamma}')\cdot\opEval(\pmb{\mu})$, and the result
now follows by Proposition \procref{CtyOfFPSPartI}.\qed
\medskip
\noindent We now obtain a proof of Lemma \procref{HomeoOfFixedPointSet}:
\medskip
{\noindent\bf Proof of Lemma \procref{HomeoOfFixedPointSet}:\ }Existence and uniqueness of $\Cal{P}$ follow
from Propositions \procref{CtyOfFPSPartI} and \procref{CtyOfFPSPartII}. By corollary 
\procref{CirclesAreNested}, for all $\pmb{\gamma}$ and for all $n$:
$$
\Cal{P}(\pmb{\gamma}) \in \opInt(D(T_n(\pmb{\gamma}))).
$$
\noindent Continuity of $\Cal{P}$ now follows by a diagonal argument.
\medskip
\noindent Let $p$ be a point in $\opFix(\Gamma)$. For all $n$, $p$ lies in the union of the interiors
of the circles in $S_n$. Since $D_n$ defines a bijection between $W^n_{\alpha,\beta}$ and $S_n$, there exists $\pmb{\gamma}_n\in W^n_{\alpha,\beta}$ such that:
$$
p\in\opInt(D(\pmb{\gamma}_n)).
$$
\noindent By compactness, there exists $\pmb{\gamma}\in W^\infty_{\alpha,\beta}$ such that
$(\pmb{\gamma}_n)_\ninn$ tends to $\pmb{\gamma}$. For $m\in\Bbb{N}$ and for all $n$ sufficiently large,
$T_m(\pmb{\gamma})=T_m(\pmb{\gamma}_n)$. Thus, by Corollary \procref{CirclesAreNested}:
$$
p\in\opInt(D(\pmb{\gamma}_n))\subseteq\opInt(T_m(\pmb{\gamma})).
$$
\noindent By taking limits, it follows that $\Cal{P}(\pmb{\gamma})=p$, and surjectivity follows.
\medskip
\noindent Let $\pmb{\gamma}$ and $\pmb{\gamma}'$ be two points in $W^\infty_{\alpha,\beta}$ such
that $\Cal{P}(\pmb{\gamma})=\Cal{P}(\pmb{\gamma})=p$. For all $n$, the interiors of $D(T_n(\pmb{\gamma}))$ and $D(T_n(\pmb{\gamma}'))$ intersect non-trivially and thus coincide. Thus, for all $n$, 
$T_n(\pmb{\gamma})=T_n(\pmb{\gamma}')$. Injectivity now follows.
\medskip
\noindent Since $\Cal{P}$ is a bijective continuous mapping between two compact sets it is a
homeomorphism, and the result now follows.\qed
\newsubhead{Homeomorphism Equivalence of Schottky Groups}
\noindent Let $(C_\alpha^\pm,C_\beta^\pm)$ be generating circles of $\Gamma$ with respect to the generators 
$(\alpha,\beta)$. We define $\Omega\subseteq\partial_\infty M$ by:
$$
\Omega = \partial_\infty M\setminus
\left(\opInt(C_\alpha^+)\munion\opInt(C_\alpha^-)\munion\opInt(C_\beta^+)\munion\opInt(C_\beta^-)\right).
$$
\noindent We define the continuous mapping $\Phi:\Omega\times\Gamma\rightarrow\partial_\infty M$ by:
$$
\Phi(x,\gamma) = x\cdot\gamma.
$$
\noindent We define the equivalence relation $\sim$ over $\Omega\times\Gamma$ such that:
$$
(x,\gamma)\sim(y,\eta)\Leftrightarrow \Phi(x,\gamma) = \Phi(y,\eta).
$$
\noindent Since $\Omega$ is a fundamental domain for the action of $\Gamma$, we find that 
$(x,\gamma)\sim(y,\eta)$ if and only if either $x=y$ and $\gamma=\eta$ or:
\medskip
\myitem{(i)} $x,y\in\partial\Omega$, and
\medskip
\myitem{(ii)} there exists $\mu\in\left\{\alpha^{\pm 1},\beta^{\pm 1}\right\}$ such that:
$$
(x,\gamma) = (y\cdot\mu,\mu^{-1}\cdot\eta).
$$
\noindent $\Omega\times\Gamma/\sim$ is a Hausdorff space and that $\Phi$ quotients down onto
a homeomorphism of $\Omega\times\Gamma/\sim$ onto its image. Using Lemma \procref{HomeoOfFixedPointSet}, it is fairly trivial to show that $\partial_\infty M = \opIm(\Phi)\munion\opFix(\Gamma)$. We now define the space $\Sigma_\Gamma$ by:
$$
\Sigma_\Gamma = (\Omega\times\Gamma/\sim)\munion W^\infty_{\alpha,\beta}.
$$
\noindent For all $n$, and for all $\pmb{\gamma}\in W^n_{\alpha,\beta}$, we define the subset 
$U_{\pmb{\gamma}}$ of $\Sigma_\Gamma$ by:
$$
U_{\pmb{\gamma}} = 
\opInt\left(
\raise 2.5ex \hbox{$\displaystyle \munion_{
\matrix
{\scriptstyle \pmb{\eta}\in W_{\alpha,\beta}\text{ s.t }}\cr
{\scriptstyle l(\pmb{\eta})\geqslant n, T_n(\pmb{\eta})=\pmb{\gamma}}\cr
\endmatrix}$}
\Omega\times\left\{\opEval(\pmb{\eta})\right\}/\sim\right)
\munion\left\{\pmb{\gamma}'\in W^\infty_{\alpha,\beta}\text{ s.t. }T_n(\pmb{\gamma}')=\pmb{\gamma}\right\}.
$$
\noindent If $\tau$ denotes the topology of $\Omega\times\Gamma/\sim$, then 
$\tau\munion\left\{U_{\pmb{\gamma}}\text{ s.t. }\pmb{\gamma}\in W_{\alpha,\beta}\right\}$ defines a
base for a topology over $\Sigma_\Gamma$. This topology is Hausdorff and restricts to the initial topologies
over $\Omega\times\Gamma/\sim$ and $W^\infty_{\alpha,\beta}$. It is fairly trivial to show that 
$\Sigma_\Gamma$ is compact. Moreover, if $\pmb{\gamma}\in W^\infty_{\alpha,\beta}$
and if $(p_n,\gamma_n)_\ninn$ is a sequence in $\Omega\times\Gamma/\sim$ converging to $\pmb{\gamma}$,
then, for all $k\in\Bbb{N}$, there exists $N\in\Bbb{N}$ such that:
$$
n\geqslant N\Rightarrow T_k(\pmb{\gamma}_n) = T_k(\pmb{\gamma}).
$$
\noindent We now define the mapping $\Psi:\Sigma_\Gamma\rightarrow\partial_\infty\tilde{M}$ by:
$$
\Psi(p)=\left\{\matrix
\Cal{P}(p)\hfill&\text{ if }p\in W^\infty_{\alpha,\beta},\hfill\cr
\Phi(p)\hfill&\text{ if }p\in\Omega\times\Gamma/\sim.\hfill\cr
\endmatrix\right.
$$
\noindent Using Lemma \procref{HomeoOfFixedPointSet} we may show that the mapping $\Psi$ is continuous and bijective. Since $\Sigma_\Gamma$ is compact, $\Phi$ is a homeomorphism. We have thus constructed a topological space homeomorphic to the sphere. Moreover, there exists a canonical action of $\Gamma$ over this space with which $\Psi$ intertwines. We now prove Lemma \procref{HomeoEquivalence}:
\medskip
{\noindent\bf Proof of Lemma \procref{HomeoEquivalence}:\ }Let $(\alpha,\beta)$ and $(\alpha',\beta')$ be generators of $\Gamma$ and $\Gamma'$ respectively. Let
$(C_\alpha^\pm,C_\beta^\pm)$ and $(C'_\alpha{}^\pm, C'_\beta{}^\pm)$ be generating circles of $\Gamma$ and
$\Gamma'$ respectively with respect to these generators. We define $\Omega\subseteq\partial_\infty M$ and
$\Omega'\subseteq\partial_\infty M'$ as before. We define $\phi:\Gamma\rightarrow\Gamma'$ such that
$\phi(\alpha)=\alpha'$ and $\phi(\beta)=\beta'$. Let $\hat{\Phi}:\Omega\rightarrow\Omega'$ be a homeomorphism such that, for all $p\in C_\alpha^-$ and for all $q\in C_\beta^-$:
$$\matrix
\hat{\Phi}(p)\cdot\alpha' \hfill&= \hat{\Phi}(p\cdot\alpha),\hfill\cr
\hat{\Phi}(q)\cdot\beta' \hfill&= \hat{\Phi}(q\cdot\alpha).\hfill\cr
\endmatrix$$
\noindent There exists a unique extension of $\hat{\Phi}$ to a mapping from $\Sigma_\Gamma$ to 
$\Sigma_{\Gamma'}$ which intertwines with $\phi$. We define
$\Psi_1:\Sigma_\Gamma\rightarrow\partial_\infty M$ and $\Psi_2:\Sigma_\Gamma\rightarrow\partial_\infty M'$
as before. For all $\gamma\in\Gamma$, the following diagram trivially commutes:
$$
\commdiag{
\partial_\infty M &\mapright^{\Psi_1^{-1}}&\Sigma_\Gamma&\mapright^{\hat{\Phi}}&\Sigma_{\gamma'}&
\mapright^{\Psi'_1}&\partial_\infty M'\cr
\mapdown_{\gamma}& &\mapdown_{\gamma}& &\mapdown_{\phi(\gamma)}& &\mapdown_{\phi(\gamma)}\cr
\partial_\infty M &\mapright^{\Psi_1^{-1}}&\Sigma_\Gamma&\mapright^{\hat{\Phi}}&\Sigma_{\gamma'}&
\mapright^{\Psi'_1}&\partial_\infty M'\cr}
$$
\noindent The mapping $\Phi=\Psi_1'\circ\hat{\Phi}\circ\Psi_1^{-1}$ is thus the desired homeomorphism, and
the result now follows.\qed
\newhead{Bibliography}
{\leftskip = 5ex \parindent = -5ex
\leavevmode\hbox to 4ex{\hfil\cite{And}}\hskip 1ex{Anderson M., The Dirichlet Problem at infinity for manifolds of negative curvature, {\sl J. Diff. Geom.} {\bf 18} (1983), 701--702}
\medskip
\leavevmode\hbox to 4ex{\hfil\cite{BallGromSch}}\hskip 1ex{Ballman W., Gromov M., Schroeder V., {\sl Manifolds of nonpositive curvature}, Progress in Mathematics, {\bf 61}, Birkh\"auser, Boston, (1985)}
\medskip
\leavevmode\hbox to 4ex{\hfil\cite{Frib}}\hskip 1ex{Friberg B., A topological proof of a theorem of Kneser, {\sl Proc. Amer. Math. Soc.} {\bf 39} (1973), 421--425}
\medskip
\leavevmode\hbox to 4ex{\hfil\cite{GallKapMard}}\hskip 1ex{Gallo D., Kapovich M., Marden A., The monodromy groups of Schwarzian equations on closed Riemann surfaces, {\sl Ann. Math.} {\bf 151} (2000), no. 2, 625--704}
\medskip
\leavevmode\hbox to 4ex{\hfil\cite{LabA}}\hskip 1ex{Labourie F., Un lemme de Morse pour les surfaces convexes, {\sl Invent. Math.} {\bf 141} (2000), 239--297}
\medskip
\leavevmode\hbox to 4ex{\hfil\cite{LeRoux}}\hskip 1ex{Le Roux F., Th\`ese doctorale, Grenoble (1997)}
\medskip
\leavevmode\hbox to 4ex{\hfil\cite{SmiB}}\hskip 1ex{Smith G., Hyperbolic Plateau problems, Preprint, Orsay (2005)}%
\medskip
\leavevmode\hbox to 4ex{\hfil\cite{SmiE}}\hskip 1ex{Smith G., Th\`ese de doctorat, Paris (2004)}%
\medskip
\leavevmode\hbox to 4ex{\hfil\cite{SmiF}}\hskip 1ex{Smith G., A Homomorphism for Braid Groups in the Sphere, in preparation}
\par}%
\enddocument